\documentclass[a4paper,reqno]{amsart}

\usepackage[all]{xy}
\usepackage{xspace}
\usepackage{amsmath}
\usepackage{amstext}
\usepackage{amsfonts}
\usepackage[mathscr]{euscript}
\usepackage{amscd}
\usepackage{latexsym}
\usepackage{amssymb}
\usepackage{diagrams}

\newtheorem{theorem}{Theorem}[section]
\newtheorem{lemma}[theorem]{Lemma}
\newtheorem{proposition}[theorem]{Proposition}
\newtheorem{definition}[theorem]{Definition}
\newtheorem{corollary}[theorem]{Corollary}
\newtheorem{remark}[theorem]{Remark}
\theoremstyle{definition}

\newcommand{\dbb}{\mathbb{D}}
\newcommand{\abb}{\mathbb{A}}
\newcommand{\bbb}{\mathbb{B}}
\newcommand{\bsy}[1]{\boldsymbol{#1}}
\newcommand{\mi}{\text{-}}
\newcommand{\enpr}{\hfill $\Box $}
\newcommand{\supar}[1]{\overset{#1}{-\!\!-\!\!\!\rightarrow}}
\newcommand{\suparle}[1]{\overset{#1}{\leftarrow\!\!\!-\!\!-}}

\newcommand{\prf}{\noindent\textbf{Proof. }}
\newcommand{\rw}{\rightarrow}
\newcommand{\lw}{\leftarrow}

\newcommand{\vep}{\varepsilon}

\newcommand{\id}{\mathrm{id}}
\newcommand{\hm}{\mathrm{Hom}}

\newcommand{\im}{\mathrm{Im}\,}

\newcommand{\hund}[1]{{\hm}_{#1}}

\newcommand{\tiund}[1]{{\times}_{#1}}

\newcommand{\ovl}[1]{\overline{#1}}

\newcommand{\til}[1]{\widetilde{#1}}
\newcommand{\pr}{\mathrm{pr}}

\newcommand{\dcl}{\mathcal{D}}

\newcommand{\ccl}{\mathcal{C}}

\newcommand{\gcl}{\mathcal{G}}
\newcommand{\lcl}{\mathcal{L}}

\newcommand{\ncl}{\mathcal{N}}
\newcommand{\fcl}{\mathcal{F}}
\newcommand{\tcl}{\mathcal{T}}
\newcommand{\ucl}{\mathcal{U}}
\newcommand{\vcl}{\mathcal{V}}
\newcommand{\wcl}{\mathcal{W}}
\newcommand{\hcl}{\mathcal{H}}
\newcommand{\kcl}{\mathcal{K}}
\newcommand{\pcl}{\mathcal{P}}
\newcommand{\rcl}{\mathcal{R}}

\newcommand{\cig}{C_i(\gcl)}
\newcommand{\cg}{C(\gcl)}
\newcommand{\cng}{\ccl^n\gcl}
\newcommand{\cngm}{\ccl^{n-1}\gcl}
\newcommand{\cngp}{\ccl^{n+1}\gcl}

\newcommand{\set}{\mathbf{Set}}

\newcommand{\ttil}{\widetilde{\times}}
\newcommand{\ob}{Ob\,}
\newcommand{\ner}{Ner\,}

\newcommand{\cat}{\mathrm{Cat}\,}

\newcommand{\ds}{ds\,}

\newcommand{\gp}{\mathrm{Gp}}
\newcommand{\gpd}{\mathrm{Gpd}}
\newcommand{\pt}{\partial}

\newcommand{\tam}{Tamsamani's }

\newcommand{\pp}{\bsy{\Pi}}

\newcommand{\cgp}{\mathrm{Cat}(\gp)}
\newcommand{\cagp}{\mathrm{Cat}^n(\gp)}
\newcommand{\cagpm}{\mathrm{Cat}^{n-1}(\gp)}
\newarrow{Dts}{}{.}{.}{.}{>}
\newcommand{\ctn}{\mathrm{cat}^{n}\mi\mathrm{group}}
\newcommand{\ctmen}{\mathrm{cat}^{n-1}\mi\mathrm{group}}

\newcommand{\iun}{i=1,\ldots,n}

\newcommand{\dop}{\Delta^{op}}
\newcommand{\dnop}{\Delta^{n^{op}}}
\newcommand{\dnopiu}{\Delta^{{n+1}^{op}}}
\newcommand{\dnomen}{\Delta^{{n-1}^{op}}}
\newcommand{\duop}{\Delta^{{2}^{op}}}

\newcommand{\sci}{SC\,}
\newcommand{\np}{{n+1}}
\newcommand{\nm}{{n-1}}
\newcommand{\ovu}{\ovl{U}}
\newcommand{\diag}{\,diag\,}
\newcommand{\spi}{Sp\,}
\newcommand{\bsim}{/\!\!\sim}
\newcommand{\vi}{,\;}
\newcommand{\tp}{\mathrm{Top}}
\newcommand{\tiu}{T^{(1)}}
\newcommand{\tid}{T^{(2)}}
\newcommand{\ta}[2]{\tau^{(#1)}_{#2}}
\newcommand{\fistr}[1]{\phi_1\tiund{\phi_0}\overset{#1}{\cdots}\tiund{\phi_0}\phi_1}
\newcommand{\fistrr}[2]{#2\phi_1\tiund{#2\phi_0}\overset{#1}{\cdots}\tiund{#2\phi_0}#2\phi_1}
\newcommand{\fistpar}[2]{(#2\phi)_1\tiund{(#2\phi)_0}\overset{#1}{\cdots}\tiund{(#2\phi)_0}(#2\phi)_1}
\newcommand{\gistr}[1]{\gcl_1\tiund{\gcl_0}\overset{#1}{\cdots}\tiund{\gcl_0}\gcl_1}

\newcommand{\gistpar}[2]{(#2\gcl)_1\tiund{(#2\gcl)_0}\overset{#1}{\cdots}\tiund{(#2\gcl)_0}(#2\gcl)_1}
\newcommand{\pro}[3]{#1\tiund{#2}\overset{#3}{\cdots}\tiund{#2}#1}
\newcommand{\fl}{\text{\begingroup\Large$\mathfrak{f}$\normalsize\endgroup}}

\newcommand{\ul}{\text{\begingroup\Large$\mathfrak{u}$\normalsize\endgroup}}
\newcommand{\tens}[2]{#1\,\tiund{#2}\,#1}

\begin{document}

\title {Semistrict Tamsamani $\mathbf{N}$-groupoids and connected $\mathbf{N}$-types}

\author{Simona Paoli}

\date{2 October 2007}

\address{Department of Mathematics\\ Macquarie University\\ NSW 2109 \\ Australia}
\email{simonap@ics.mq.edu.au}



\begin{abstract}
Tamsamani's weak $n$-groupoids are known to model $n$-types. In this paper we
show that every Tamsamani weak $n$-groupoid representing a connected $n$-type
is equivalent in a suitable way to a semistrict one. We obtain this result by
comparing Tamsamani's weak $n$-groupoids and cat$^{n-1}$-groups as models of
connected $n$-types.
\end{abstract}

\maketitle



\section{Introduction}\label{int}
A very important and non-trivial problem in higher category theory is that of
finding coherence theorems for weak higher categorical structures. Broadly
speaking, one way to formulate a coherence theorem for a weak higher
categorical structure consists of saying that it is, in a suitable sense,
equivalent to one in which some of the associativity and identity laws hold
strictly. These structures which are ``less weak" than  general weak
$n$-categories are often called ``semistrict" \cite{be}.

The first example of a coherence theorem is the one for monoidal categories due
to Mac Lane \cite{macl}. In turn, this is a special case of a more general
coherence theorem for bicategories  \cite{mp} which says that every bicategory
is biequivalent to a strict 2-category. In dimension 3, Gordon, Power and
Street \cite{gps} showed that every tricategory is triequivalent to a Gray
category. The latter is a type of semistrict tricategory. These low-dimensional
coherence results also have analogues in the
 groupoid case.

The fundamental information carried by a weak $n$-groupoid is its homotopy
type. Indeed, it is thought that any good definition of weak higher category
should give a model of $n$-types in the weak $n$-groupoid case. For Tamsamani's
model of weak $n$-categories this property was proved in \cite{tam}
\cite{tam1}. More precisely, Tamsamani proved that there is an equivalence of
categories between the homotopy category of $n$-types and the localization of
weak $n$-groupoids with respect to some suitably defined $n$-equivalences.

In low dimensions, it is known that strict $2$-groupoids model $2$-types
\cite{ms} and Gray groupoids model $3$-types \cite{be} \cite{jt} \cite{le}.
These low dimensional results led several people (see for instance \cite{bd})
to formulate the \emph{semistrictification hypothesis} for homotopy types: In
every model of weak $n$-category, a weak $n$-groupoid should be suitably
equivalent to a semistrict one, and, further weak $n$-groupoids should model
$n$-types. A more general form of the semistrictification hypothesis states
that every weak $n$-category should be suitably equivalent to a semistrict one.

Until now, no proof of this hypothesis exists in dimension higher than $3$. In
this paper we establish a semistrictification result for Tamsamani weak
$n$-groupoids, valid for any $n$ in the path-connected case. Our main result,
Theorem \ref{result.the1}, states that every Tamsamani weak $n$-groupoid such
that its classifying space is path-connected can be linked by a zig-zag of
$n$-equivalences to a weak $n$-groupoid which is semistrict.

The category $\hcl_n$ of semistrict Tamsamani $n$-groupoids is defined in
Section \ref{semtam} and is the full subcategory of the category $\tcl_n$ of
Tamsamani weak $n$-groupoids consisting of
$\phi\in\tcl_n\subset[\dop,\tcl_\nm]$ such that $\phi_0$ is the trivial
$(n-1)$-groupoid and the Segal maps $\phi_k\rw\fistr{k}$ are isomorphisms (see
Section \ref{simtec} for the definition of Segal maps). As explained in Remark
\ref{tam3.rem1} c), $\hcl_n$ is isomorphic to a full subcategory of the
category of monoids in the category of Tamsamani weak $(\nm)$-groupoids with
respect to the cartesian monoidal structure. From the definition of $\hcl_n$,
the Segal maps in the top simplicial direction are isomorphisms. Since in
Tamsamani's model the Segal maps encode the composition law, this means that
the composition of cells along that direction is strict. In the other
directions, however, composition is in general weak. This is why objects of
$\hcl_n$ are called ``semistrict".

Semistrictification is a theme which has given rise to a great deal of interest
recently. While classically Gray groupoids and Gray categories were the primary
examples of semistrict structures, more recently Carlos Simpson formulated some
conjectures \cite{simp} about semistrict structures where everything is strict
except units which are weak. A precise formulation of this conjecture has been
given by Joachim Kock \cite{kock} in terms of so called ``fair $n$-categories".
In particular, Joachim Kock conjectured that fair $n$-groupoids should model
$n$-types. A proof of this conjecture was given by Joyal and Kock \cite{joko}
in the special case of 1-connected 3-types.

As was pointed out by Tom Leinster, another definition of higher category in
which semistrict structures naturally arise is the one due to Trimble
\cite{lein} \cite{che1} \cite{cheg}. Recently Eugenia Cheng \cite{che2}
established a comparison result between Trimble's and Batanin's definitions
\cite{bata} leading to an interesting class of semistrict $n$-categories.

It is natural to ask what  the relationship is between our result and these
recent works on semistrictification. We present some conjectures on this issue
in Section \ref{furth} a). As discussed there, we think that the type of
semistrictification studied in this paper ought to correspond to the ``weak
units" approach rather than the Gray groupoids (and their possible
higher-dimensional generalizations). Making this precise is however a
substantial project which goes outside the scope of the present work, and we
will tackle it in a subsequent paper.

Our semistrictification result is very significant homotopically. In fact we
know by \cite{tam} that the homotopy category of connected $n$-types $\hcl
o(\tp^{(n)}_*)$ is equivalent to $\tcl_n^*\bsim^n$ (where objects of $\tcl_n^*$
are Tamsamani weak $n$-groupoids whose classifying space is path-connected);
Theorem \ref{result.the1} implies immediately (Corollary \ref{result.cor1})
that $\hcl o(\tp^{(n)}_*)$ is in fact equivalent to a category which is smaller
than $\tcl_n^*\bsim^n$, more precisely it is the full subcategory whose objects
are in $\hcl_n$. A further discussion and open question about this point can be
found in Section \ref{furth} c).

From a higher categorical perspective, it would be desirable that the zig-zag
of $n$-equivalences of Theorem \ref{result.the1} could be made into a single
map. This is at present an open question: as discussed further in Section
\ref{furth} b) the problem of inverting $n$-equivalences in the Tamsamani model
is very delicate and may require some non trivial modifications of the
morphisms in this category. Although very interesting, this investigation is
outside the scope of the present work.

On the other hand, we believe the result of Theorem \ref{result.the1} is very
significant from a higher categorical perspective because no version of
semistrictification, even in a weaker or special case, appears to have been
proved in the literature in dimension higher that 3. Another reason why we
think the present work is of interest to higher category theory is because of
the method we use to prove our result, which relies on a comparison between two
very different types of higher categorical structures: the Tamsamani model and
the $\ctn$ model. The latter is a model of connected $(n+1)$-types which was
introduced by Loday \cite{lod} and later developed in \cite{bcd}, \cite{por},
generalizing earlier work of Whitehead \cite{whit} on crossed modules. Our
semistrictification functor is given by the composite
\begin{equation*}
    S:\tcl_\np\supar{B}\tp\supar{\pcl_n}\cagp\supar{F}\hcl_\np.
\end{equation*}
Here $B$ is the classifying space functor, $\pcl_n$ is the fundamental $\ctn$
functor of a topological space from \cite{bcd} and $F$ is the functor whose
construction occupies the main part of this paper. An overview of the method
used to construct the functor $F$ can be found in Section \ref{orga}. The
functor $F$ preserves the homotopy type; that is, for each $\ctn$ $\gcl$ there
is an isomorphism $B\gcl\cong BF\gcl$ in $\hcl o(\tp^{(\np)}_*)$. The existence
of $F$ and the fact that it preserves the homotopy type (Corollary
\ref{semist.cor01}) easily imply the semistrictification result, Theorem
\ref{result.the1}. The functor $F$ establishes a direct comparison between two
higher categorical structures which have very different features: in fact
$\ctn$s, being $n$-fold categories internal to groups, are strict cubical
structures while Tamsamani's model is a weak globular model. In category theory
$n$-fold structures have been studied extensively in the work of Ehresmann and
more recently, for $n=2$, (that is, for double categories), by several authors
such as Dawson, Grandis, Par\'{e} and Pronk. These works, however, did not
establish any comparison with the weak globular structures of higher category
theory. Our work establishes, in the special case of $n$-fold structures
internal to groups, a first bridge between these two areas. Even if the
extension to the non-groupoidal setting is challenging (see Section \ref{furth}
d) for a further discussion on this point) we believe this work is an important
first step in understanding how $n$-fold structures interact with higher
category theory.

\bigskip

\textbf{Acknowledgements} This work was supported by an Australian Research
Council Postdoctoral Fellowship (Project No. DP0558598). I wish to thank Ross
Street for helpful comments and for his constant encouragement. I also thank
Clemens Berger, Joachim Kock, Dorette Pronk, Carlos Simpson for useful feedback
on this work. I finally extend my gratitude to the referee for useful
suggestions on improving the presentation of the paper.


\section{Organization of the paper and overview of the main results.}\label{orga}
The main result of the paper (Theorem \ref{result.the1}) asserts the existence
of a ``semi-strictification" functor $S$ from the category $\tcl_\np$ of
Tamsamani weak $(\np)$-groupoids to the category $\hcl_\np$ of semistrict
Tamsamani $(\np)$-groupoids, with the property that, if $\gcl\in\tcl_\np$ is
such that $B\gcl$ is path-connected, there exists a zig-zag of
$(\np)$-equivalences in $\tcl_\np$ connecting $\gcl$ and $S(\gcl)$. The
relevance of this result in relation to the modelling of homotopy types is
discussed in Corollary \ref{result.cor1}. The semistrictification functor $S$
is given by the composite
\begin{equation*}
    S:\tcl_\np\supar{B}\tp\supar{\pcl_n}\cagp\supar{F}\hcl_\np.
\end{equation*}
Here $B$ is the classifying space functor and $\pcl_n$ is the fundamental
$\ctn$ functor; the latter is well known in the literature \cite{brr} for
inducing an equivalence of categories $\cagp/\!\!\!\sim\;\cong\hcl
o(\tp^{(\np)}_*)$ between the localization of $\ctn$s with respect to weak
equivalences and the homotopy category of connected $(\np)$-types. The functor
$F$ preserves the homotopy type. That is, for each $\ctn$ $\gcl$ there is an
isomorphism $B\gcl\cong BF\gcl$ in $\hcl o(\tp^{(\np)}_*)$. The proof of its
existence and of the fact that it preserves the homotopy type (Corollary
\ref{semist.cor01}) occupy the main part of the paper. Together with the
properties of $\pcl_n$ (Theorem \ref{catn.the1}) and of Tamsamani's model
(Proposition \ref{group.pro1}), Corollary \ref{semist.cor01} immediately
implies the semistrictification result Theorem \ref{result.the1}.

The functor $F$ is defined via three intermediate steps, namely it is given by
the composite
\begin{equation*}
    F:\cagp\supar{\spi}\cagp_S\supar{\dcl_n}\dbb_n\supar{V_n}\hcl_\np.
\end{equation*}
The category $\cagp_S$ is a full subcategory of $\cagp$ whose objects are
called \emph{special cat$^n$-groups} (Definition \ref{defspec.def1}). The name
of the functor $\spi$ thus stands for ``specialization". The definition of
special $\ctn$s needs the notion of \emph{strongly contractible cat$^n$-groups}
(Definition \ref{contra.def1}). Roughly speaking a strongly contractible $\ctn$
$\gcl$ is one for which both $\gcl$ and all its subfaces are ``homotopically
discrete" in a very strong form; a special $\ctn$ is one for which certain
prescribed faces are strongly contractible. More precisely, these faces
correspond to those that in a globular object (that is, a strict $n$-category
internal to groups) are discrete. A more lengthy informal discussion of the
ideas behind these notions can be found at the beginnings of Sections
\ref{contra} and \ref{defspec}.

The category $\dbb_n$ is called the category of \emph{internal weak
$n$-groupoids} (Definition \ref{intweak.def1}). This can be thought of as an
intermediate model between $\ctn$s and the Tamsamani model. A more lengthy
informal discussion about $\dbb_n$ can be found at the beginning of Section
\ref{intweak}. The existence of the functors $\spi$, $\dcl_n$, $V_n$ and the
fact that each of them preserves the homotopy type is established respectively
in Theorems \ref{specgen.the1}, \ref{globgen.the1} and \ref{semist.the01}.

We are now going to describe the main steps leading up to each of these
theorems, starting with Theorem \ref{specgen.the1}. We need two main results
about the category of $\ctn$s, which are Theorem \ref{adj.the1} and Proposition
\ref{twolem.pro2}. Section \ref{catcat} is devoted to their proof, which
necessitates an equivalent definition of $\ctn$s as a category of groups with
operations. We refer the reader to the beginning of Section \ref{catcat} for an
overview of the method used to prove these results. Theorem \ref{adj.the1}
asserts the existence of an adjunction $\fl_n\vdash\ul_n:\cagp \rw\set$ and,
for any set $X$, gives an explicit description of $\fl_n(X)$ in terms of
$\fl_\nm(X)$. Proposition \ref{twolem.pro2} establishes that if a morphism $f$
in $\cagp$ is such that $\ul_n f$ is surjective then the corresponding map of
multinerves $\ncl f$ is levelwise surjective.

These two results about $\ctn$s are then used as follows. Theorem
\ref{adj.the1} is the main ingredient used in proving that, for any set $X$,
$\fl_n(X)$ is strongly contractible (see Lemma \ref{catn.lem3} and Theorem
\ref{spec.the0}): The proof of this fact uses precisely the explicit
description of $\fl_n(X)$ in terms of $\fl_\nm(X)$ provided by Theorem
\ref{adj.the1}. The strong contractibility of $\fl_n(X)$ together with
Proposition \ref{twolem.pro2}  are two of the three main ingredients in the
proof of Proposition \ref{spec.pro4}: the latter asserts the existence, for
each $1\leq i \leq n$, of a functor $C_i:\cagp\rw\cagp$ such that, for any
$\gcl\in\cagp$, as an internal category in $\cagpm$ in the $i^{th}$ direction,
$C_i(\gcl)$ has object of objects a strongly contractible $\ctmen$; further, it
shows there is a weak equivalence $\alpha^{(i)}_\gcl:C_i(\gcl)\rw\gcl$, natural
in $\gcl$. The third main ingredient in the proof of Proposition
\ref{spec.pro4} is a standard construction on internal categories described in
Section \ref{cons}. We call the functor $C_i$ the \emph{cofibrant replacement
on cat$^n$-groups in the $i^{th}$ direction} (Definition \ref{apcat.def1}). The
reason for this terminology is explained in Remark \ref{apcat.rem1}.

The functors $C_i$ are used to construct $\spi$. For $n=2$ it is $\spi=C_1$,
while for general $n$ we need to take $\nm$ cofibrant replacements; that is,
$\spi=C_1C_2\ldots C_\nm$. The reason for this is carefully explained in
Sections \ref{spelow} and \ref{specgen}. In Section \ref{spelow} we illustrate
some low dimensional cases in order to gain intuition for general $n$; in
Section \ref{specgen} the formal treatment of the general case (Lemma
\ref{specgen.lem1} and Theorem \ref{specgen.the1}) follows a more informal
discussion of the main ideas used in the proof. In addition to the strong
contractibility of $\fl_n(X)$ and of Proposition \ref{twolem.pro2}, the proof
of Theorem \ref{specgen.the1} requires the fact that strong contractibility
behaves well with respect to certain pullbacks: this is proved in Lemma
\ref{spec.lem1}, which also plays an important role in the construction of
$\dcl_n$. The main steps in the construction of $\spi$ are:
\begin{itemize}
  \item[a)]The adjunction $\fl_n\dashv \ul_n$ (Theorem \ref{adj.the1}) implies,
  (with Lemma \ref{spec.lem1} and Lemma \ref{catn.lem3}), that
  \item[b)]$\fl_n(X)$ is strongly contractible (Theorem \ref{spec.the0}). This
  implies, (with Proposition \ref{twolem.pro2} and Lemma \ref{spec.lem3})
  \item[c)]functors $C_i:\cagp\rw\cagp$ (Proposition \ref{spec.pro4}). This
  implies (with Lemma \ref{spec.lem1} and Proposition \ref{twolem.pro2})
  \item[d)]Lemma \ref{specgen.lem1} and hence functor $\spi$ (Theorem \ref{specgen.the1}
  ).
\end{itemize}
We now come to the functor $\dcl_n$. The idea of this functor is to ``squeeze"
the strongly contractible faces of a special $\ctn$ to become discrete ones in
such a way that the homotopy type is preserved. This will produce a globular
object from a cubical one but will destroy the strictness of the Segal maps,
which will become only weak equivalences. The definition of $\dcl_n$ is
obtained by iteration of a basic construction, called the \emph{discrete nerve}
$\ds\ncl$, which we describe in Section \ref{disner} in a general context. The
properties of $\dcl_2$ are described in Theorem \ref{glob2.the1}, which is the
first step in the inductive argument for general $n$ given in Theorem
\ref{globgen.the1}. A key property that allows us to define $\dcl_n$ for
general $n$ is the fact (Lemma \ref{globgen.lem0}) that the discrete nerve
construction of Section \ref{disner}, when applied to the case of special
$\ctn$s, yields a functor
\begin{equation}\label{orga.eq1}
    \ds\ncl:\cagp_S\rw[\dop,\cagpm_S].
\end{equation}
This is an immediate consequence of the fact, proved in Corollary
\ref{nerspec.cor1}, that the usual nerve functor $\ner:\cagp\rw[\dop,\cagpm]$
restricts to a functor $\ner:\cagp_S\rw[\dop,\cagpm_S]$. In turn, this relies
on a property of special $\ctn$s with respect to certain pullbacks (Lemma
\ref{nerspec.lem1}) which is a direct consequence of the already mentioned
Lemma \ref{spec.lem1}. Thanks to  (\ref{orga.eq1}) we can define
$\dcl_2=\ds\ncl$ and for $n>2$ $\dcl_n:\cagp_S\rw[\dop,\dbb_\nm]\,$
$\dcl_n=\ovl{\dcl_\nm}\circ\ds\ncl$, where $\ovl{\dcl_\nm}$ is the functor
induced by applying $\dcl_\nm$ levelwise (in the sense of Definition
\ref{simtec.def1}).

The main point of Theorem \ref{globgen.the1} is to show that the image of
$\dcl_n$ is in fact in $\dbb_n\subset[\dop,\dbb_\nm]$ and that $\dcl_n$
preserves the homotopy type. An overview of the main ideas used in this proof
can be found at the beginning of Section \ref{globgen}, before the formal
treatment. The main steps leading to Theorem \ref{globgen.the1} are:
\begin{itemize}
  \item[a)] Lemma \ref{spec.lem1} implies (with Lemma \ref{nerspec.lem1}) that
  \item[b)] the functor $\ner:\cagp\rw[\dop,\cagpm]$ restricts to
  $\ner:\cagp_S\rw[\dop,\cagpm_S]$ (Corollary \ref{nerspec.cor1}); this implies
  that
  \item[c)] there is a functor $\ds\ncl:\cagp_S\rw[\dop,\cagpm_S]$ (Lemma \ref{globgen.lem0}).
  Hence define $\dcl_2=\ds\ncl$ and
  $\dcl_n=\ovl{\dcl_\nm}\circ\ds\ncl$ for $n>2$.
  \item[d)] Use the above definitions together with Lemmas \ref{weak.lem01}
  and \ref{globgen.lem1} to show Theorem \ref{globgen.the1}.
\end{itemize}
We finally come to the functor $V_n$. We first define a functor
\begin{equation*}
    V_n:\underset{\nm}{\underbrace{[\dop,[\dop,\ldots,[\dop}},\cgp]\cdots]
    \rw\underset{n}{\underbrace{[\dop,[\dop,\ldots,[\dop}},\gpd]\cdots].
\end{equation*}
Informally, the definition of $V_n$ simply amounts to taking the nerve of each
group occurring in its domain. It is necessary for its use that this definition
of $V_n$ is made precise and the notation set up carefully: we have done so in
Section \ref{deloop}. As explained in Remarks \ref{teclem.rem1} and
\ref{intweak.rem1}, $\dbb_n$ and $\tcl_\np$ are full subcategories respectively
of the domain and codomain of $V_n$. What Theorem \ref{semist.the01} asserts is
that $V_n$ restricts to a functor $V_n:\dbb_n\rw\hcl_\np$ and, further, that it
preserves the homotopy type and sends weak equivalences to
$(n+1)$-equivalences.

The case $n=2$ (which is the first step of the inductive argument for general
$n$) is treated in Section \ref{tam3}, which starts with an overview of the
main ideas used in the proof. The general case is treated in Section
\ref{semist}. Because of the centrality of Theorem \ref{semist.the01}, we have
given in Section \ref{result} a thorough overview of the method used in the
proof of the general case, so we refer the reader to that section. Here we just
point out that the only facts from Sections \ref{pre} to \ref{weak} used in the
proof of Theorem \ref{semist.the01} are the definition of the category
$\dbb_n$, the definition and properties of Tamsamani model and in particular
Lemma \ref{teclem.lem1} (which is a simple calculation from the definition) and
the elementary fact, (recalled in Lemma \ref{catn.lem1}), that the functor
$\pi_0:\cgp\rw\gp$ preserves fibre products over discrete objects and sends
weak equivalences to isomorphisms. An absolutely essential ingredient in the
proof of Theorem \ref{semist.the01} is the fact, shown in Lemma
\ref{geom.lem1}, that a map $f$ of Tamsamani weak $n$-groupoids is a
$n$-equivalence if and only if the corresponding map $Bf$ of classifying spaces
is a weak equivalence.

We have provided a concise expository account of the main simplicial techniques
used in the paper in Section \ref{simtec}; an account of the Tamsamani model
following \cite{toe} is in Section \ref{pre}. We have concluded the paper with
a discussion of some conjectures and themes for further investigations (Section
\ref{furth}).

\section{Simplicial techniques.}\label{simtec}
In this section we recall some basic background on simplicial techniques. The
material in this section is well known. Further details can be found for
instance in \cite{goja}, \cite{may}, \cite{weib}, \cite{hove}, \cite{sega},
\cite{hirs}, \cite{qui1}. Throughout this paper, we will use the following
notation.
\begin{definition}\label{simtec.def1}
    Let $F:\ccl\rw\dcl$ be a functor, $\kcl$ a small category and
    $[\kcl,\ccl]$, $[\kcl,\dcl]$ the functor categories. We denote by
    $\ovl{F}=[\kcl,F]:[\kcl,\ccl]\rw[\kcl,\dcl]$ the functor defined as follows. Given
    a functor $\phi:\kcl\rw\ccl$, $\ovl{F}\phi:\kcl\rw\dcl$ is
    defined by $(\ovl{F}\phi)(x)=F(\phi(x))$ for all $x\in\ob \kcl$. Given a
    natural transformation $f:\phi\rw\psi$
    with components $f_x:\phi(x)\rw\psi(x)$, $x\in\ob\kcl$, define
    $\ovl{F}f:\ovl{F}\phi\rw\ovl{F}\psi$ to be the natural transformation with
    components $(\ovl{F}f)_x=F(f_x)$ for all $x\in\ob\kcl$. We will call
    $\ovl{F}$ the functor ``induced by applying $F$ levelwise", or simply
    ``induced by $F$".
\end{definition}
The structures used in this paper can be regarded as multi-simplicial
structures, thus we recall some elementary facts about multi-simplicial
objects. Let $\Delta$ be the simplicial category and let $\dnop$ denote the
product of $n$ copies of $\dop$. Given a category $\ccl$, $[\dnop,\ccl]$ is
called the category of \emph{$n$-simplicial objects} in $\ccl$
(\emph{simplicial objects} in $\ccl$ when $n=1$). If $\phi\in[\dnop,\ccl]$ and
$([p_1]\ldots[p_n])\in\dnop$, we shall often denote $\phi([p_1]\ldots[p_n])$ by
$\phi(p_1\ldots p_n)$ or by $\phi_{p_1\ldots p_n}$.

Every $n$-simplicial object in $\ccl$ can be regarded as a simplicial object in
$[\dnomen,\ccl]$ in $n$ possible ways. That is, for each $1\leq k\leq n$, there
is an isomorphism $\xi_k:[\dnop,\ccl]\rw[\dop,[\dnomen,\ccl]]$ as follows:
given $\phi\in[\dnop,\ccl]$, $\xi_k\phi$ is the simplicial object taking
$[r]\in\dop$ to $(\xi_k\phi)_r\in[\dnomen,\ccl]$ defined by
$(\xi_k\phi)_r(p_1\ldots p_\nm)=\phi(p_1\ldots p_{k-1}\,r\,p_k\ldots p_\nm)$.
The inverse $\xi'_k$ associates to $\psi\in[\dop,[\dnomen,\ccl]]$ the object
$\xi'_k\psi\in[\dnop,\ccl]$ given, for all $([p_1]\ldots[p_n])\in\dnop$, by
$(\xi'_k\psi)(p_1\ldots p_n)=\psi_{p_k}(p_1\ldots p_{k-1}\,p_{k+1}\ldots p_n)$.
These isomorphisms will be important in analysing the simplicial directions in
a $\ctn$ (see Section \ref{simpdir}).

An object $\phi\in[\dnop,\ccl]$ is said to satisfy the \emph{globularity
condition} (or to be a \emph{globular object}) if $\phi(0,\mi):\dnomen\rw\ccl$
is constant and if $\phi(m_1\ldots m_r\,0\,\mi):\Delta^{n-r-1^{op}}\rw\ccl$ is
constant for all $[m_i]\in\dop$, $i=1,\ldots, r$, $1\leq r \leq n-2$.

Simplicial (resp. multi-simplicial) objects often arise as nerves (resp.
multi-nerves) of internal categories (resp. internal $n$-fold categories).
Suppose $\ccl$ has finite limits. An \emph{internal category} $\abb$ in $\ccl$
is a diagram in $\ccl$
\begin{equation}\label{simtec.eq1}
    \tens{A_1}{A_0}\rTo^c A_1\pile{\rTo^{d_0}\\ \rTo^{d_1}\\ \lTo_s}A_2
\end{equation}
where $c,d_0,d_1,s$ satisfy the usual axioms of a category (see for instance
\cite{borc} for details). An \emph{internal functor} is a morphism of diagrams
like (\ref{simtec.eq1}) which commutes in the obvious way. We denote by
$\cat\ccl$ the category of internal categories and internal functors.

There is a \emph{nerve} functor $\ner:\cat\ccl\rw[\dop,\ccl]$ such that, given
$\abb$ as in (\ref{simtec.eq1}), $(\ner\abb)_0= A_0$, $(\ner\abb)_1=A_1$ and,
for each $k\geq 2$, $(\ner\abb)_k=\pro{A_1}{A_0}{k}$ is the limit of the
diagram $A_1\supar{d_1}A_0\suparle{d_0}\cdots\supar{d_1}A_0\suparle{d_0}A_1$.
The image of the nerve functor can be characterized as follows. Given any
$\phi\in[\dop,\ccl]$, it is easy to see from the definition of limit that for
each $k\geq 2$ there are unique maps $\eta_k:\phi_k\rw\fistr{k}$ such that
$\pr_j\eta_k=\nu_j$, $1\leq j\leq k$, where $\pr_j$ is the $j^{th}$ projection
and $\nu_j:A_k\rw A_1$ is induced by the map in $\Delta\;$ $[1]\rw[k]$ sending
$0$ to $j-1$ and $1$ to $j$. The maps $\eta_k$ are called \emph{Segal maps} and
they play an important role in the Tamsamani as well as in the $\ctn$ model.
Segal maps can be used to characterize nerves of internal categories, as in the
following proposition.
\begin{proposition}\label{simtec.pro1}
    A simplicial object in $\ccl$ is the nerve of an internal category in
    $\ccl$ if and only if the Segal maps $\eta_k$ are isomorphisms for all $k\geq
    2$.
\end{proposition}
The category $\cat^n(\ccl)$ of \emph{$n$-fold categories in $\ccl$} is defined
inductively by iterating $n$ times the internal category construction. That is,
$\cat^1(\ccl)=\cat\ccl$ and, for $n>1$, $\cat^n(\ccl)=\cat(\cat^\nm(\ccl))$.
When $\ccl=\gp$ this gives the category of $\ctn$s. By iterating the nerve
construction one obtains a \emph{multinerve} functor
$\ncl:\cat^n(\ccl)\rw[\dnop,\ccl]$. A characterization of the image of $\ncl$
in terms of Segal maps is possible (see Lemma \ref{simpdir.lem1} in the case
$\ccl=\gp$). If $\phi\in\cat^n(\ccl)$, $\ncl\phi$ does not in general satisfy
the globularity condition. It can be shown that if it does then $\phi\in n
\mi\cat(\ccl)$ where $n\mi\cat$ denotes strict $n$-categories.

We will use extensively the notions of classifying space of multi-simplicial
sets and multi-simplicial groups, thus we recall the main definitions and
constructions. The classifying space of an $n$-simplicial set is defined by the
composite $B:[\dnop,\set]\supar{\diag}[\dop,\set]\supar{|\cdot |}\tp$, where
$|\cdot |$ is the geometric realization and $\diag$ is the diagonal functor,
defined by $(\!\diag\phi)_r=\phi(r,r,\ldots,r)$. The following observation will
be used when discussing the classifying space of multi-simplicial groups.
\begin{remark}\rm\label{simtec.rem1}
    For each $\phi\in[\dop,[\dnomen,\set]]$, $\diag\xi'_k\phi$ is independent
    of $k$. In fact, from the expression of $\xi'_k$ given above, we have
    $\diag(\xi'_k\phi)[r]=(\xi'_k\phi)(r,...,r)=\phi_r(r,\ldots,r)$.
\end{remark}
The classifying space of a $n$-simplicial group is the composite
\begin{equation*}
    B:[\dnop,\gp]\supar{\diag}[\dop,\gp]\supar{\ovl{W}}[\dop,\set]\supar{|\cdot
    |}\tp
\end{equation*}
where $\ovl{W}$ is right adjoint to the Kan loop group functor.

The classifying space of a $n$-simplicial group can also be computed, up to
homotopy, in two alternative ways, which are the ones we mostly use in this
paper, and which are as follows. Recall that a group can be considered as a
category with just one object, hence there is a nerve functor
$\ner:\gp\rw[\dop,\set]$. Then the classifying space of a $n$-simplicial group
can be computed, up to homotopy, by taking the nerve of the group at each
dimension and then taking the classifying space of the resulting
$(\np)$-simplicial set. Formally we have the composite
\begin{equation*}
    B:[\dnop,\gp]\supar{\ovl{\ner}}[\dnop,[\dop,\set]]\supar{\xi'_k}
    [\Delta^{\np^{op}},\set]\supar{\diag}[\dop,\set]\supar{|\cdot
    |}\tp.
\end{equation*}
This method of computing the classifying space will be used in the proofs of
Theorems \ref{glob2.the1}, \ref{globgen.the1}, \ref{tam3.the1} and Proposition
\ref{semist.pro01}. Notice that, by Remark \ref{simtec.rem1}, this composite is
independent of $k$.

The second method to compute the classifying space of a $n$-simplicial group is
as follows. Given $\phi\in[\dnop,\gp]$, we first regard $\phi$ as a simplicial
object in $[\dnomen,\gp]$ along a fixed direction; that is, we consider
$\xi_k\phi$. Then we apply the classifying space functor
$B:[\dnomen,\gp]\rw\tp$ levelwise; that is, we consider $\ovl{B}\xi_k\phi$; we
obtain in this way a simplicial space. There is a standard notion of geometric
realization of a simplicial space. Then the classifying space of $\phi$ is, up
to homotopy, the geometric realization of this simplicial space; that is,
$|\ovl{B}\xi_k\phi|$. In connection with this we will also use the following
fact. Recall that a \emph{weak equivalence} in $\tp$ is a map $f:X\rw Y$
inducing isomorphisms $\pi_n(X,x)\cong\pi_n(Y,f(x))$ for all $x\in X$, $n\geq
0$. If a map of simplicial spaces is levelwise a weak equivalence then this map
induces a weak equivalence of geometric realizations. This fact, together with
this second method of computing the classifying space, will be used in the
proofs of Lemmas \ref{catn.lem3} and \ref{spec.lem3}.

We will also use extensively the notion of homotopy groups $\pi_n(G_*)$ of a
simplicial group $G_*$. These are independent of the basepoint and coincide
with the homotopy group of the underlying simplicial set, which is always a
\emph{fibrant} simplicial set; that is, it satisfies the Kan condition.
Further, the homotopy groups of a simplicial group can be calculated via an
algebraic device called the \emph{Moore complex}. This is the (not necessarily
abelian) chain complex $N_*$
\begin{equation*}
    \cdots\supar{}N_2\supar{\pt_2} N_1\supar{\pt_1}N_0\supar{}1
\end{equation*}
where $N_n(G_*)=\{x\in G_n\;|\;\pt_i x=1\text{ for all }i\neq n\}$. Then
$\pi_n(G_*)$ is the $n^{th}$ homology group of $N_*$. This will be useful, for
instance, in calculating the homotopy groups of the nerve of a cat$^1$-group.

Finally we will use a model structure on simplicial groups due to Quillen, in
which a map $f:G_*\rw H_*$ is a weak equivalence if and only if it induces
isomorphisms $\pi_n(G_*)\cong\pi_n(H_*)$ for all $n\geq 0$. We  recall that if
$G_*\in[\dop,\gp]$, $\pi_0 BG_*=\{\cdot\}$ and $\pi_i BG_*=\pi_{i-1}G_*$ for
$i\geq 1$. Thus $f$ is a weak equivalence if and only if $Bf$ is a weak
equivalence of spaces. In this model structure $f$ is a fibration if and only
if $N_q f:N_q G_*\rw N_q H_*$ is surjective for $q>0$. In particular every
simplicial group is fibrant in this model structure, and every surjection of
simplicial groups is a fibration. This model structure will be used in the
proofs of Lemma \ref{spec.lem1} and Theorems \ref{glob2.the1} and
\ref{globgen.the1}.

\section{Tamsamani's model.}\label{pre}

\subsection{Tamsamani weak $n$-categories.}\label{formal}
We recall the definition of Tamsamani weak $n$-category and Tamsamani weak
$n$-groupoid. We follow closely the treatment given in \cite{toe}; see also
\cite{tam}, \cite{tam1}.

We first describe informally the idea behind this notion. Recall that the
category $n$-cat of strict $n$-categories is by definition the category
$((\nm)\mi\mathrm{cat})\mi\mathrm{cat}$ of categories enriched in $(\nm)$-cat
with respect to the cartesian monoidal structure. There is a general result of
Ehresmann \cite[Appendix]{ehr} that shows that, under sufficiently general
conditions on a category $\ccl$, categories enriched in $\ccl$ with respect to
the cartesian monoidal structure correspond to internal categories in $\ccl$ in
which the object of objects is discrete (that is, it is the coproduct of copies
of the terminal object). It is not hard to see that the hypotheses of
\cite{ehr} are satisfied when $\ccl=(\nm)$-cat. It follows that strict
$n$-categories correspond to internal categories in $(\nm)$-cat in which the
object of objects is discrete. From the characterization of nerves of internal
categories of Proposition \ref{simtec.pro1}, the latter correspond to objects
$\phi\in[\dop,(\nm)\mi\mathrm{cat}]$ such that $\phi_0$ is discrete and the
Segal maps are isomorphisms. This is how strict $n$-categories are described in
the Tamsamani model. For the weak case, the idea is to relax this structure by
demanding that the Segal maps are no longer isomorphisms but some suitably
defined $(n-1)$-equivalences. Thus denoting by $\wcl_n$ the category of
Tamsamani weak $n$-categories, $\wcl_1=\cat$ and 1-equivalences are
equivalences of categories; inductively, having defined $\wcl_\nm$ and
$(\nm)$-equivalences, we say $\phi\in\wcl_n\subset[\dop,\wcl_\nm]$ if $\phi_0$
is discrete and the Segal maps are $(n-1)$-equivalences. The point of the
formal definition is then to define $n$-equivalences, so as to complete the
inductive step.

We now give the formal definition of $\wcl_n$. In what follows, the ``bar"
notation is as in Definition \ref{simtec.def1}. For $n=1$, let $\wcl_1=\cat$. A
morphism in $\wcl_1$ is a 1-equivalence if it is an equivalence of categories.
Let $\tau^{(1)}_0:\cat\rw\set$ be the functor which associates to a category
the set of isomorphism classes of its objects. Let $\delta^{(1)}:\set\rw\wcl_1$
be the functor which associates to a set $X$ the discrete category with
object-set $X$. Let $\tau^{(1)}_1:\wcl_1\rw\cat$ be the identity functor. We
notice the following properties:
\begin{itemize}
  \item [i)] $\delta^{(1)}$ is fully faithful and finite product-preserving and
  $\tau^{(1)}_0\delta^{(1)}=\id_{Set}$,
  \item [ii)] $\tau^{(1)}_0$ sends 1-equivalences to bijections,
  \item [iii)] $\tau^{(1)}_0$ preserves fibre products over discrete objects,
  \item [iv)] if $\ccl\rw\dcl$ is a morphism in $\wcl_1$ with $\dcl$ discrete, then
  $\ccl\cong\underset{x\in\dcl}{\coprod}\ccl_x$.
\end{itemize}
Suppose, inductively, that we have defined the category $\wcl_\nm$ as well as
$(\nm)$-equivalences between objects of $\wcl_\nm$. Suppose also that there
exists a fully faithful and finite product-preserving functor
$\delta^{(\nm)}:\set\rw\wcl_\nm$. Objects in the image of $\delta^{(\nm)}$ are
called \emph{discrete}. Suppose there exists $\tau^{(\nm)}_0:\wcl_\nm\rw\set$
such that
\begin{itemize}
  \item [i)] $\delta^{(\nm)}$ is fully faithful and finite product-preserving and
  $\tau^{(\nm)}_0\delta^{(\nm)}\cong\id_{Set}$,
  \item [ii)] $\tau^{(\nm)}_0$ sends $(\nm)$-equivalences to bijections,
  \item [iii)] $\tau^{(\nm)}_0$ preserves fibre products over discrete objects,
  \item [iv)] if $\phi\rw\psi$ is a morphism in $\wcl_\nm$ with $\psi$ discrete, there is an
  isomorphism $\phi\cong\underset{x\in\psi}{\coprod}\phi_{x}$.
\end{itemize}
We define the category $\wcl_n$ as the full subcategory of functors
$\phi\in[\dop,\wcl_\nm]$ such that $\phi_0$ is discrete and the Segal maps
$\phi_k\rw\phi_1\tiund{\phi_0}\overset{k}{\cdots}\tiund{\phi_0}\phi_1$,
$\;k\geq 2$, are $(\nm)$-equivalences.

To complete the inductive step, we need to define $n$-equivalences as well as
functors $\delta^{(n)}$, $\tau_0^{(n)}$ and verify (i) - (iv) at step $n$. For
this we consider the functor
$\ovl{\tau^{(\nm)}_0}:[\dop,\wcl_\nm]\rw[\dop,\set]$ induced by applying
$\tau^{(\nm)}_0$ levelwise as in Definition \ref{simtec.def1}; that is, for
each $[k]\in\dop$, $(\ovl{\tau^{(\nm)}_0}\phi)_k=\tau^{(\nm)}_0\phi_k$. Notice
that the restriction of this functor to $\wcl_n$ yields a functor
$\ovl{\tau^{(\nm)}_0}:\wcl_n\rw\ner(\cat)$ where $\ner(\cat)$ is the full
subcategory of those simplicial sets which are in the image of the nerve
functor $\ner:\cat\rw[\dop,\set]$. This is because, if $\phi\in\wcl_n$, the
Segal maps
$\phi_k\rw\phi_1\tiund{\phi_0}\overset{k}{\cdots}\tiund{\phi_0}\phi_1$ are
$(\nm)$-equivalences, hence, by (ii) and (iii), the maps
\begin{equation*}
    \tau^{(\nm)}_0\phi_k\rw
    \tau^{(\nm)}_0(\phi_1\tiund{\phi_0}\overset{k}{\cdots}\tiund{\phi_0}\phi_1)\cong
    \tau^{(\nm)}_0\phi_1\tiund{\tau^{(\nm)}_0\phi_0}\overset{k}{\cdots}
    \tiund{\tau^{(\nm)}_0\phi_0}\tau^{(\nm)}_0\phi_1
\end{equation*}
are bijections. Composing $\ovl{\tau^{(\nm)}_0}$ with the isomorphism
$\nu:\ner(\cat)\cong\cat$ yields a functor
\begin{equation}\label{formal.eq1}
\begin{split}
   & \hspace{5cm}\tau_1^{(n)}:\wcl_n\rw\cat \\
    & \text{with}\hspace{12cm}\\
    &\hspace{2cm}\ner\tau_1^{(n)}=\ovl{\tau^{(\nm)}_0}_{|_{\wcl_n}}, \quad
    \tau_1^{(n)}=\nu\,\ovl{\tau^{(\nm)}_0}_{|_{\wcl_n}}.
\end{split}
\end{equation}
Consider the morphism $(\pt_0,\pt_1):\phi_1\rw\phi_0\times\phi_0$. Since
$\phi_0\times\phi_0$ is discrete (as $\delta^{(\nm)}$ preserves finite
products), by iv), $\phi_1\cong\underset{x,y\in\phi_0}{\coprod}\phi_{(x,y)}$. A
morphism $f:\phi\rw\psi$ in $\wcl_n$ is defined to be an $n$-equivalence if,
for all $x,y\in\phi_0$, $\phi_{(x,y)}\rw\psi_{(fx,fy)}$ is an
$(n-1)$-equivalence and $\tau_1^{(n)}\phi\rw\tau_1^{(n)}\psi$ is an equivalence
of categories.

Let $d:\wcl_\nm\rw[\dop,\wcl_\nm]$ take every object $\phi\in\wcl_\nm$ to the
constant simplicial object $d\phi\in[\dop,\wcl_\nm]$; that is,
$(d\phi)_k=\phi$. Then $d\delta^{(\nm)}:\set\rw\wcl_n$. Let
\begin{align*}
   & \tau^{(n)}_0:\wcl_n\rw\set \qquad\qquad  \tau^{(n)}_0=\tau^{(1)}_0\tau^{(n)}_1\\
   & \delta^{(n)}:\set\rw\wcl_n \qquad\qquad  \delta^{(n)}=d\delta^{(\nm)}.
\end{align*}
We check that conditions (i) -- (iv) hold for $n$ in place of $\nm$. Given a
set $X$, by the induction hypothesis (i), $(\ner\tau^{(n)}_1 d \delta^{(\nm)}
X)_k=\tau^{(\nm)}_0\delta^{(\nm)} X = X$ for all $k$, so that $\tau^{(n)}_1 d
\delta^{(\nm)}X=\delta^{(1)}X$. Hence
$\tau^{(n)}_0\delta^{(n)}X=\tau^{(1)}_0\tau^{(n)}_1 d \delta^{(\nm)}
X=\tau^{(1)}_0\delta^{(1)}X=X$. Also, $\delta^{(n)}$ is fully faithful and
finite product-preserving, since $\delta^{(\nm)}$ and $d$ are. Thus (i) holds
for $n$. If $f$ is an $n$-equivalence, $\tau^{(n)}_1 f$ is an equivalence of
categories, hence $\tau^{(n)}_0 f=\tau^{(1)}_0 \tau^{(n)}_1 f$ is a bijection,
proving (ii). Notice that $\tau^{(n)}_1$ preserves fibre products over discrete
objects. For instance, consider $\phi\rTo^f\dcl\lTo^g\phi$ in $\wcl_n$, with
$\dcl$ discrete. By the induction hypothesis (iii), for all $k$,
\begin{align*}
  & (\ner\tau^{(n)}_1(\phi\tiund{\dcl}\psi))_k=\tau^{(\nm)}_0(\phi_k\tiund{\dcl_k}\psi_k)=\\
  & \tau^{(\nm)}_0\phi_k\tiund{\tau^{(\nm)}_0\dcl_k}\tau^{(\nm)}_0\psi_k=
  (\ner\tau^{(n)}_1\phi)_k\tiund{(\ner\tau^{(n)}_1\dcl)_k}(\ner\tau^{(n)}_1\psi)_k.
\end{align*}
This shows that $\ner\tau^{(n)}_1$, and therefore $\tau^{(n)}_1$, preserve
fibre products over discrete objects. Since $\tau^{(n)}_0 =\tau^{(1)}_0
\tau^{(n)}_1$ and $\tau^{(1)}_0$ also preserves fibre products over discrete
objects, it follows that (iii) holds for $n$. Finally, if $\phi\rw\psi$ is a
morphism in $\wcl_n$ with $\psi$ discrete, then for each $k$, $\phi_k\rw\psi_k$
is a morphism in $\wcl_\nm$ and $\psi_k=\psi_0$ is discrete. Thus, by the
induction hypothesis iv),
$\phi_k\cong\underset{x\in\psi_0}{\coprod}(\phi_k)_{x}$. It follows that
$\phi\cong\underset{x\in\psi_0}{\coprod}\phi_{x}$ where
$(\phi_{x})_k=(\phi_k)_{x}$.

\subsection{Weak ${n}$-groupoids}\label{group}
The definition of \tam \emph{weak $n$-groupoids} $\tcl_n$ can now be given
inductively as follows. We set $\tcl_1=\gpd$, the category of groupoids in the
ordinary sense; that is, categories in which every arrow is invertible. Suppose
we have defined the category $\tcl_\nm$ of weak $(\nm)$-groupoids. Then an
object $\phi$ of $\wcl_n$ is a weak $n$-groupoid if, for all $x,y\in\phi_0$,
$\phi_{(x,y)}$ is a weak $(\nm)$-groupoid and $\tau^{(n)}_1\phi$ is a groupoid.

Tamsamani proved that weak $n$-groupoids model $n$-types. Let $\tp^{(n)}$ be
the category of $n$-types, that is, objects of $\tp^{(n)}$ are topological
spaces $X$ with $\pi_i(X,x)=0$ for all $i>n$ and $x\in X$. By unravelling the
inductive definition of $\tcl_n$ and taking the nerve functor
$\ner:\gpd\rw[\dop,\set]$ one obtain a full and faithful embedding
$\ncl:\tcl_n\rw[\dnop,\set]$. The classifying space functor
$B:\tcl_n\rw\tp^{(n)}$ is obtained by composing $\ncl$ with the classifying
space functor $B:[\dnop,\set]\rw\tp$ as in Section \ref{simtec}. Let $\hcl
o\;(\tp^{(n)})$ be the category obtained by formally inverting the weak
equivalences, that is the morphisms in $\tp^{(n)}$ which induce isomorphisms of
homotopy groups. Then Tamsamani constructed a fundamental weak $n$-groupoid
functor $\pp_n:\tp\rw\tcl_n$ and proved the following theorem:
\begin{theorem}\label{s1.the1}\rm\cite[Theorem 8]{tam}\it
    The functors $B$ and $\pp_n$ induce an equivalence of categories
\begin{equation*}
  \tcl_n\bsim^{n}\;\;\simeq\;\;\hcl o\;(\tp^{(n)})
\end{equation*}
between the localization of $\tcl_n$ with respect to the $n$-equivalences and
the homotopy category of $n$-types.
\end{theorem}
We also record the following fact, which is an immediate consequence of the
properties of the functor $\pp_n$ proved in \cite{tam}, \cite{tam1}. We will
use the following proposition in the proof of our semistrictification result
Theorem \ref{result.the1}.
\begin{proposition}\label{group.pro1}
    Let $\gcl,\gcl'\in\tcl_n$ and suppose that there is a zig-zag of weak equivalences in $\mathrm{Top}$ between $B\gcl$ and $B\gcl'$. Then
there is a zig-zag of $n$-equival-ences in $\tcl_n$ between $\gcl$ and $\gcl'$.
\end{proposition}
\prf Let us denote the zig-zag of weak equivalences in $\mathrm{Top}$ between
$B\gcl$ and $B\gcl'$
\begin{equation}\label{group.eq1}
    B\gcl\rw X_1\lw X_2\rw\cdots\rw X_i\lw\cdots\rw B\gcl'.
\end{equation}
It is proved in \cite[Proposition 11.4 (i)]{tam} that the functor $\pp_n$ sends
weak equivalences to $n$-equivalences. Hence applying $\pp_n$ in
(\ref{group.eq1}) we obtain a zig-zag of $n$-equivalences in $\tcl_n$
\begin{equation}\label{group.eq2}
    \pp_n B\gcl\rw \pp_n X_1\lw \pp_n X_2\rw\cdots\rw \pp_n X_i\lw\cdots\rw \pp_n B\gcl'.
\end{equation}
It is also proved in \cite[Proposition 11.4 (ii)]{tam} that for each
$\gcl\in\tcl_n$ there is an $n$-equivalence $\gcl\rw\pp_n B\gcl$, natural in
$\gcl$. From (\ref{group.eq2}) we therefore obtain the following zig-zag of
$n$-equivalences between $\gcl$ and $\gcl'$
\begin{equation*}
   \gcl\rw\pp_n B\gcl\rw \pp_n X_1\lw \pp_n X_2\rw\cdots\rw
   \pp_n X_i\lw\cdots\rw \pp_n B\gcl'\lw\gcl'. \hspace{9mm}\Box
\end{equation*}

\subsection{Some basic lemmas}\label{teclem}
We prove some basic technical lemmas about the Tamsamani model. These will not
be needed until Section \ref{tamwe}, so they can be skipped at a first reading.
Lemma \ref{teclem.lem1} is a simple calculation using the definition of
$\wcl_n$ and will be used in the proof of Proposition \ref{semist.pro01},
leading to one of our main theorem (Theorem \ref{semist.the01}). Lemma
\ref{teclem.lem2} is used to prove Corollary \ref{teclem.cor1}, which is a
basic fact used in Section \ref{scontra}, as well as in Remark
\ref{teclem.rem1}. The latter is used in the proof of Lemma \ref{semist.lem01}.

Let $\ovl{\tau^{(\nm)}_1}:[\dop,\wcl_\nm]\rw[\dop,\cat]$ be induced by
$\tau^{(\nm)}_1:\wcl_\nm\rw\cat$ as in Definition \ref{simtec.def1}.
\begin{lemma}\label{teclem.lem1}
    For each $\phi\in\wcl_n$, $n>2$, $\tau^{(n)}_1\phi=\tau^{(2)}_1\ovl{\tau^{(\nm)}_1}\phi$.
\end{lemma}
\prf For each $r$, $(\ovl{\tau^{(\nm)}_1}\phi)_r=\tau^{(\nm)}_1\phi_r$. We have
$\phi_0=\delta^{(\nm)}X=d\delta^{(n-2)}X$, so that, for each $k$,
\begin{equation*}
  (\ner\tau^{(\nm)}_1\phi_0)_k=(\ovl{\tau^{(n-2)}_0}\phi_0)_k=\tau^{(n-2)}_0(d\delta^{(n-2)}X)_k=
    \tau^{(n-2)}_0\delta^{(n-2)}X=X.
\end{equation*}
Hence $\tau^{(\nm)}_1\phi_0$ is a discrete category. Since the Segal maps
$\phi_k\!\rw\!\phi_1\!\tiund{\phi_0}\overset{k}{\cdots}\!\tiund{\phi_0}\phi_1$
are $(\nm)$-equivalences and, as observed earlier, $\tau^{(\nm)}_1$ preserves
fibre products over discrete objects, we have an equivalence of categories
\begin{equation*}
  \tau^{(\nm)}_1\phi_k\simeq\tau^{(\nm)}_1(\fistr{k})
  =\tau^{(\nm)}_1\phi_1\tiund{\tau^{(\nm)}_1\phi_0}\overset{k}{\cdots}\tiund{\tau^{(\nm)}_1
  \phi_0}\tau^{(\nm)}_1\phi_1;
\end{equation*}
this proves that $\overline{\tau_1^{(\nm)}}\phi\in\wcl_2$. It follows that

$\tau^{(n)}_1\phi=\nu\,\ovl{\tau^{(\nm)}_0}\phi=\nu\,\ovl{\tau^{(1)}_0\tau^{(n-1)}_1\phi}=
\nu\,\ovl{\tau^{(1)}_0}\;\,\ovl{\tau^{(\nm)}_1}\,\phi={\tau}^{(2)}_1\,\ovl{\tau^{(\nm)}_1}\,\phi$.
\enpr

\begin{lemma}\label{teclem.lem2}
    Let $f:\phi\rw\psi$ be an $n$-equivalence in $\wcl_n$ and suppose that
    $\psi\in\tcl_n$. Then $\phi\in\tcl_n$.
\end{lemma}
\prf We use induction on $n$. For $n=1$ the lemma amounts to saying that if
$f:\phi\rw\psi$ is an equivalence of categories and $\psi$ is a groupoid, then
$\phi$ is itself a groupoid. To show this, let $g:\psi\rw\phi$ be a pseudo
inverse and let $\alpha$ be the natural isomorphism $\alpha:\id\cong gf$. For
each arrow $u:x\rw y$ in $\phi$ we have by naturality that $\alpha_y
u=(gfu)\alpha_x$. Since $\psi$ is a groupoid, $gfu$ has an inverse $h$. The
inverse for $u$ is then provided by the composite $\alpha_x^{-1}h\alpha_y$, as
is easily checked.

Suppose the lemma holds for $(n-1)$ and let $f:\phi\rw\psi$ be an
$n$-equivalence in $\wcl_n$ with $\psi\in\tcl_n$. By definition, for each
$x,y\in\phi_0$, the map $\phi_{(x,y)}\rw\psi_{(fx,fy)}$ is an
$(n-1)$-equivalence, and the map $\tau^{(n)}_1\phi\rw\tau^{(n)}_1\psi$ is an
equivalence of categories. Since $\psi\in\tcl_n$, $\psi_{(fx,fy)}$ is an
$(n-1)$-groupoid and $\tau^{(n)}_{1}\psi$ is a groupoid. From the induction
hypothesis and the case $n=1$ it follows that $\phi_{(x,y)}\in\tcl_{n-1}$ and
$\tau^{(n)}_{1}\phi\in\mathrm{Gpd}$. Hence $\phi\in\tcl_n$. \enpr
\begin{corollary}\label{teclem.cor1}
    For each $n\geq 2$ the embedding $\wcl_n\subset[\dop,\wcl_{\nm}]$ restricts
    to an embedding $\tcl_n\subset[\dop,\tcl_{\nm}]$.
\end{corollary}
\prf Let $\phi\in\tcl_n\subset\wcl_n\subset[\dop,\wcl_{\nm}]$ and put
$\phi_k=\phi[k]$. Clearly $\phi_0\in\tcl_{\nm}$ since $\phi_0$ is discrete.
Also, $\phi_1\cong\underset{x,y\in\phi_0}{\amalg}\phi_{(x,y)}\in\tcl_\nm$ since
each $\phi_{(x,y)}\in\tcl_\nm$ as $\phi\in\tcl_n$. For each $k\geq 2$, the
Segal map $\phi_k\rw\fistr{k}$ is an $(\nm)$-equivalence. It is immediate that,
since $\phi_1,\phi_0\in\tcl_\nm$, we have $\fistr{k}\in\tcl_\nm$. Lemma
\ref{teclem.lem2} then implies that $\phi_k\in\tcl_\nm$ for $k\geq 2$. We
conclude that $\phi_k\in\tcl_\nm$ for all $k\geq 0$.  \enpr
\begin{remark}\rm\label{teclem.rem1}
    By Corollary \ref{teclem.cor1}, $\tcl_n\subset[\dop,\tcl_\nm]$
    and $\tcl_1=\gpd$. It follows that $\tcl_n$ is a full subcategory
    of $\underset{\nm}{\underbrace{[\dop,[\dop,\ldots,[\dop}},\gpd]\cdots]$.
    Since the latter is isomorphic to $[\dnomen,\gpd]$, we also have a fully
    faithful embedding $\tcl_n\rw[\dnomen,\gpd]$.
\end{remark}

\section{The category of cat$\bsy{^{n}}$-groups}\label{catcat}
The goal of this section is to prove two facts about the category of $\ctn$s,
which are Theorem \ref{adj.the1} and Proposition \ref{twolem.pro2}. The first
one asserts the existence of an adjunction $\fl_n\dashv\ul_n$ between $\ctn$s
and $\set$ and gives an inductive explicit description of $\fl_n(X)$. The
second says that if $\ul_n f$ is surjective then the map of multinerves $\ncl
f$ is levelwise surjective. As outlined in Section \ref{orga}, both these facts
are essential in the construction of the functor $\spi$ in Section
\ref{specia}.

Theorem \ref{adj.the1} and Proposition \ref{twolem.pro2} are deduced from
analogous results (Theorem \ref{adj.the0} and Proposition \ref{catn.pro1}
respectively) for the category $\cng$ which is well known to be equivalent to
$\cagp$ (Corollary \ref{cng.cor1}). The reason for passing to the category
$\cng$ is that this is a category of groups with operations. In this context we
have the formalism of extensions and semidirect products, and internal
categories in categories of groups with operations are particularly easy to
describe, as they are equivalent to reflexive graphs satisfying an extra
condition (Lemma \ref{oper.lem1}). We have recalled these basic notions at the
beginning of Section \ref{oper}.

In the general context of a category of groups with operations of the type
described at the beginning of Section \ref{oper}, we prove (Theorem
\ref{oper.the1}) the existence of an adjunction $\lcl:\dcl\pile{\lTo\\
\rTo}\cat\dcl:\vcl$, $\lcl\dashv\vcl$. This adjunction is one of the main tools
used in constructing the adjunction $\fcl_n\dashv\ucl_n$ of Theorem
\ref{adj.the0} between $\cng$ and $\set$. The idea is to apply Theorem
\ref{oper.the1} to the case $\dcl=\cngm$ to build inductively
$\ucl_n:\cng\rw\set$ as the composite
\begin{equation*}
    \cng\supar{\beta^{(k)}_n}\cat\cngm\supar{\vcl_n}\cngm\supar{\ucl_\nm}\set.
\end{equation*}
Here $\beta^{(k)}_n$ is the equivalence of categories of Corollary
\ref{cgn.cor0}, which is a standard fact proved from the basic Lemmas
\ref{oper.lem1} and \ref{cgn.lem1}.

We notice that the adjunction between $\cng$ and $\set$ of Theorem
\ref{adj.the0} is closely related to some results already in the literature.
Namely, for $n=1$ this adjunction is studied in \cite{ccg} using the language
of crossed modules in groups (which is a category equivalent to $\ccl^1\gcl$).
For $n\geq 1$, the adjunction between $\cng$ and $\set$ is studied in
\cite{celp}. However, in the case $n>1$ the result of \cite[Lemma 9]{celp}
gives a description of $\fcl_n(X)$ which is less precise than the one we give
in Theorem \ref{adj.the0}. On the other hand our description of $\fcl_n(X)$,
and the corresponding one for $\fl_n(X)$ in Theorem \ref{adj.the1} is needed in
the proof of Theorem \ref{spec.the0} where we show that $\fl_n(X)$ is strongly
contractible. We also mention that there exists in the literature another
equivalent description of $\cng$ in terms of the category $Crs^n$ of crossed
$n$-cubes \cite{elst}. In \cite{cihl} they study an adjunction between $Crs^n$
and $\set$ which corresponds precisely to the above adjunction between $\cng$
and $\set$. In this paper we have chosen to avoid the complex formalism of
crossed $n$-cubes, so we give an independent self-contained account of the
adjunction $\fcl_n\dashv\ucl_n$ using theorem \ref{oper.the1}.

We start this section with an account of well known facts about $\ctn$s
(Sections \ref{catn} and \ref{simpdir}) which will be needed later. We have
provided in this section many of the proofs of background facts concerning
$\ctn$s. Further references may be found for instance in \cite{brr},
\cite{bcd}, \cite{ccg}, \cite{celp}, \cite{cihl}, \cite{elst}, \cite{lod},
\cite{orz}, \cite{pasu}, \cite{por1}, \cite{por}.

\subsection{Cat${^n}$-groups as homotopy models}\label{catn}
The notion of $\ctn$ is obtained starting from the category of groups and
iterating $n$ times the internal category construction:
\begin{definition}\label{catn.def1}
    The category $\cagp$ of cat$^{n}$-groups is defined inductively by $\cat^0(\gp)=\gp$,
    and, for $n\geq 1$, $\cagp=\cat(\cagpm)$.
\end{definition}
As recalled in Section \ref{simtec}, by iterating the nerve construction one
obtains the \emph{multinerve functor}
\begin{equation*}
  \ncl:\cagp\rw[\dnop,\gp].
\end{equation*}
The functor $\ncl$ has a left adjoint and is fully faithful \cite{bcd} . The
\emph{classifying space} of a $\ctn$ is the classifying space of its
multinerve, where the classifying space of a $n$-simplicial group is as in
Section \ref{simtec}. We say a morphism $f$ in $\cagp$ is a \emph{weak
equivalence} if $Bf$ is a weak equivalence of topological spaces.

It can be shown \cite{bcd} that the classifying space of a $\ctn$ is a
connected $(n+1)$-type. Further, there is a functor $\pcl_n:\tp\rw \cagp$
called the fundamental $\ctn$ functor, such that the following theorem holds:
\begin{theorem}\label{catn.the1}\cite{bcd}
    The functors $B$ and $\pcl_n$ induce an equivalence of categories
\begin{equation*}
  \cagp\bsim\;\;\simeq\;\;\hcl o\;(\tp_{*}^{(n+1)}).
\end{equation*}
between the localization of $\cagp$ with respect to the weak equivalences and
the homotopy category of connected $(n+1)$-types.
\end{theorem}
Theorem \ref{catn.the1} will be used in the proof of the semistrictification
result Theorem \ref{result.the1}. We will also need the following well known
elementary facts about $\cgp$. In particular, Lemma \ref{catn.lem1} (iii) will
be important in Section \ref{tam3} in defining the functor $\tiu:\dbb_2\rw\cgp$
and in establishing its properties (Corollary \ref{tam3.cor1}). Given
$\gcl\in\cgp$ let $\pi_i(\gcl)=\pi_i(\ner\gcl)$, where $\ner:\cgp\rw[\dop,\gp]$
is the nerve functor. An object $\gcl\in\cgp$ is \emph{discrete} if $\ner\gcl$
is constant.
\begin{lemma}\label{catn.lem1}
    Let $\gcl\in\cgp$ be given by
    \begin{equation*}
        \tens{G_1}{G_0}\rTo G_1\pile{\rTo^{\pt_0}\\ \rTo^{\pt_1}\\
        \lTo_{\sigma_0}} G_0
    \end{equation*}
then
\begin{itemize}
  \item[(i)] $\pi_0(\gcl)=G_0/\pt_1(\ker\pt_0)$,
  $\pi_1(\gcl)=\ker{\pt_1}_{|_{\ker\pt_0}}$, $\pi_i(\gcl)=0$ for $i>1$;
  \item[(ii)] a map $f:\gcl\rw\gcl'$ in $\cgp$ is a weak equivalence if and
  only if it induces isomorphisms $\pi_0(\gcl)\cong\pi_0(\gcl')$,
  $\pi_1(\gcl)\cong\pi_1(\gcl')$;
  \item[(iii)] the functor $\pi_0:\cgp\rw\gp$ preserves fibre products over
  discrete objects and sends weak equivalences to isomorphisms.
\end{itemize}
\end{lemma}
\prf (i) As recalled in Section \ref{simtec}, the homotopy groups of the
simplicial group $\ner\gcl$ are the homology groups of the Moore complex. The
result follows from a straightforward computation using this fact.

(ii) As recalled at the end of Section \ref{simtec}, a map $g$ of simplicial
groups is a weak equivalence, that is it induces isomorphisms of homotopy
groups, if and only if $Bg$ is a weak equivalence of spaces. Applying this fact
to $\ner f$, the result follows.

(iii) Consider the diagram $\gcl\supar{f}\hcl\suparle{g}\hcl'$, where $\hcl$ is
discrete; we have
\begin{equation*}
    \begin{split}
       & \pi_0(\gcl\tiund{\hcl}\gcl')=G_0\tiund{H_0}G'_0/(\pt_1\ker\pt_0, \pt'_1\ker\pt'_0)\\
        & \pi_0\gcl\tiund{\pi_0\hcl}\pi_0\gcl'=G_0/\pt_1(\ker\pt_0)\tiund{H_0}G'_0/\pt'_1(\ker\pt'_0)
    \end{split}
\end{equation*}
The isomorphism
$\pi_0\gcl\tiund{\pi_0\hcl}\pi_0\gcl'\rw\pi_0(\gcl\tiund{\hcl}\gcl')$ takes
$([x],[y])$ to $([(x,y)])$ for each
$([x],[y])\in\pi_0\gcl\tiund{\pi_0\hcl}\pi_0\gcl'$. The rest is immediate from
part (ii).\enpr

\subsection{Simplicial directions in a cat$^n$-group}\label{simpdir}
A basic fact about $\ctn$s, which will be used throughout Sections
\ref{catcat}, \ref{scontra}, \ref{specia} is that a $\ctn$ can be considered as
an internal category in $\ctmen$s in $n$ possible ways, corresponding to the
$n$ simplicial directions of its multinerve (Proposition \ref{simpdir.pro1}).
To show this, we use the following Lemma \ref{simpdir.lem1} (and its Corollary
\ref{simpdir.cor1}) which are a straightforward consequence of the definition
of $\ctn$s. Proposition \ref{simpdir.pro1} will be used in the proof of Theorem
\ref{adj.the1} as well as in the proof of Lemma \ref{spec.lem3}. The $\ctmen$s
$\gcl^{(k)}_i$ associated to a $\ctn$ $\gcl$ as in Corollary \ref{simpdir.cor1}
will be used throughout Sections \ref{catcat}, \ref{scontra}, \ref{specia}.
\begin{lemma}\label{simpdir.lem1}
    An $n$-simplicial group $\psi\in[\dnop,\gp]$ is in the image of the
    multinerve functor $\ncl:\cagp\rw[\dnop,\gp]$ if and only if, for all $1\leq r\leq
    n$ and $[p_1]\ldots [p_r]\in\dop$, $p_r\geq 2$, there is an isomorphism in
    $[\Delta^{n-r^{op}},\gp]$
    \begin{equation}\label{simpdir.eq1}
        \psi(p_1... p_r\,\mi)\cong\pro{\psi(p_1... p_{r-1}\,1\,\mi)}
        {\psi(p_1... p_{r-1}\,0\,\mi)}{p_r}
    \end{equation}
\end{lemma}
\prf We use induction on $n$. By Proposition \ref{simtec.pro1} the lemma is
true for $n=1$. Suppose it holds for $n-1$ and let $\gcl\in\cat(\cagpm)$ with
object of objects and object of arrows $\gcl_0\;\gcl_1\in\cagpm$ respectively.
For each $p\geq 2$ put $\gcl_p=\gistr{p}$. Then by the definition of the
multinerve, for each $[p_1]\ldots [p_r]\in\dop$, $1\leq r\leq n$,
$\ncl\gcl(p_1... p_r\,\mi)=\ncl\gcl_{p_1}(p_2... p_r\,\mi)$. Hence, using the
induction hypothesis, we obtain
\begin{equation*}
    \begin{split}
        & \ncl\gcl(p_1...p_r\,\mi)=\ncl\gcl_{p_1}(p_2... p_r\,\mi)\cong \\
         & \cong \pro{\ncl\gcl_{p_1}(p_2... p_{r-1}\,1\,\mi)}
         {\ncl\gcl_{p_1}(p_2... p_{r-1}\,0\,\mi)}{p_r}=\\
          & =\pro{\ncl\gcl(p_1... p_{r-1}\,1\,\mi)}
         {\ncl\gcl(p_1... p_{r-1}\,0\,\mi)}{p_r}.
     \end{split}
\end{equation*}
Conversely, suppose that $\psi\in[\dnop,\gp]$ satisfies (\ref{simpdir.eq1}).
Regarding $\psi$ as a simplicial object in $(\nm)$-simplicial groups along
direction 1, we obtain the following simplicial object:
\begin{equation}\label{simpdir.eq2}
    \cdots\;\pile{\rTo\\ \rTo\\ \rTo}\psi(2\;\mi)\pile{\rTo\\ \rTo\\
    \rTo}\psi(1\;\mi)\pile{\rTo\\ \rTo\\ \lTo}\psi(0\;\mi).
\end{equation}
For each $[p]\in\dop$, $\psi(p\;\mi)$ is a $(\nm)$-simplicial group satisfying
(\ref{simpdir.eq1}), hence by the induction hypothesis it is the multinerve of
a $\ctmen$; that is, $\psi(p\;\mi)=\ncl\gcl_p$ for $\gcl_p\in\cagpm$. Also,
since $\psi$ satisfies (\ref{simpdir.eq1}), in particular taking $r=1$, we have
$\psi(p\;\mi)\cong\pro{\psi(1\;\mi)}{\psi(0\;\mi)}{p}$ for $p\geq 2$. Hence
(\ref{simpdir.eq2}) is isomorphic to
\begin{equation*}
    \cdots\;\pile{\rTo\\ \rTo\\ \rTo}\tens{\ncl\gcl_1}{\ncl\gcl_0}
    \pile{\rTo\\ \rTo\\ \rTo}\ncl\gcl_1\pile{\rTo\\ \rTo\\ \lTo}\ncl\gcl_0.
\end{equation*}
Since $\pro{\ncl\gcl_1}{\ncl\gcl_0}{}=\ncl(\gistr{})$ (as $\ncl$, being right
adjoint, preserves limits), we finally have a $\ctn$ $\gcl$
\begin{equation*}
    \tens{\gcl_1}{\gcl_0}\rw\gcl_1\pile{\rTo\\ \rTo\\ \lTo}\gcl_0
\end{equation*}
with $\ncl\gcl=\psi$, as required.\enpr
\begin{corollary}\label{simpdir.cor1}
    Let $\gcl\in\cagp$. Then for each $1\leq k\leq n$ and $[i]\!\in\!\!\dop$ there
    exists a cat$^{n-1}\!$-group $\gcl_i^{(k)}\!\!$ with $\ncl\gcl_i^{(k)}(p_1...
    p_{\nm})\!=\!\ncl\gcl(p_1... p_{k-1}\,i\,p_k... p_\nm)$.
\end{corollary}
\prf Let $1\leq r\leq n$ $r\neq k$. By Lemma \ref{simpdir.lem1} there is an
isomorphism for $p_r\geq 2$
\begin{equation*}
    \begin{split}
       & \ncl\gcl(p_1... p_{r}... p_n)\cong \\
        &\! \cong\!\pro{\ncl\gcl(p_1... p_{r-1} 1... p_n)}
    {\ncl\gcl(p_1... p_{r-1} 0... p_n)}{p_r}.
    \end{split}
\end{equation*}
In particular, evaluating this isomorphism at $p_k=i$, this is saying that the
$(n\!-\!1)$-simplicial group taking $(p_1... p_\nm)$ to $\ncl\gcl(p_1...
p_{k-1}\, i... p_\nm)$\! satisfies condition\! (\ref{simpdir.eq1}) in Lemma
\ref{simpdir.lem1}. Hence by \ref{simpdir.lem1} there exists a $\ctmen$
$\gcl_i^{(k)}$ with $\ncl\gcl_i^{(k)}(p_1... p_\nm)\!=\ncl\gcl(p_1...
p_{k-1}\;i\;p_k... p_\nm)$, as required.\enpr
\begin{proposition}\label{simpdir.pro1}
    For each $1\!\leq\! k\!\leq n$ there is an isomorphism $\xi_k:\!\cagp\!\rw\cagp$ which
    associates to each cat$^n$-group an object of\, $\cat(\cagpm)$ $\xi_k\gcl$ with
    object of objects and object of arrows given by
    $(\xi_k\gcl_0)_0\!=\!\gcl_0^{(k)}$, $(\xi_k\gcl_1)_1\!=\!\gcl_1^{(k)}$  where
    $\gcl_i^{(k)}$ is as in Corollary \ref{simpdir.cor1}.
\end{proposition}
\prf From Section \ref{simtec} the $n$-simplicial group $\ncl\gcl$ can be
regarded as a simplicial object in $(\nm)$-simplicial groups along the $k^{th}$
simplicial direction, taking each $[i]\in\dop$ to the object of $[\dnomen,\gp]$
which associates to $(p_1... p_\nm)\in\dnomen$ the group $\ncl\gcl(p_1...
p_{k-1}\;i\;p_{k}... p_\nm)$. By Corollary \ref{simpdir.cor1} the latter is the
multinerve of a $\ctmen$ $\gcl_i^{(k)}$. Further, from Lemma \ref{simpdir.lem1}
for each $i\geq 2$ we have
$\ncl\gcl_i^{(k)}\cong\pro{\ncl\gcl_1^{(k)}}{\ncl\gcl_0^{(k)}}{i}$. Hence
$\ncl\gcl$, as a simplicial object in $[\dnomen,\gp]$ along the $k^{th}$
direction, has the form
\begin{equation*}
    \cdots\ncl(\tens{\gcl_1^{(k)}}{\gcl_0^{(k)}})
    \pile{\rTo\\ \rTo\\ \rTo}\ncl\gcl_1^{(k)}\pile{\rTo\\ \rTo\\ \lTo}\ncl\gcl_0^{(k)}.
\end{equation*}
We define $\xi_k\gcl$ to be the $\ctn$
\begin{equation*}
    \tens{\gcl_1^{(k)}}{\gcl_0^{(k)}}\rw\gcl_1^{(k)}\pile{\rTo\\ \rTo\\ \lTo}\gcl_0^{(k)}.
\end{equation*}
We are now going to define the inverse for $\xi_k$. Let $\gcl\in\cat(\cagpm)$,
with object of objects and object of arrows given by $\gcl_0$ and $\gcl_1$
respectively. For each $i\geq 2$ put $\gcl_i=\gistr{i}$. Consider the
$n$-simplicial group $\psi_k$ taking $(p_1\ldots p_n)$ to $\psi_k(p_1\ldots
p_n)=\ncl\gcl_{p_k}(p_1\ldots p_{k-1}p_{k+1}\ldots p_n)$. It is easy to see
that $\psi_k$ satisfies (\ref{simpdir.eq1}) in Lemma \ref{simpdir.lem1}. In
fact, for $r\neq k$ this follows from the fact that $\ncl\gcl_{p_k}$ satisfies
(\ref{simpdir.eq1}). For $r=k$, notice that when $p_k\geq 2$ we have
\begin{equation*}
     \begin{split}
        & \psi_k(p_1\ldots p_n)=\ncl(\gistr{p_k})(p_1\ldots p_{k-1}p_{k+1}\ldots p_n)= \\
         & =\pro{\psi_k(p_1... p_{k-1}\,1... p_n)}{\psi_k(p_1... p_{k-1}\,0... p_n)}
         {p_k}.
     \end{split}
\end{equation*}
Hence by Lemma \ref{simpdir.lem1} there exists $\xi'_k\gcl\in\cagp$ such that
$\ncl(\xi'_k\gcl)=\psi_k$. It is immediate to check that $\xi_k$ and $\xi'_k$
are inverse bijections.\enpr

\subsection{Internal categories in categories of groups with operations}\label{oper}
The goal of this section is to prove Theorem \ref{oper.the1} asserting the
existence of a certain adjunction between $\cat\dcl$ and $\dcl$ for a class of
categories of groups with operations $\dcl$. Theorem \ref{oper.the1} is one of
the main tools used in the proof of Theorem \ref{adj.the0} where it will be
applied to the case where $\dcl$ is the category of groups with operations
$\cngm$, described in the next section. This will yield an adjunction between
$\cng$ and $\set$ (Theorem \ref{adj.the0}) and then a corresponding adjunction
between $\cagp$ and $\set$ (Theorem \ref{adj.the1}). We first need to recall
some basic facts about internal categories in categories of groups with
operations. With the exception of Theorem \ref{oper.the1}, the material in this
section is well known.

There exists in the literature a well known notion of category of groups with
operations, see for instance \cite{orz}, \cite{por1}. These are particular
types of universal algebras defined by sets $\Omega_i$ of operations of order
$i=0,1,2$ satisfying certain axioms. For simplicity we will restrict our
attention to categories of groups with operations where the only operation of
order 2 is the group multiplication. This will be sufficient for the
application of this notion to $\ctn$s which we pursue in Sections \ref{cgn},
\ref{adj}, \ref{twolem}.

Thus, following \cite{por1}, a \emph{category of groups with operations} is for
us a variety of universal algebras which is first of all a variety of groups.
Also, group identity is the only operation of order 0 and group multiplication
is the only operation of order 2. Denoting by $\Omega'_1$ the set of operations
of order 1 different from group inverse, we require
$\omega(a\;b)=\omega(a)\omega(b)$ for all $\omega\in\Omega'_1$, $a,b\in\dcl$. A
morphism in a category of groups with operations of the type above consists of
a group homomorphism which commutes with all $\omega\in\Omega'_1$.

In a category $\dcl$ of groups with operations we have notions of action and of
semidirect product which generalize well known notions in the category of
groups. Given $A,B\in \ob \dcl$, a \emph{B-structure A} is a split extension of
$B$ by $A$ in $\dcl$:
\begin{equation}\label{oper.eq1}
    A\;\overset{i}{\rightarrowtail}\;E\;\;\pile{\overset{p}{\twoheadrightarrow}\\
    \underset{s}{\leftarrow}}\;\; B\qquad\qquad ps=\id_B.
\end{equation}
If $A$ is a $B$-structure we define the semidirect product $A\ttil B$ as the
object of $\dcl$ which is $A\times B$ as a set, with operations
\begin{equation*}
    \begin{split}
        & \omega(a,b)=(\omega(a),\omega(b)) \\
         & (a',b')(a,b)=(a's(b')as(b')^{-1},b'b)
     \end{split}
\end{equation*}
for all $a,a'\in A$, $b,b'\in B$, $\omega\in\Omega'_1$. Given the split
extension (\ref{oper.eq1}), there is an isomorphism $A\ttil B\cong E$ taking
$(a,b)$ to $i(a)s(b)$.

Internal categories in categories of groups with operations are particularly
easy to describe. In fact, given $\abb\in\cat\dcl$:
\begin{equation*}
    \tens{A_1}{A_0}\rTo^m A_1\;\pile{\rTo^{\pt_0}\\ \rTo_{\pt_1}\\
     \lTo_{\sigma_0}}\;A_0,
\end{equation*}
the composition $m$ is uniquely determined by
\begin{equation}\label{oper.eq2}
    m(x,y)=x(\sigma_0\pt_1 x^{-1})y
\end{equation}
for $(x,y)\in \tens{A_1}{A_0}$. This follows easily from the identity
\begin{equation*}
    (x,y)=(x(\sigma_0\pt_1 x^{-1}),1)(\sigma_0\pt_1 x,\sigma_0\pt_0
    y)(1,(\sigma_0\pt_0 y^{-1})y)
\end{equation*}
and the fact that $m$ is a group homomorphism. The identity (\ref{oper.eq2})
has two immediate consequences. One is that every arrow $z$ in $A_1$ has an
inverse for composition given by $(\sigma_0\pt_0 z)z^{-1}(\sigma_0\pt_1 z)$, so
every internal category in $\dcl$ is an internal groupoid. The second
consequence of (\ref{oper.eq2}) is the identity
\begin{equation}\label{oper.eq3}
    [\ker\pt_0,\ker\pt_1]=1.
\end{equation}
In fact, if $x\in\ker\pt_0$, $y\in\ker\pt_1$, since $m$ is a group homomorphism
we have:
\begin{equation*}
\begin{split}
 y  & =m(\sigma_0\pt_1 y,y)=m(1,y)=m((x,1)(1,y)(x^{-1},1))= \\
    & =m(x,\sigma_0\pt_0 x)m(1,y)m(x^{-1},\sigma_0\pt_0 x^{-1})=xyx^{-1}
\end{split}
\end{equation*}
which proves (\ref{oper.eq3}). Thus $\abb$ gives rise to a reflexive graph
\begin{equation*}
    A_1\;\pile{\rTo^{\pt_0}\\ \rTo_{\pt_1}\\
     \lTo_{\sigma_0}}\;A_0
\end{equation*}
satisfying condition (\ref{oper.eq3}). Conversely, given such a reflexive graph
satisfying (\ref{oper.eq3}), we can build an object of $\cat\dcl$ by putting
$m(x,y)=x(\sigma_0\pt_1 x^{-1})y$. Using (\ref{oper.eq3}) it easy to check that
$m$ is a group homomorphism; it is also clear that $m$ commutes with every
$\omega\in\Omega'_1$, as the same is true for $\sigma_0,\pt_1$ and
$\omega(a,b)=\omega(a)\,\omega(b)$. Thus $m$ is a morphism in $\dcl$. The above
discussion motivates the following definition.
\begin{definition}\label{oper.def1}
    The category $\ccl^1\dcl$ has for objects the triples $(D,d,t)$ where $D$ is
    an object of $\dcl$ and $d,t:D\rw D$ are morphisms in $\dcl$ such that
    $dt=t$, $td=d$ and $[\ker d, \ker t]=1$.
\end{definition}
The following basic fact is well known, and will be used throughout this and
the subsequent sections.
\begin{lemma}\label{oper.lem1}
    The categories $\cat\dcl$ and $\ccl^1\dcl$ are equivalent.
\end{lemma}
\prf Given $\abb\in\cat\dcl$, from the above discussion  $[\ker \pt_0, \ker
\pt_1]=1$. Since $\sigma_0$ is injective, it follows that $[\ker \sigma_0\pt_0,
\ker \sigma_0\pt_1]=1$. Hence $(A_1, \sigma_0\pt_0,\sigma_0\pt_1)$ is an object
of $\ccl^1\dcl$. Conversely, given an object $(D,d,t)$ of $\ccl^1\dcl$, we have
a reflexive graph
\begin{equation*}
    D\;\pile{\rTo^{d}\\ \rTo_{t}\\
     \lTo_{i}}\;\im d
\end{equation*}
where $i$ is the inclusion, and $[\ker d, \ker t]=1$. From the discussion
before Definition \ref{oper.def1}, this determines an object of $\cat\dcl$. It
is easily seen that this gives an equivalence of categories.\enpr
\begin{remark}\rm\label{oper.rem1}
There is a well known notion of crossed module in any category of groups with
operations, see for instance \cite{por1}. Although we will not use crossed
modules in this paper, we recall that the equivalence of Lemma \ref{oper.lem1}
extends to an equivalence of categories between $\cat \dcl$ and crossed modules
in $\dcl$. More details can be found for instance in \cite{por1}.
\end{remark}
The next two lemmas are straightforward facts which will be used in the proof
of Theorems \ref{oper.the1} and \ref{adj.the1}.
\begin{lemma}\label{oper.lem2}
    Let $f:A\rw A'$ be a morphism in $\dcl$.
\begin{itemize}
  \item[a)] Then the following is an object of
    $\cat\dcl$, which we denote by $A_f$:
    \begin{equation}\label{oper.eq4}
        \tens{(\ker f\ttil A)}{A}\rTo^{m}\ker f\ttil A\;\pile{\rTo^{\pt_0}\\
        \rTo_{\pt_1}\\ \lTo_{\sigma_0}}\;A
    \end{equation}
where $\ker f\ttil A$ is the object of $\cat\dcl$ which is $\ker f\times A$ as
a set, with operations
\begin{equation*}
    \begin{split}
        & (x',y')(x,y)=(x'y'xy'^{-1},y'y) \\
         & \omega(x',y')=(\omega x',\omega y')
     \end{split}
\end{equation*}
for all $\omega \in \Omega'_1$, $(x,y),(x',y')\in\ker f\times A$. Also,
$\pt_0(x,y)=y$, $\;\pt_1(x,y)=xy$, $\;\sigma_0(x)=(1,x)$,
$\;c((x,y)(x',xy))=(x'x,y)$.

\noindent Given morphisms in $\dcl$, $f:A\rw A'$, $g:B\rw B'$, $r:A\rw B$,
$s:A'\rw B'$ such that $sf=gr$, then the maps $\ker f \ttil
A\supar{(r_{|},r)}\ker g\ttil B$ and $r$ determine a morphism $A_f\rw B_g$ in
$\cat\dcl$.
  \item[b)] $A_f$ is isomorphic to the following object of $\cat\dcl$, which we
  denote by $A^f$:
  \begin{equation*}
    \tens{(\tens{A}{f})}{A}\rTo^{c}\tens{A}{f}\pile{\rTo^{p_1}\\ \rTo_{p_2}\\
    \lTo_{s}} A
  \end{equation*}
where $\tens{A}{f}$ is the kernel pair of $f$, $p_1(x,y)=x$, $p_2(x,y)=y$,
$s(z)=(z,z)$, $c((x,y)(y,z))=(x,z)$.
\end{itemize}
\end{lemma}
 \prf It is immediate to verify that $\pt_0,\pt_1,\sigma_0, m$ are
morphisms in $\dcl$. Part a) follows from a straightforward verification of the
axioms of internal categories and internal functors. For part b), consider the
map $\alpha:\tens{A}{f}\rw\ker f\ttil A$ given by $\alpha(x,y)=(yx^{-1},x)$. It
is straightforward to verify that $\alpha$ is a morphism in $\dcl$. Clearly
$\alpha$ is injective. Given $(z,y)\in\ker f\ttil A$ we have
$(z,y)=\alpha(y,zy)$, so $\alpha$ is also surjective, hence an isomorphism. It
is straightforward to verify that $A^f$ is an internal category and that
$(\alpha,\id)$ is an internal functor. \enpr
\begin{lemma}\label{oper.lem3}
    Let $A,B$ be two objects of $\dcl$ and let $u_1:A\rw A\amalg B$ and $u_2:B\rw A\amalg B$
     be the two coproduct
    injections. Then every element of $A\amalg B$ has the form
    $u_1(a_1)u_2(b_1)u_1(a_2)u_2(b_2)...$ or
    $u_2(b_1)u_1(a_1)u_2(b_2)u_1(a_2)...$ where $a_i\in A$, $b_i\in B$.
\end{lemma}
\prf See Appendix.

We now come to the main result of Section \ref{oper}. The following theorem
will be applied in Section \ref{adj} to the proof of Theorem \ref{adj.the0} in
the case where $\dcl$ is the category of groups with operations $\cngm$,
described in the next Section \ref{cgn}. In the case where $\dcl$ is the
category of groups, Theorem \ref{oper.the1} corresponds to the result of
\cite[Proposition 2]{ccg} which is formulated using the language of crossed
modules in groups rather than the equivalent category $\cgp$.
\begin{theorem}\label{oper.the1}
    Let $\vcl:\cat\dcl\rw\dcl$ be given by $\vcl(\abb)=\ker\pt_0\times A_0$.
    Given $H\in\dcl$, let $p:H\amalg H\rw H$ be determined by $pu_1=1$
    $pu_2=\id$ where $u_1, u_2:H\rw H\amalg H$ are the two coproduct
    injections. Let $(H\amalg H)_p\in\cat\dcl$ be as in lemma \ref{oper.lem2} a).
    Then the functor
    $\lcl:\dcl\rw\cat\dcl$ $\lcl(H)=(H\amalg H)_p$ is left adjoint to $\vcl$.
\end{theorem}
\prf See Appendix.

\subsection{An equivalent description of cat$^n$-groups}\label{cgn}
In this section we introduce the category $\cng$ and we show it is equivalent
to the category of $\ctn$s. The material in this section is well known.

Historically, $\cng$ was the first category to be described in \cite{lod} and
called the category of $\ctn$s. In this paper we mainly work with the
equivalent category $\cagp$ so that we refer to $\cagp$ as the category of
$\ctn$s. However, we need the category $\cng$ in Section \ref{adj} in order to
describe an important adjunction between $\ctn$s and sets.

We observed in Lemma \ref{oper.lem1} that for any category $\dcl$ of groups
with operations $\cat\dcl$ is equivalent to $\ccl^1\dcl$. If $\dcl$ is the
category of groups, we shall denote $\ccl^1\gp$ by $\ccl^1\gcl$. The idea of
the category $\cng$ when $n>1$ is to have $n$ internal categorical structures
in groups which are mutually independent.
\begin{definition}\label{cgn.def1}
The category $\cng$ is defined as follows. An object of $\cng$ consists of a
group $G$ together with $2n$ endomorphisms $t_i,d_i:G\rw G, \;\;1\leq i \leq
n$, such that, for all $1 \leq i,j \leq n$,
\begin{itemize}
  \item [(i)] $d_i t_i=t_i,\quad t_i d_i=d_i,$
  \item [(ii)] $d_it_j=t_jd_i,\quad d_id_j=d_jd_i,\quad t_it_j=t_jt_i,\qquad i \neq j,$
  \item [(iii)] $[\ker d_i, \ker t_i]=1$.
\end{itemize}
A morphism  $(G,d_i,t_i)\rw (G',d'_i,t'_i)$ in $\cng $ is a group homomorphism
$f:G\rw G'$ such that $f d_i=d'_if$, $\,ft_i=t'_if$ $\,1\leq i \leq n$.
\end{definition}
In the above definition, axioms (i) and (iii) say that for each fixed $i$,
$(G,d_i,t_i)$ is an object of $\ccl^1\gcl$, while (ii) says that for $i\neq j$
$(G,d_i,t_i)$ and $(G,d_j,t_j)$ are independent. Notice that the identity (iii)
in Definition \ref{cgn.def1} is equivalent to
\begin{itemize}
  \item[(iii)'] $d_i(x)x^{-1}t_i(x)x^{-1}=t_i(x)x^{-1}d_i(x)x^{-1}$, $\quad x\in
  G$.
\end{itemize}
Using this fact we can view $\cng$ as a category of groups with operations as
follows. The only operations of order 0 and 2 are, respectively, group identity
and multiplication. The set $\Omega'_1$ of operations of order 1 different from
group inverse consists of $2n$ operations $d_i,t_i:G\rw G,\;\,\iun$; there is
also a set of identities which comprises the group laws and the identities (i)
(ii) and (iii)' above.

Since $\cng$ is a category of groups with operations, from Section \ref{oper}
we can construct the category $\ccl^1(\cng)$. The following lemma is a
straightforward consequence of the definitions.
\begin{lemma}\label{cgn.lem1}
    The categories $\ccl^1(\cng)$ and $\cngp$ are isomorphic.
\end{lemma}
\prf By definition, an object of $\ccl^1(\cng)$ consists of a triple
$(\gcl,d,t)$ where $\gcl=(G,d_1,\ldots,d_n,t_1,\ldots,t_n)$ is an object of
$\cng$ and $d,t:\gcl\rw\gcl$ are morphisms in $\cng$ satisfying $[\ker d, \ker
t]=1$, $dt=t$, $td=d$. Fix $1\leq k\leq \np$; since $d,t$ are morphisms in
$\cng$, they commute with all $d_i,t_i$, $\iun$. It follows that
$\gcl'=(G,d_1,\ldots, d_{k-1}, \,d,\,d_{k},\ldots, d_n,\,t_1,\ldots, t_{k-1},\,
t,\,t_{k},\ldots, t_n)$ is an object of $\cngp$. Conversely, let
$\gcl=(G,d_1,\ldots,d_\np, t_1,\ldots,t_\np)\in\cngp$. Then if we leave out the
operators $d_k,t_k$, we find an object of $\cng\,$ $\gcl'=(G,d_1,\ldots,
d_{k-1}, \,d_{k+1},$ $\ldots, d_n,\,t_1,\ldots, t_{k-1},\,t_{k+1},\ldots
,t_n)$. Further, since $d_kt_k=t_k$, $t_kd_k=d_k$ and $[\ker d_k,$ $\ker
t_k]=1$, it follows that $(\gcl',d_k,t_k)$ is an object of $\ccl^1\cng$. It is
trivial to check that these correspondences give an isomorphism of categories.
\enpr

\bigskip

The following corollaries will be needed in the proof of Theorems
\ref{adj.the0} and \ref{adj.the1}.
\begin{corollary}\label{cgn.cor0}
For each $1\leq k\leq n$, there is an equivalence of categories
$\alpha^{(k)}_n\dashv\beta^{(k)}_n$:
\begin{equation*}
   \alpha^{(k)}_n:\cat\cngm\pile{\lTo\\ \rTo}\cng:\beta^{(k)}_n.
\end{equation*}
\end{corollary}
\prf By Lemma \ref{cgn.lem1} for each $1\leq k\leq n$ there is an isomorphism
$\cng\cong \ccl^1(\cngm)$; by Lemma \ref{oper.lem1} there is an equivalence
$\ccl^1(\cngm)\simeq\cat(\cngm)$. Hence the result.\enpr
\begin{corollary}\label{cng.cor1}
    The categories $\cng$ and $\cagp$ are equivalent.
\end{corollary}
\prf By induction on $n$. The case $n=1$ is an instance of Lemma
\ref{oper.lem1}, taking for $\dcl$ the category of groups. Suppose,
inductively, that $\cngm$ and $\cagpm$ are equivalent. By Lemma \ref{cgn.lem1},
$\cng$ is isomorphic to $\ccl^1\cngm$ and, by Lemma \ref{oper.lem1} this is
equivalent to $\cat(\cngm)$. From the induction hypothesis, $\cat(\cngm)$ is
equivalent to $\cat(\cagpm)=\cagp$. Thus $\cng$ and $\cagp$ are
equivalent.\enpr

\subsection{An adjunction between cat$^n$-groups and Set}\label{adj}
In this section we discuss an important adjunction $\fl_n\dashv\ul_n$ between
the category $\cagp$ and $\set$. This adjunction will play a crucial role in
the construction of the functor $\spi$ in Section \ref{specia}. The adjunction
$\fl_n\dashv\ul_n$ is deduced from an analogous adjunction $\fcl_n\dashv\ucl_n$
between $\cng$ and $\set$ as in the following theorem.
\begin{theorem}\label{adj.the0}
    Let $1\leq k\leq n$ and let $U:\gp\rw\set$ be the forgetful functor
    and $F:\set\rw\gp$ its left adjoint. There is a functor
    $\ucl_n:\cng\rw\set$ whose left adjoint $\fcl_n$ is given by $\fcl_0=F$
    and, for each $n\geq 1$ and any set $X$,
    \begin{equation*}
        \fcl_n(X)=\alpha^{(k)}_n(\fcl_\nm(X)\amalg\fcl_\nm(X))_p
    \end{equation*}
where $(\fcl_\nm(X)\amalg\fcl_\nm(X))_p$ is as in Lemma \ref{oper.lem2} a),
$p:\fcl_\nm(X)\amalg\fcl_\nm(X)\!\!\rw\fcl_\nm(X)$ is defined by $pu_1=1$,
$pu_2=\id$, $u_1,u_2:\fcl_\nm(X)\rw\fcl_\nm(X)\amalg\fcl_\nm(X)$ are the
coproduct injections and $\alpha^{(k)}_n$ is the equivalence of categories
$\alpha^{(k)}_n:\cat(\cngm)\rw\cng$ as in Corollary \ref{cgn.cor0}.
\end{theorem}
\prf We define $\fcl_n$ inductively. Recall from Theorem \ref{oper.the1} that
there is a functor $\vcl_1:\cat(\gp)\rw\gp$ with a left adjoint $\lcl_1$ given
by $\lcl_1 H=(H\amalg H)_p$ where $p:H\amalg H\rw H$ is defined by $pu_1=1$,
$pu_2=\id$ and $(H\amalg H)_p$ is as in Lemma \ref{oper.lem2}. Thus for $n=1$
we define $\ucl_1$ to be the composite
\begin{equation*}
    \ccl^1\gcl\supar{\beta_1}\cat(\gp)\supar{\vcl_1}\gp\supar{U}\set,
\end{equation*}
where $\beta_1$ is as in Corollary \ref{cgn.cor0}. Therefore $\ucl_1$ has a
left adjoint $\fcl_1=\alpha_1\lcl_1 F$ and $\fcl_1(X)=\alpha_1(F(X)\amalg
F(X))_p$.

Inductively, suppose we have defined $\ucl_\nm:\cngm\rw\set$ with left adjoint
$\lcl_\nm$. Applying Theorem \ref{oper.the1} to the case where $\dcl$ is the
category of groups with operations $\cngm$ we obtain a functor
$\vcl_n:\cat(\cngm)\rw\cngm$ with a left adjoint
$\lcl_n\hcl=(\hcl\amalg\hcl)_p$ where $p:\hcl\amalg\hcl\rw\hcl$, $pu_1=1$,
$pu_2=\id$. We define $\ucl_n$ to be the composite
\begin{equation*}
    \ucl_n:\cng\supar{\beta^{(k)}_n}\cat(\cngm)\supar{\vcl_n}\cngm\supar{\ucl_\nm}\set,
\end{equation*}
where $\beta^{(k)}_n$ is as in Corollary \ref{cgn.cor0}. Then $\ucl_n$ has left
adjoint $\fcl_n=\alpha^{(k)}_n\lcl_n\fcl_\nm$; that is,
\begin{equation*}
    \fcl_n(X)=\alpha^{(k)}_n(\fcl_\nm(X)\amalg\fcl_\nm(X))_p.
\end{equation*}
\enpr

\bigskip

The following lemma is a standard fact about internal categories which
generalizes Lemma \ref{oper.lem2} b) from the case of a category of groups with
operations to the case of an arbitrary category $\ccl$ with finite limits. This
will be needed in the next theorem where it will be applied to the case where
$\ccl=\cagpm$.
\begin{lemma}\label{adj.lem2}
    Let $f:A\rw B$ be a morphism in a category $\ccl$ with finite limits; the
    identity map $\id_A$ defines a unique map $s:A\rw\tens{A}{f}$ such that $p_1 s=\id_A=p_2
    s$ where $p_1,p_2:\tens{A}{f}\rw A$ are the two projections. The commutative diagram
    \begin{diagram}[h=2.6em]
        \tens{A}{f} & \rTo^{p_1}& A & \lTo^{p_2} & \tens{A}{f}\\
        \dTo^{p_2}&& \dTo^{f} && \dTo_{p_1}\\
        A & \rTo_{f}& B & \lTo_f & A
    \end{diagram}
    defines a unique morphism $m:\tens{(\tens{A}{f})}{A}\rw\tens{A}{f}$
    such that $p_2 m=p_2\pi_2,\;p_1 m=p_1\pi_1$ where $\pi_1,\pi_2:
    \tens{(\tens{A}{f})}{A}\rw\tens{A}{f}$ are the two projections. Then the
    following is an object of $\cat\ccl$ which we denote by $A^f$:
    \begin{equation*}
       A^f : \tens{(\tens{A}{f})}{A}\rTo^m\tens{A}{f}\;\pile{\rTo^{p_1}\\ \rTo^{p_2}\\
        \lTo_{s}}\;A.
    \end{equation*}
\end{lemma}
\prf It is a straightforward verification of the axioms of internal
categories.\enpr

\bigskip

In the following theorem the map $1:\gcl\rw\gcl$ where $\gcl\in\cagp$ is as
follows. Since $\ncl\gcl\in[\dnop,\gp]$ there is a map $e:\ncl\gcl\rw\ncl\gcl$
whose component $e_{p_1\ldots p_n}:\ncl\gcl(p_1\ldots p_n)\rw\ncl\gcl(p_1\ldots
p_n)$ for $({[p_1]}\ldots{[p_n]})\in\dnop$ is given by $e_{p_1\ldots p_n}(x)=1$
for all $x\in\ncl\gcl({p_1}\ldots{p_n})$. Since $\ncl$ is fully faithful there
is a unique map $1:\gcl\rw\gcl$ such that $\ncl(1)=e$.

Given $\gcl\in\cagp$ and $1\leq k \leq n$, in what follows $\gcl_1^{(k)},
\gcl_0^{(k)}\in\cagpm$ are as in Corollary \ref{simpdir.cor1}.
\begin{theorem}\label{adj.the1}
    There is a functor $\ul_n:\cagp\rw\set$ whose left adjoint $\fl_n$ is given
    by $\fl_0=F$ and, for each $n\geq 1$, $1\leq k\leq n$ and any set $X$, as an
    internal category in $\cagpm$ in the $k^{th}$ direction, $\fl_n(X)$ is isomorphic to
    \begin{equation*}
        (\fl_\nm(X)\amalg\fl_\nm(X))^p
    \end{equation*}
as in Lemma \ref{adj.lem2} where $p:\fl_\nm(X)\amalg\fl_\nm(X)\rw\fl_\nm(X)$ is
defined by $pu_1=1$, $p u_2=\id$ and
$u_1,u_2:\fl_\nm(X)\rw\fl_\nm(X)\amalg\fl_\nm(X)$ are the coproduct injections.
\end{theorem}
\prf Denote by $\gamma_n$ the equivalence of categories $\gamma_n:\cng\rw\cagp$
of Corollary \ref{cng.cor1} and by $\eta_n$ its right adjoint pseudo-inverse.
Let $\ul_n=\ucl_n\eta_n:\cagp\rw\set$. Then, by Theorem \ref{adj.the0}, $\ul_n$
has a left adjoint $\fl_n=\gamma_n\fcl_n$. Recalling from the proof of Theorem
\ref{adj.the0} that $\ucl_n=\ucl_\nm\vcl_n\beta^{(k)}_n$ and
$\fcl_n=\alpha^{(k)}_n\lcl_n\fcl_\nm$, we have a commutative diagram:
\begin{diagram}
    \cagp & \pile{\lTo^{\gamma_n}\\ \rTo_{\eta_n}}& \cng & \pile{\lTo^{\fcl_n}\\
    \rTo_{\ucl_n}}& \set\\
    \dTo^{\xi_k}\uTo_{\xi'_k}&& \dTo^{\beta^{(k)}_n}\uTo_{\alpha^{(k)}_n} &&
    \uTo^{\ucl_\nm}\dTo_{\fcl_\nm}\\
    \cat(\cagpm) & \pile{\lTo^{\cat\gamma_\nm}\\ \rTo_{\cat\eta_\nm}}& \cat(\cngm) & \pile{\lTo^{\lcl_n}\\
    \rTo_{\vcl_n}}& \cngm
\end{diagram}
where $\xi_k$ and $\xi'_k$ are as in Proposition  \ref{simpdir.pro1} and
$\alpha^{(k)}_n,\;\beta^{(k)}_n$ as in Corollary \ref{cgn.cor0}. Then for any
set $X$ we have
\begin{equation}\label{adj.eq2}
\begin{split}
   &  \xi_k\fl_n(X)=\xi_k\gamma_n\fcl_n(X)=\\
    & (\cat\gamma_\nm)\beta^{(k)}_n\alpha^{(k)}_n\lcl_n\fcl_\nm(X)\cong
    (\cat\gamma_\nm)\lcl_n\fcl_\nm(X).
\end{split}
\end{equation}
From Theorem \ref{oper.the1},
$\lcl_n\fcl_\nm(X)=(\fcl_\nm(X)\amalg\fcl_\nm(X))_p$; and by Lemma
\ref{oper.lem2} b) the latter is isomorphic to
$(\fcl_\nm(X)\amalg\fcl_\nm(X))^p$ so that
\begin{equation}\label{adj.eq2bis}
    \begin{split}
       & (\lcl_n\fcl_\nm(X))_0\cong\fcl_\nm(X)\amalg\fcl_\nm(X)\;\,\text{and} \\
        & (\lcl_n\fcl_\nm(X))_1\cong\tens{(\fcl_\nm(X)\amalg \fcl_\nm(X))}{p}.
    \end{split}
\end{equation}
Also notice that for each $n$, $\gamma_n$ preserves kernel pairs. In fact, it
is immediate to check that for each $n$, $\beta^{(k)}_n$ preserves kernel
pairs; since $\gamma_n=\xi_k\cat\gamma_\nm\beta^{(k)}_n$ a simple inductive
argument shows that $\gamma_n$ preserves kernel pairs. Since $\gamma_n$ is a
left adjoint, it also preserves coproducts. From (\ref{adj.eq2}) and
(\ref{adj.eq2bis}) we therefore obtain
\begin{equation*}
    \begin{split}
        \fl_n(X)^{(k)}_0 & =(\xi_k\fl_n(X))_0\cong\gamma_\nm(\lcl_n\fcl_\nm(X))_0\cong
        \gamma_\nm(\fcl_\nm(X)\amalg\fcl_\nm(X))= \\
        &
        =\gamma_\nm\fcl_\nm(X)\amalg\gamma_\nm\fcl_\nm(X)=\fl_\nm(X)\amalg\fl_\nm(X).\\
      \fl_n(X)^{(k)}_1 &= (\xi_k\fl_n(X))_1\cong \gamma_\nm(\lcl_n\fcl_\nm(X))_1\cong\\
      &=\gamma_n(\tens{(\fcl_\nm(X)\amalg\fcl_\nm(X))}{p}) =\\
        & =\tens{\gamma_\nm(\fcl_\nm(X)\amalg\fcl_\nm(X))}{p}=\\
        &= \tens{(\fl_\nm(X)\amalg\fl_\nm(X))}{p}.
    \end{split}
\end{equation*}
The structural maps are straightforward to identity, hence $\fl_n(X)$, as
internal category in $\cagpm$ in the $k^{th}$ direction, has the form
$(\fl_\nm(X)\amalg\fl_\nm(X))^p$ as in Lemma \ref{adj.lem2}.\enpr

\subsection{Surjective maps of multinerves.}\label{twolem}
We say a morphism $f:A\rw B$ in $[\dnop,\gp]$ is \emph{levelwise surjective}
if, for all $([x_1]\ldots [x_n])\in\dnop$ the map of groups $f([x_1]\ldots
[x_n]):A([x_1]\ldots [x_n])\rw B([x_1]\ldots [x_n])$ is surjective. The goal of
this section is to prove Proposition \ref{twolem.pro2}, which asserts that if
$\ul_n f$ is surjective, then $\ncl f$ is levelwise surjective. This fact will
be an essential tool in the construction of the functors $C_i$ (Proposition
\ref{spec.pro4}) and in the proof of Lemma \ref{specgen.lem1}, which leads to
the functor $\spi$ of Theorem \ref{specgen.the1}. We deduce Proposition
\ref{twolem.pro2} from an analogous result for the category $\cng$ (Proposition
\ref{catn.pro1}). The latter is stated and used in \cite{celp}. Since, however,
\cite{celp} gives explicit details of the proof only for $n=1$, we have
provided a proof in the Appendix. This proof uses the following lemma and its
corollary.
\begin{lemma}\label{adj.lem1}
    Let $U:\gp\rw\set$ be the forgetful functor and let $R_n:\cng\rw\set$ be the functor
    $R_n(G,d_1,...,d_n,t_1,...,t_n)=UG$. Then for each $\gcl\in\cng$
    there is an isomorphism $\ucl_n\gcl\cong R_n\gcl$, natural in $\gcl$.
\end{lemma}
\prf See Appendix.
\begin{corollary}\label{adj.cor1}
    The functor $\ucl_n:\cng\rw\set$ reflects regular
    epimorphisms.
\end{corollary}
\prf See Appendix.
\begin{proposition}\cite{celp}.\label{catn.pro1}
    Let $\ovl{\ncl}$ be the composite functor
\begin{equation*}
  \ovl{\ncl}:\cng\supar{\gamma_n}\cagp\supar{\ncl} [\dnop,\gp].
\end{equation*}
    If a morphism $f$ in $\cng$ is a regular epi, then $\ovl{\ncl}f$ is levelwise surjective.
\end{proposition}
\prf See Appendix.
\begin{proposition}\label{twolem.pro2}
    Let $f$ be a morphism in $\cagp$ such that $\ul_n f$ is surjective. Then $\ncl
    f$ is levelwise surjective.
\end{proposition}
\prf Recall that in $\set$ the regular epimorphisms are precisely the
surjective maps, thus $\ul_n f$ is a regular epi. Since $\ul_n=\ucl_n\eta_n$
and, by Corollary \ref{adj.cor1}, $\ucl_n$ reflects regular epis, $\eta_n f$ is
a regular epi in $\cng$. Hence, by Proposition \ref{catn.pro1},
$\ovl{\ncl}\eta_n f=\ncl\gamma_n\eta_n f$ is levelwise surjective. But
$\gamma_n\eta_n\cong\id$, hence $\ncl f\cong\ncl\gamma_n\eta_n f$, so $\ncl f$
is also levelwise surjective.\enpr

\section{Strongly contractible cat${^n}$-groups}\label{scontra}
In this section we introduce strongly contractible $\ctn$s and we show (Theorem
\ref{spec.the0}) that for any set $X$, $\fl_n(X)$ is strongly contractible,
where $\fl_n$ is as in Theorem \ref{adj.the1}. As outlined in Section
\ref{orga}, the strong contractibility of $\fl_n(X)$ is an essential ingredient
in the proof of the existence of the functors $C_i:\cagp\rw\cagp$ in
Proposition \ref{spec.pro4} which in turn are used in the construction of the
functor $\spi$ of Theorem \ref{specgen.the1}. We also prove in this section the
important Lemma \ref{spec.lem1}, which asserts the good behaviour of strong
contractibility with respect to certain pullbacks.

\subsection{Definition of strongly contractible cat${^n}$-groups}\label{contra}
In this section we introduce strongly contractible $\ctn$s: these will be used
in Section \ref{specia} to define special $\ctn$s. The idea behind this notion,
which is defined inductively on dimension, is as follows. First of all, we say
that a $\ctn$ is \emph{discrete} if $\ncl\gcl$ is constant. Discrete $\ctn$s
are strongly contractible. A strongly contractible cat$^1$-group is a
cat$^1$-group $\gcl$ which is weakly equivalent to a discrete cat$^1$-group
$\gcl^d$ through a map $d_\gcl:\gcl\rw\gcl^d$ which has a section
$t_\gcl:\gcl^d\rw\gcl$, $\;d_\gcl t_\gcl=\id$; further, to a map $\gcl\rw\gcl'$
of strongly contractible cat$^1$-groups there corresponds functorially a map
$\gcl^d\rw\gcl'^{d}$ and the weak equivalences $d_\gcl$ and $t_\gcl$ are
natural in $\gcl$. Thus strongly contractible cat$^1$-groups are objects which
are ``homotopically discrete in a strong sense" (the map $d_\gcl$
 is a weak equivalence and has a section) and some obvious naturality conditions are satisfied.

For $n>1$, for a $\ctn$ $\gcl$ to be strongly contractible, we first require a
similar condition of being ``homotopically discrete in a strong sense" to hold
for $\gcl$; that is, we require $\gcl$ to be weakly equivalent to a discrete
object $\gcl^d$ through a map which has a section, and we ask some naturality
conditions to be satisfied. Further, we require the same condition to hold for
any ``face" the $\ctn$ is made of. More precisely, in the inductive definition
it will be enough to require that for each direction $1\leq k\leq n$,
$\gcl^{(k)}_0$ and $\gcl^{(k)}_1$ are strongly contractible $\ctmen$s (where
$\gcl^{(k)}_i\in\cagpm$ is as in Corollary \ref{simpdir.cor1}); we will show in
Remark \ref{procon.rem1} that this implies that $\gcl^{(k)}_i$ is strongly
contractible for any $i\geq 0$. The formal definition is as follows
\begin{definition}\label{contra.def1}
    The category $\,\sci\,\cat^1(\,\gp\,)$ of \,strongly contractible cat$^1$-groups is
    the full subcategory of $\cat^1(\gp)$ whose objects $\gcl$ are such that:
    \begin{itemize}
      \item[(i)] If $\gcl$ is discrete, $\gcl\in\sci\cat^1(\gp)$.
      \item[(ii)] There is a functor $D:\sci\cat^1(\gp)\rw\sci\cat^1(\gp)$ such
      that $D(\gcl)$ is discrete. Also, $D(\gcl)=\gcl$ if $\gcl$ is discrete.
      We shall put $D(\gcl)=\gcl^d$ and $D(f)=f^d$.
      \item[(iii)] There are natural transformations $d:\id_{\sci\cat^1(\gp)}\Rightarrow
      D$, $\;t: D\Rightarrow \id_{\sci\cat^1(\gp)}$ satisfying $dt=\id$ such that for each
      $\gcl\in SC\,\cat^1(\gp)$ $\,d_\gcl$ is a weak equivalence. Also,
      $d_\gcl=\id$ if $\gcl$ is discrete.
    \end{itemize}
Suppose, inductively, that we have defined the category $\sci\cat^{n\mi
1}(\!\gp\!)\subset\cat^{n\mi 1}(\!\gp\!)$ of strongly contractible
cat$^{n-1}$-groups. The category $\sci\cagp$ of strongly contractible
cat$^n$-groups is the full subcategory of cat$^n$-groups whose objects $\gcl$
are such that
\begin{itemize}
  \item[(i)] If $\gcl$ is discrete, $\gcl\in\sci\cagp$.
  \item[(ii)] There is a functor $D:\sci\cagp\rw\sci\cagp$ such
      that $D(\gcl)$ is discrete. Also, $D(\gcl)=\gcl$ if $\gcl$ is discrete.
      We shall put $D(\gcl)=\gcl^d$ and $D(f)=f^d$.
  \item[(iii)] There are natural transformations $d:\id_{\sci\cagp}\Rightarrow
      D$, $\;t: D\Rightarrow \id_{\sci\cagp}$ satisfying $dt=\id$ such that, for each
      $\gcl\in SC\,\cagp$, $\,d_\gcl$ is a weak equivalence. Also,
      $d_\gcl=\id$ if $\gcl$ is discrete.
  \item[(iv)] For each $1\leq k \leq n$,
  $\;\gcl^{(k)}_0,\;\gcl^{(k)}_1\in\sci\cagpm$.
\end{itemize}
\end{definition}

\begin{remark}\rm\label{contra.rem1}
    The following fact is immediate from the definition of strong
    contractibility and will be useful in the proof of Lemma
    \ref{specgen.lem1}. Suppose that $\gcl$ is a strongly contractible $\ctn$
    and consider $r$ distinct simplicial directions $k_1,\ldots,k_r$, $\;1\leq k_j\leq
    n$. If we evaluate the multinerve $\ncl\gcl(x_1,\ldots,x_n)$ at each
    $x_{k_j}$ taking either $x_{k_j}=0$ or $x_{k_j}=1$, then we obtain the
    multinerve of a strongly contractible cat$^{n-r}$-group. In fact, from the
    definition, $\gcl^{(k_1)}_0$, $\gcl^{(k_1)}_1$ are strongly contractible,
    hence $(\gcl^{(k_1)}_0)^{(k_2)}_0$, $(\gcl^{(k_1)}_0)^{(k_2)}_1$,
    $(\gcl^{(k_1)}_1)^{(k_2)}_0$, $(\gcl^{(k_1)}_1)^{(k_2)}_1$ are strongly
    contractible; continuing in this way, the claim follows.
\end{remark}

\subsection{A property of strongly contractible cat${^n}$-groups}\label{procon}
The following lemma shows that the notion of strong contractibility behaves
well with respect to certain types of pullbacks. This lemma will be used in the
next section in the proof that $\fl_n(X)$ is a strongly contractible $\ctn$
(Theorem \ref{spec.the0}). Further, it is an essential tool in the proof of
Proposition \ref{spec.pro4} and Lemma \ref{specgen.lem1}, which leads to one of
our main results (Theorem \ref{specgen.the1}) about the functor $\spi$. This
lemma will also be very important to prove Lemma \ref{nerspec.lem1}, which
leads (via Corollary \ref{nerspec.cor1}) to another main theorem (Theorem
\ref{globgen.the1}) about the functor $\dcl_n$.

The proof of this lemma uses the definition of strong contractibility, the fact
that the multinerve functor preserves limits and is fully faithful and the
Quillen model structure on simplicial groups recalled at the end of Section
\ref{simtec}.
\begin{lemma}\label{spec.lem1}$\quad\\$
 a) Consider the pullback of cat$^n$-groups
\begin{equation*}
\begin{diagram}
    P & \rTo^{p_2} & B\\
  \dTo^{p_1} & & \dTo_{g} \\
  A & \rTo_f & C.
\end{diagram}
\end{equation*}
Suppose that $A,B,C$ are strongly contractible cat$^n$-groups and either $\ncl
f$ or $\ncl g$ is levelwise surjective. Then $P$ is strongly contractible, and
$P^d$ is the pullback of $A^d\supar{f^d}C^d\suparle{g^d}B^d$.

b) Let $A,B,C,P\in\cagp$ and suppose $A,B,C$ are strongly contractible. Suppose
there is a pullback in $[\dnop,\gp]$ of the form
\begin{diagram}
    \ncl P & \rTo^{r_2} & \ncl B\\
  \dTo^{r_1} & & \dTo_{G} \\
  \ncl A & \rTo_F & \ncl C,
\end{diagram}
with either $F$ or $G$ levelwise surjective. Then $P$ is strongly contractible.
\end{lemma}
\textbf{Proof of a)}. We start with a preliminary observation. We notice that
there is a commutative diagram
\begin{equation*}
  \xy
    0;/r.17pc/:
    (-45,20)*+{P}="1";
    (-15,20)*+{B}="2";
    (-45,-10)*+{A}="3";
    (-15,-10)*+{C}="4";
    (15,0)*+{P^d}="5";
    (45,0)*+{B^d}="6";
    (15,-30)*+{A^d}="7";
    (45,-30)*+{C^d}="8";
    {\ar^{p_2}"1";"2"};
    {\ar_{p_1}"1";"3"};
    {\ar^{g}"2";"4"};
    {\ar^{f}"3";"4"};
    {\ar^{d_B}"2";"6"};
    {\ar^(0.6){d_P}"1";"5"};
    {\ar^(0.4){d_C}"4";"8"};
    {\ar_{d_A}"3";"7"};
    {\ar^{\til{p}_1}"5";"7"};
    {\ar^{\til{p}_2}"5";"6"};
    {\ar_{f^d}"7";"8"};
    {\ar^{g^d}"6";"8"};
  \endxy
\end{equation*}
where $P^d$ is the pullback of $(f^d,g^d)$. In fact, from naturality of $d$  we
have $d_Cf=f^dd_A,\;\; d_Cg=g^dd_B$. Hence $f^d d_A p_1=d_C f p_1=d_C g p_2=g^d
d_B p_2$. It follows that there is $d_P:P\rw P^d$ with $\til{p}_1 d_P=d_A p_1$,
$\til{p}_2 d_P=d_B p_2$. We claim that $d_P$ is a weak equivalence and that it
has a section $t_P$. In fact, applying the multinerve
$\ncl:\cagp\rw[\dnop,\gp]$ and then the diagonal functor
$diag:[\Delta^{n^{op}}, \gp]\rw[\dop, \gp]$, one obtains a similar diagram of
pullbacks of simplicial groups ($\ncl$ is a right adjoint, so preserves limits;
$diag$ has a left adjoint, given by Kan extension, hence it preserves limits).
Thus we have a commutative diagram of simplicial groups
\begin{diagram}[h=2.2em]
    \diag \ncl A & \rTo^{\diag \ncl f} & \diag \ncl C & \lTo^{\diag \ncl g} & \diag \ncl B\\
    \dTo^{diag\,\ncl d_A} && \dTo_{diag\,\ncl d_C} && \dTo_{diag\,\ncl d_B}\\
    \diag \ncl A^d & \rTo^{\diag \ncl f^d} & \diag \ncl C^d & \lTo^{\diag \ncl g^d} & \diag \ncl B^d
\end{diagram}
in which, by hypothesis, the vertical arrows are weak equivalences and at least
one map in each row, being surjective, is a fibration of simplicial groups.
Since the model structure on simplicial groups is right proper as every
simplicial group is fibrant \cite{qui1}, it follows from \cite[Proposition
13.3.9]{hirs} that the map of pullbacks $\diag \ncl d_P:\diag\ncl P\rw\diag\ncl
P^d$ is a weak equivalence of simplicial groups. Hence, $d_P$ is a weak
equivalence in $\cagp$. We construct a section $t_P$ by considering the
commutative diagram:
\begin{equation*}
  \xy
    0;/r.18pc/:
    (-45,20)*+{P^d}="1";
    (-15,20)*+{B^d}="2";
    (-45,-10)*+{A^d}="3";
    (-15,-10)*+{C^d}="4";
    (15,0)*+{P}="5";
    (45,0)*+{B}="6";
    (15,-30)*+{A}="7";
    (45,-30)*+{C}="8";
    {\ar^(0.4){p_2}"5";"6"};
    {\ar^{p_1}"5";"7"};
    {\ar^{g}"6";"8"};
    {\ar^(0.4){f}"7";"8"};
    {\ar^{t_B}"2";"6"};
    {\ar^(0.6){t_P}"1";"5"};
    {\ar^(0.4){t_C}"4";"8"};
    {\ar_{t_A}"3";"7"};
    {\ar_{\til{p}_1}"1";"3"};
    {\ar^{\til{p}_2}"1";"2"};
    {\ar^{f^d}"3";"4"};
    {\ar^{g^d}"2";"4"};
  \endxy
\end{equation*}
Since, by naturality of $t$,
$ft_A\til{p}_1=t_Cf^d\til{p}_1=t_Cg^d\til{p}_2=gt_B\til{p}_2$, there is a map
$t_P:P^d\rw P$ such that $p_1t_P=t_A\til{p}_1$, $p_2t_P=t_B\til{p}_2$. Since
$d_At_A=\id$, $d_Ct_C=\id$, $d_Bt_B=\id$ it follows, by universality of
pullbacks, that $d_Pt_P=\id$. Finally we observe that all these constructions
are natural. Given a map of diagrams
\begin{diagram}[h=2.2em]
    A & \rTo & C & \lTo & B\\
    \dTo && \dTo && \dTo\\
    A' & \rTo & C' & \lTo & B'
\end{diagram}
this induces a map of the respective pullbacks $h:P\rw P'$, and therefore a map
$h^d: P^d\rw P'^d$. From the universality of pullbacks, it is easily seen that
the assignments $P\rightsquigarrow P^d$, $\;h\rightsquigarrow h^d$ are
functorial and that $d_P$ and $t_P$ are natural in $P$.

We now prove the lemma by induction on $n$. For $n=1$, by the previous
discussion the lemma is true. Suppose it is true for $\nm$ and let $P,A,B,C$ be
as in the hypothesis. From above, we know there is a weak equivalence $d_P:P\rw
P^d$ with a section $t_P$, that $d_P,t_P$ are natural in $P$ and the assignment
$P\rightsquigarrow P^d$, $h\rightsquigarrow h^d$ is functorial. By definition,
to show that $P$ is strongly contractible it remains to prove that, for each
$1\leq k\leq n$, the $\ctmen$s $P^{(k)}_i,\;\;i=0,1$, are strongly
contractible. We claim that each $P^{(k)}_i$ is given by the pullback in
$\cat^{\nm}(\gp)$ as follows
\begin{equation*}
\begin{diagram}[h=2.2em]
    P^{(k)}_i & \rTo & B^{(k)}_i\\
    \dTo && \dTo_{g_i^{(k)}}\\
    A^{(k)}_i & \rTo_{f^{(k)}_i} & C^{(k)}_i.
\end{diagram}
\end{equation*}
In fact since $\ncl$, being right adjoint, preserves pullbacks, $\ncl P$ is
given by the pullback in $[\dnop,\gp]$ of $\ncl A\supar{\ncl f}\ncl
C\suparle{\ncl g}\ncl B$. Since pullbacks in $[\dnop,\gp]$ are computed
pointwise, evaluating at $x_k=i$ we find a corresponding pullback in
$[\dnomen,\gp]$. But recall from Corollary \ref{simpdir.cor1} that $\ncl
P^{(k)}_i(p_1... p_\nm)=\ncl P(p_1... p_{k-1}\, i\, p_{k+1}... p_\nm)$ and
similarly for the other terms. Hence $\ncl P^{(k)}_i$ is the pullback of $\ncl
A^{(k)}_i\!\supar{\ncl f^{(k)}_i}\!\ncl C^{(k)}_i\!\suparle{\ncl
g^{(k)}_i}\!\ncl B^{(k)}_i$. Since $\!\ncl\!$ is fully faithful, this implies
that $P^{(k)}_i$ is the pullback of
$A^{(k)}_i\supar{f^{(k)}_i}C^{(k)}_i\suparle{g^{(k)}_i} B^{(k)}_i$ as claimed.

Since, by hypothesis, $A,B,C$ are strongly contractible $\ctn$s, by definition
$A^{(k)}_i$, $C^{(k)}_i$, $B^{(k)}_i$ are strongly contractible $\ctmen$s for
$i=0,1$ and $1\leq k\leq n$; also from the hypothesis either $\ncl f^{(k)}_i$
or $\ncl g^{(k)}_i$ is levelwise surjective. Thus by the induction hypothesis,
$P^{(k)}_i$ is strongly contractible for $i=0,1$ and $1\leq k\leq n$. This
proves the inductive step.

\bigskip

\textbf{Proof of b)}. Since $\ncl$ is fully faithful, there exist unique
morphisms $p_1:P\rw A$, $p_2:P\rw B$, $f:A\rw C$, $g:B\rw C$ such that $\ncl
p_1=r_1$, $\ncl p_2=r_2$, $\ncl f=F$, $\ncl g=G$ and $f p_1=g p_2$. Since
$\ncl$ is fully faithful, it reflects pullbacks, therefore $P$ is the pullback
of $A\supar{f}C\suparle{g}B$. By part a), $P$ is strongly contractible. \enpr
\begin{remark}\rm\label{procon.rem1}
   If $\gcl$ is a strongly contractible $\ctn$, by definition for each $1\leq k\leq
   n$, $\gcl_0^{(k)}$ and $\gcl_1^{(k)}$ are strongly contractible $\ctmen$s.
   We notice that Lemma \ref{spec.lem1} easily implies that, for each $i\geq
   0$,
   $\gcl_i^{(k)}$ is strongly contractible. In fact, if $i\geq 2$
   we have $\gcl_i^{(k)}=\pro{\gcl_1^{(k)}}{\gcl_0^{(k)}}{i}$. Since
   the source and target maps $\pt_0,\pt_1:\gcl_1^{(k)}\pile{\rw\\
   \rw}\gcl_0^{(k)}$ have a section, $\ncl\pt_0$ and $\ncl\pt_1$ are levelwise
   surjective. So by Lemma \ref{spec.lem1}, $\tens{\gcl_1^{(k)}}{\gcl_0^{(k)}}$
   is strongly contractible. A simple inductive argument shows similarly that
   $\gcl_i^{(k)}$ is strongly contractible for every $i\geq 0$.
\end{remark}

\subsection{A source of strongly contractible cat$^n$-groups}
In this section we show that, for any set $X$, $\fl_n(X)$ is strongly
contractible. This will be very important in the construction of the functors
$C_i:\cagp\rw\cagp$ of Proposition \ref{spec.pro4}, which are used to construct
the functor $\spi$ in Theorem \ref{specgen.the1}. We first need the following
lemma. Recall that, given a morphism $f:\gcl\rw\gcl'$ in $\cagpm$, we denoted
by $\gcl^f$ the object of $\cat(\cagpm)$ as in Lemma \ref{adj.lem2}. Also, for
any $\ctmen$ $\gcl$ we shall denote by $c(\gcl)$ the discrete internal category
in $\cagpm$ with object of objects $\gcl$.
\begin{lemma}\label{catn.lem3}
    Let $f:\gcl\rw\gcl'$ be a morphism in $\cagpm$ with a section
    $v:\gcl'\rw\gcl$, $fv=\id$. Then there is a map $(h,f):\gcl^f\rw
    c(\gcl')$ in $\cat(\cagpm)$, with a section $(sv,v)$
\begin{equation}\label{twolem.eq1bis}
    \begin{diagram}[h=2.7em]
        \tens{(\tens{\gcl}{f})}{\gcl} & \rTo &  \tens{\gcl}{f}& \pile{\rTo^{p_1}\\ \rTo^{p_2}\\
        \lTo_{s}} & \gcl\\
        \dTo&&\dTo_h\uTo_{sv} && \dTo_f\uTo_v\\
        \gcl' & \rTo^{\id} & \gcl' & \pile{\rTo^{\id}\\ \rTo^{\id}\\ \lTo_{\id}} & \gcl'
    \end{diagram}
\end{equation}
where $h=fp_1=fp_2$. Further, $(h,f)$ is a weak equivalence in $\cagp$.
\end{lemma}
\prf It is immediate that $(h,f)$ and $(sv,v)$ are maps in $\cat(\cagpm)$ and
that $(h,f)(sv,v)=\id$. To show that $(h,f)$ is a weak equivalence we proceed
by induction on $n$. In the case $n=1$, by Lemma \ref{catn.lem1} (i),
$\pi_0=\gcl'$ and $\pi_i=0$ for $i>0$. These are also the homotopy groups of
the nerve of $c(\gcl')$. Hence by Lemma \ref{catn.lem1} (ii) the map $(h,f)$ is
a weak equivalence in $\cat(\gp)$.

Suppose, inductively, that the lemma holds for $n-1$ and let $f:\gcl\rw\gcl'$
be as in the hypothesis. Recall that the classifying space of a $\ctn$ is the
classifying space of its multinerve, which is an $n$-simplicial group. As
recalled in Section \ref{simtec} the latter can be computed by regarding the
$n$-simplicial group as a simplicial object in $(\nm)$-simplicial groups, then
computing the classifying space of the $(\nm)$-simplicial groups at each
dimension, thus obtaining a simplicial space, and then taking the geometric
realization of this simplicial space. We are going to use this method to prove
that $B\gcl^f\simeq Bc(\gcl')$.

The idea is simply to ``slice" the diagram (\ref{twolem.eq1bis}) along a
simplicial direction, obtaining for each $[k]\in\dop$ a diagram
(\ref{twolem.eq2}) to which we can apply the inductive hypothesis. This will
yield a map of simplicial spaces which is a levelwise weak equivalence, and
hence induces a weak equivalence of geometric realizations, giving
$B\gcl^f\simeq Bc(\gcl')$. To formally obtain  the diagram (\ref{twolem.eq2}),
we need to pass first to multinerves and use the fact that limits in presheaf
categories are computed pointwise to get for each $k$ a corresponding diagram
of multinerves of cat$^{n-2}$-groups; the fully faithfulness of $\ncl$ will
then give the diagram (\ref{twolem.eq2}). The formal details follow.

The nerve of $\gcl^f\in\cat(\cagpm)$ gives the simplicial object in $\cagpm$
\begin{diagram}
    \cdots & \pile{\rTo\\ \rTo\\ \rTo\\ \rTo\\ } & \tens{(\tens{\gcl}{f})}{\gcl} &
    \pile{\rTo\\ \rTo\\ \rTo\\} &  \tens{\gcl}{f}& \pile{\rTo\\ \rTo\\
        \lTo} & \gcl.
\end{diagram}
Applying the multinerve $\ncl:\cagpm\rw[\Delta^{n-1^{op}},\gp]$ dimensionwise
and recalling that $\ncl$, being a right adjoint, preserves limits, we find a
simplicial object in $[\Delta^{\nm^{op}},\gp]$
\begin{equation}\label{twolem.eq1}
\begin{diagram}[w=2.1em]
    \cdots \!& \pile{\rTo\\ \rTo\\ \rTo\\ \rTo\\} & \tens{(\tens{\ncl\gcl}{\ncl f})}{\ncl\gcl} &
    \pile{\rTo\\ \rTo\\ \rTo\\} &  \tens{\ncl\gcl}{\ncl f}& \pile{\rTo\\ \rTo\\ \lTo} & \ncl\gcl.
\end{diagram}
\end{equation}
At each dimension we have an object of $[\Delta^{\nm^{op}},\gp]$, which we
regard as a simplicial object in $[\Delta^{{n-2}^{op}},\gp]$ along a fixed
direction. We shall denote $(\ncl\gcl)_k=(\ncl\gcl)[k]\in
[\Delta^{{n-2}^{op}},\gp]$ for each $k\geq 0$ with corresponding maps $(\ncl
f)_k:(\ncl\gcl)_k\rw(\ncl\gcl')_k$.

Thus for each $\,[k]\,\in\,\dop$ we obtain the object of \;
$[\dop,\,[\Delta^{{n-2}^{op}},\,\gp]]\cong$ $[\Delta^{\nm^{op}},\gp]$ given by
\small
\begin{equation*}
    ...\pile{\rTo\\ \rTo\\ \rTo\\ \rTo }\!\!  \tens{(\tens{\ncl\gcl}{\ncl
    f})_k}{(\ncl\gcl)_k}\!\!
    \pile{\rTo\\ \rTo\\ \rTo}\!\!  (\tens{\ncl\gcl}{\ncl f})_k \!\! \pile{\rTo\\ \rTo\\
        \lTo} \!\! (\ncl\gcl)_k.
\end{equation*}\normalsize
Since limits in $[\Delta^{{n-2}^{op}},\gp]$ are computed pointwise, this is the
same as
\begin{equation}\label{twolem.eq2ter}
\begin{split}
    \cdots&\pile{\rTo\\ \rTo\\ \rTo\\ \rTo} \tens{(\tens{(\ncl\gcl)_k}{(\ncl f)_k})}{(\ncl\gcl)_k}
   \pile{\rTo\\ \rTo\\ \rTo} \\
   & \pile{\rTo\\ \rTo\\ \rTo}  \tens{(\ncl\gcl)_k}{(\ncl f)_k}  \pile{\rTo\\ \rTo\\
        \lTo}  (\ncl\gcl)_k.
\end{split}
\end{equation}
Since $\gcl\in\cagpm$, $(\ncl\gcl)_k$ is the multinerve of a cat$^{n-2}$-group,
which we denote by $\gcl_k$; that is, $(\ncl\gcl)_k=\ncl\gcl_k$. Since $\ncl$
is fully faithful, (\ref{twolem.eq2ter}) gives rise to a simplicial object in
cat$^{n-2}$-groups
\begin{diagram}
    \cdots & \pile{\rTo\\ \rTo\\ \rTo\\ \rTo\\ } & \tens{(\tens{\gcl_k}{f_k})}{\gcl_k} &
    \pile{\rTo\\ \rTo\\ \rTo\\} &  \tens{\gcl_k}{f_k}& \pile{\rTo\\ \rTo\\
        \lTo} & \gcl_k
\end{diagram}
where $f_k:\gcl_k\rw\gcl'_k$ is such that $\ncl f_k=(\ncl f)_k$. In turn, this
is the nerve of the object $\;\gcl^{f_k}_k$ of
$\cat(\cat^{n-2}(\gp))$, and there are morphisms $\gcl^{f_k}_k\;\pile{\lTo\\
\rTo}\;c(\gcl'_k)$ given by
\begin{equation}\label{twolem.eq2}
\begin{diagram}[h=3.5em]
        \tens{(\tens{\gcl_k}{f_k})}{\gcl_k} & \rTo &  \tens{\gcl_k}{f_k}& \pile{\rTo^{p^{(k)}_1}\\ \rTo^{p^{(k)}_2}\\
        \lTo_{s_k}} & \gcl_k\\
        \dTo && \dTo^{h_k}\dTo_{s_k v_k} && \dTo^{f_k}\uTo_{v_k}\\
        \gcl'_k & \rTo & \gcl'_k & \pile{\rTo^{\id}\\ \rTo^{\id}\\ \lTo_{\id}} & \gcl'_k
\end{diagram}
\end{equation}
with $f_kv_k=\id$, $h_k=f_kp_1^{(k)}=f_kp_2^{(k)}$ where $p_1^{(k)},
p_2^{(k)}:\tens{\gcl_k}{f_k}\rw\gcl_k$ are the two projections. We can
therefore apply the induction hypothesis and conclude that $(h_k,f_k)$ is a
weak equivalence in $\cat(\cat^{n-2}(\gp))=\cagpm$ for each $k\geq 0$; that is,
$B(h_k,v_k)$ is a weak equivalence. Thus we have two simplicial spaces,
$\phi,\phi'\in[\dop,\tp]$ with $\phi[k]=B\gcl^{f_k}_k$, $\phi'[k]=Bc(\gcl'_k)$
and a simplicial map $\phi\rw\phi'$ which is levelwise a weak equivalence.
Hence this map induces a weak equivalence of geometric realizations,
$|\phi|\simeq|\phi'|$. On the other hand, $B\gcl^f=|\phi|$ and
$Bc(\gcl')=|\phi'|$.  Hence $B\gcl^f\simeq Bc(\gcl')$, that is $(h,v)$ is a
weak equivalence. This completes the inductive step. \enpr
\begin{theorem}\label{spec.the0}
    For any set $X$, $\fl_n(X)$ is a strongly contractible cat$^{n}$-group.
\end{theorem}
\prf By definition we need to show that there is a weak equivalence
$g_n:\fl_n(X)\rw\fl_n(X)^d$, where $\fl_n(X)^d$ is discrete, with a section
$t_n$, $g_nt_n=\id$, that for each $1\leq k\leq n$, $(\fl_n(X))_0^{(k)}$ and
$(\fl_n(X))_1^{(k)}$ are strongly contractible, and that all these
constructions are natural. We proceed by induction on $n$.

The two main ingredients in this proof are the form of $\fl_n(X)$ given in
Theorem \ref{adj.the1}, where it is shown that $\fl_n(X)$ as internal category
in $\cagpm$ in the $k^{th}$ direction is isomorphic to
$(\fl_\nm(X)\amalg\fl_\nm(X))^p$, where $p$ has a section; and Lemma
\ref{catn.lem3} which shows that internal categories in $\ctmen$s of this form
are weakly equivalent to discrete internal categories in $\ctmen$s via a map
which has a section. The formal details follow.

In the case $n=1$ recall, by Theorem \ref{adj.the1}, that
$\fl_1(X)\cong(F(X)\amalg F(X))^p$, where $p:F(X)\amalg F(X)\rw F(X)$,
$pu_1=1$, $pu_2=\id$. Since $p$ has a section, we can apply Lemma
\ref{catn.lem3} to conclude that there is a weak equivalence $\fl_1(X)\rw
c(F(X))$ with a section; also $c(F(X))$ is discrete. Further, every morphism
$\fl_1(X)\rw\fl_1(Y)$ corresponds functorially to a morphism $c(F(X))\rw
c(F(Y))$ and the weak equivalence $\fl_1(X)\rw c(F(X))$ is clearly natural in
$\fl_1(X)$. Hence, by definition, $\fl_1(X)$ is a strongly contractible
cat$^1$-group.

Suppose, inductively, that for any set $X$, $\fl_\nm(X)$ is strongly
contractible. From Theorem \ref{adj.the1}, $\fl_n(X)$, as an internal category
in $\cagpm$ in the $k^{th}$ direction, is isomorphic to
$(\fl_\nm(X)\amalg\fl_\nm(X))^p$ where
$p:\fl_\nm(X)\amalg\fl_\nm(X)\rw\fl_\nm(X)$, $pu_1=1$, $pu_2=\id$. Since $p$
has a section, we can apply Lemma \ref{catn.lem3} and conclude that there is a
weak equivalence $r:(\fl_\nm(X)\amalg\fl_\nm(X))^p\rw c(\fl_\nm(X))$ with a
section $l$, $rl=\id$.

By the inductive hypothesis $\fl_\nm(X)$ is strongly contractible, hence by
definition there is a weak equivalence $g_\nm:\fl_\nm(X)\rw\fl_\nm(X)^d$ with a
section $t_\nm$. In turn, this induces a weak equivalence
$c(g_\nm):c(\fl_\nm(X))\rw c(\fl_\nm(X)\!)^d$ with a section $c(t_\nm\!)$.
Since $\fl_\nm(X)^d$ is discrete, so is $c(\fl_\nm(X)\!)^d$. We shall denote
$c(\fl_\nm(X)\!)^d$ by $\fl_n(X)^d$. Summarizing, we have maps:
\begin{equation*}
    (\fl_\nm(X)\amalg\fl_\nm(X))^p\;\pile{\rTo^r\\ \lTo_l}\;c(\fl_\nm(X))
    \;\pile{\rTo^{c(g_\nm)}\\ \lTo_{c(t_\nm)}}\;\fl_n(X)^d.
\end{equation*}
We let $g_n= c(g_\nm)r$ and $t_n=lc(t_\nm)$. Then $g_nt_n=\id$ and $g_n$ is a
weak equivalence, as are $r$ and $c(g_\nm)$.

Every\, morphism \;$\fl_n(X)\;\rw\;\fl_n(Y)$\, corresponds\, functorially\,
to\, a\, morphism $\;\fl_\nm(X)\rw\fl_\nm(Y)$ and therefore, by the induction
hypothesis, to a morphism $\fl_\nm(X)^d\rw\fl_\nm(Y)^d$ and so to a morphism
$\fl_n(X)^d\rw\fl_n(Y)^d$. Also it is clear that $r$ and $l$ are natural in
$\fl_n(X)$; by the induction hypothesis $g_\nm$ and $t_\nm$ are natural in
$\fl_\nm(X)$, hence $c(g_\nm)$ and $c(t_\nm)$ are also natural in
$c(\fl_\nm(X))$. Thus $g_n$ and $t_n$ are natural in $\fl_n(X)$.

To prove that $\fl_n(X)$ is strongly contractible it remains to show that, for
each $1\leq k\leq n$, $(\fl_n(X))_0^{(k)}$ and $(\fl_n(X))_1^{(k)}$ are
strongly contractible. From Theorem \ref{adj.the1}, for each $1\leq k\leq n$,
\begin{equation*}
    \begin{split}
      (\fl_n(X))_0^{(k)} &\cong \fl_\nm(X)\amalg\fl_\nm(X), \\
      (\fl_n(X))_1^{(k)} & \cong\tens{(\fl_\nm(X)\amalg\fl_\nm(X))}{p}.
    \end{split}
\end{equation*}
Since $\fl_\nm$ is a left adjoint, it preserves coproducts, so
$\fl_\nm(X)\amalg\fl_\nm(X)=\fl_\nm(X\amalg X)$; therefore, by the inductive
hypothesis $(\fl_n(X))_0^{(k)}$ is strongly contractible. As for
$(\fl_n(X))_1^{(k)}$, notice that the diagram
\begin{equation*}
    \fl_\nm(X)\amalg\fl_\nm(X)\supar{p}\fl_\nm(X)\suparle{p}\fl_\nm(X)\amalg\fl_\nm(X)
\end{equation*}
satisfies\, the\, hypotheses\, of\, Lemma\; \ref{spec.lem1}; this\, is\,
because\; $\fl_\nm(X)\amalg\fl_\nm(X)=\fl_\nm(X\amalg X)$ and $\fl_\nm(X)$ are
strongly contractible (by the inductive hypothesis) and $\ncl p$ is levelwise
surjective since $pu_2=\id$. It follows by Lemma \ref{spec.lem1} that the
pullback of the above diagram is strongly contractible, hence
$(\fl_n(X))_1^{(k)}$ is strongly contractible. This completes the inductive
step. \enpr

\section{Special cat$\bsy{^n}$-groups}\label{specia}
In this section we define the category $\cagp_S$ of special $\ctn$s and we
prove one of our main theorems (Theorem \ref{specgen.the1}) asserting the
existence of a functor $\spi:\cagp\rw\cagp_S$ and, for each $\gcl\in\cagp$ of a
weak equivalence $\alpha_\gcl:\spi\gcl\rw\gcl$. The functor $\spi$ is built
from the functors $C_i:\cagp\rw\cagp$ of Proposition \ref{spec.pro4}. In
Section \ref{spelow} we illustrate the construction of $\spi$ for $n=2$ and
$n=3$ to gain intuition on the process for general $n$, which is treated in
Section \ref{specgen}. In this section we also prove an important property of
the nerve of special $\ctn$s (Corollary \ref{nerspec.cor1}) which will be
needed in Section \ref{weak} to define the functor $\dcl_n:\cagp_S\rw\dbb_n$
when $n>2$.

\subsection{Definition of special cat${^n}$-groups}\label{defspec}
In this section we use the notion of strong contractibility of Section
\ref{scontra} to define special $\ctn$s.

We start by illustrating the idea behind this notion. Recall from Section
\ref{simtec} that an $n$-simplicial object $\psi\in[\dnop,\ccl]$ in a category
$\ccl$ is said to satisfy the globularity condition if
\begin{itemize}
  \item[(i)] $\psi(0,\mi):\Delta^{\nm^{op}}\rw\ccl$ is constant
  \item[(ii)] $\psi(m_1\ldots m_r\;0\;\mi):\Delta^{{n-r-1}^{op}}\;\rw\;\ccl$ is
  constant for all $\;[m_i]\in\dop$, $\quad i=1,\cdots,r$, $1\leq r\leq n-2$.
\end{itemize}
Given a $\ctn$ $\gcl$, its multinerve $\ncl\gcl$ is a $n$-simplicial object in
groups satisfying the Segal condition (\ref{simpdir.eq1}) in Lemma
\ref{simpdir.lem1}. As easily seen,  it then follows that having conditions (i)
and (ii) for $\ncl\gcl$ is equivalent to requiring:
\begin{itemize}
  \item[(i)'] $\ncl\gcl(0\;\mi):\Delta^{\nm^{op}}\rw\gp$ is constant; that is,
  it is the multinerve of a discrete $\ctmen$.
  \item[(ii)'] $\ncl\gcl(\overset{r}{\overbrace{1\;\ldots\;1}}\;0\;\mi):\Delta^{{n-r-1}^{op}}\;\rw\gp$
  is constant for all $1\leq r\leq n-2$; that is, it is the multinerve of a
  discrete cat$^{n-r-1}$-group.
\end{itemize}
If $\ncl\gcl$ satisfies (i)' and (ii)' then $\gcl\in n\mi\cat(\gp)$. Clearly
this is not the case for every $\ctn$ $\gcl$. The idea of a special $\ctn$ is
to replace in (i)' and (ii)' ``discrete" by ``strongly contractible".

We are now going to give the formal definition. We first fix a terminology.
Given a $\ctn$ $\gcl$ we shall say that its multinerve $\ncl\gcl$ is strongly
contractible if $\gcl$ is strongly contractible in the sense of Definition
\ref{contra.def1}. Given a $\ctn$ $\gcl$, recall that $\ncl\gcl(0,\mi)$ is the
multinerve of a $\ctmen$ and
$\ncl\gcl(\overset{r}{\overbrace{1\;\ldots\;1}}\;0\;\mi)$ is the multinerve of
a cat$^{n-r-1}$-group (this can be seen by repeatedly applying Corollary
\ref{simpdir.cor1}). Thus, according to this terminology, in the next
definition, $\ncl\gcl(0\;\mi)$ strongly contractible means that the $\ctmen$
whose multinerve is $\ncl\gcl(0\;\mi)$ is strongly contractible; similarly
$\ncl\gcl(\overset{r}{\overbrace{1\;\ldots\;1}}\;0\;\mi)$ strongly\,
contractible\, means\, that\, the \,cat$^{n-r-1}$-group\, whose\, multinerve\,
is\; $\ncl\gcl(\overset{r}{\overbrace{1\;\ldots\;1}}\;0\;\mi)$ is strongly
contractible.
\begin{definition}\label{defspec.def1}
    We say that a cat$^2$-group $\gcl$ is special if $\ncl\gcl(0,\mi)$ is strongly contractible.
    For $n>2$ we say that a cat$^n$-group $\gcl$ is special if
    \begin{itemize}
      \item[(i)] $\ncl\gcl(0\;\mi)$ is strongly contractible.
      \item[(ii)] For each $1\leq r\leq n-2$,
      $\ncl\gcl(\overset{r}{\overbrace{1\;\ldots\;1}}\;0\;\mi)$ is strongly
      contractible
    \end{itemize}
    We denote by $\cagp_S$ the full subcategory of $\cagp$ whose objects are
    special cat$^n$-groups.
\end{definition}
\begin{remark}\rm\label{defspec.rem1}
Given $\gcl\in\cagp$, let $\gcl^{(k)}_i\in\cagpm$ be as in Corollary
\ref{simpdir.cor1}. The following facts are immediate from the definitions:
\begin{itemize}
  \item[a)] Every strongly contractible $\ctn$ is special. In fact, recall from
  Remark \ref{contra.rem1} that if $\gcl$ is strongly contractible and we
  evaluate $\ncl\gcl$ at $x_{k_j}=0$ or $x_{k_j}=1$ for $r$ distinct directions
  $k_1,\ldots,k_r$ we find the multinerve of a strongly contractible
  cat$^{n-r}$-group. In particular, conditions (i) and (ii) in the definition
  of special $\ctn$ are satisfied.
  \item[b)] A cat$^2$-group $\gcl$ is special if and only if $\gcl_0^{(1)}$ is
  a strongly contractible cat$^1$-group. If $n>2$, a $\ctn$ $\gcl$ is special
  if and only if $\gcl^{(1)}_0$ is a strongly contractible $\ctmen$ and $\gcl^{(1)}_1$
  is a special $\ctmen$.
\end{itemize}
\end{remark}
\begin{remark}\rm\label{defspec.rem2}
    The following remark can be skipped at a first reading, as it is not
    needed later. We include it for completeness to illustrate the notion of
    special $\ctn$.
    We saw at the beginning of this section that, for the multinerve $\ncl\gcl$ of a $\ctn$ $\gcl$,
    having conditions (i)' and (ii)' is equivalent to having (i)
    and (ii). Similarly, we could ask if, given a special $\ctn$ $\gcl$, for each
    $[m_i]\in\dop$, $i=1,\ldots,r$, $1\leq r\leq n-2$, $\;\ncl\gcl(m_1 \ldots
    m_r\;0\;\mi)$ is also strongly contractible. This is indeed the case and is
    easily proved using Lemma \ref{spec.lem1}. In fact, from Lemma
    \ref{simpdir.lem1}, for every $m_r\geq 2$ we have
    \begin{equation}\label{defspec.eq1}
    \begin{split}
       & \ncl\gcl(1...1\;m_r\;0\;\mi)\cong \\
        & \cong\pro{\ncl\gcl(\overset{r}{\overbrace{1...1\;1}}\;0\;\mi)}
        {\ncl\gcl(\underset{\overset{\underbrace{}}{r-1}}{1...1}\;0\;0\;\mi)}{m_r}.
    \end{split}
    \end{equation}
    Since $\gcl$ is special,
    $\ncl\gcl(\overset{r}{\overbrace{1..1\;1}}\,0\,\mi)$ is strongly
    contractible as well as $\ncl\gcl(\overset{r-1}{\overbrace{1...1}}\;$ $0\;0\;\mi)$. Each
    of the maps $\ncl\gcl(\overset{r}{\overbrace{1...1\;1}}\;0\;\mi)\rw
    \ncl\gcl(\overset{r-1}{\overbrace{1...1}}\;$ $0\;0\;\mi)$ has a section, so is
    levelwise surjective. By repeated use of Lemma \ref{spec.lem1} we conclude
    that the right hand side of (\ref{defspec.eq1}) is strongly contractible,
    hence the same is true for $\ncl\gcl(1...1\;m_r\;0\;\mi)$. By Lemma \ref{simpdir.lem1} we
    also have, for each $m_{r-1}\geq 2$,
    \begin{equation}\label{defspec.eq2}
    \begin{split}
       & \ncl\gcl(1...1\;m_{r-1}\;m_r\;0\;\mi)= \\
        & =\pro{\ncl\gcl(1...1\;m_r\;0\;\mi)}
        {\ncl\gcl(1...1\;0\;m_r\;0\;\mi)}{m_{r-1}}.
    \end{split}
    \end{equation}
    From the previous step, $\ncl\gcl(1...1\;m_r\;0\;\mi)$ is strongly
    contractible. By hypothesis,
    $\ncl\gcl(\overset{r-1}{\overbrace{1...1}}\;0\;\mi)$ is strongly
    contractible hence so is
    $\ncl\gcl(\overset{r-1}{\overbrace{1...1}}\;0\;m_r\;0\;\mi)$ (recall Remark
    \ref{procon.rem1}). Hence by iterated use of Lemma \ref{spec.lem1} we
    conclude from (\ref{defspec.eq2}) that
    $\ncl\gcl(1...1\;m_{r-1}\;m_r\;0\;\mi)$ is strongly contractible.
    Continuing in this way, one concludes that $\ncl\gcl(m_1...m_r\;0\;\mi)$ is strongly
    contractible.
\end{remark}

\subsection{A construction on internal categories}\label{cons}
We describe a well known construction on the category $\cat\ccl$ of internal
categories in a category $\ccl$ with finite limits. This construction will be
applied in Section \ref{apcat} to the case where $\ccl=\cagpm$ to construct the
functors $C_i$ of Proposition \ref{spec.pro4}. The proof of what follows is
quite straightforward, and can be found for instance in \cite{jost}. Let
$\abb\in\cat\ccl$:
\begin{equation*}
    A_1\tiund{A_0}A_1\rTo^m A_1\;\pile{\rTo^{d_0}\\ \rTo_{d_1}\\ \lTo_s}\;A_0
\end{equation*}
and let $p_0:P_0\rw A_0$ be a morphism in $\ccl$. Consider the pullback in
$\ccl$:
\begin{diagram}[h=2.3em]
    P_1 & \rTo^\nu & P_0\times P_0\\
    \dTo^{p_1} && \dTo_{p_0\times p_0}\\
    A_1 & \rTo_{(d_0,d_1)} & A_0\times A_0 .
\end{diagram}
There is a unique object $C(\abb)\in\cat\ccl$ with $C(\abb)_0=P_0$,
$C(\abb)_1=P_1$, such that $(p_1,p_0):C(\abb)\rw\abb$ is an internal functor.
The source and target maps are given by $\pr_1\nu$, $\pr_2\nu$. Further,
suppose that the map $p_0:P_0\rw A_0$ is natural in $\abb$; then given an
internal functor $F:\abb\rw\bbb$ this induces an internal functor
$C(F):C(\abb)\rw C(\bbb)$. From the universality of pullbacks, $C(F)$ is
functorial in $F$ and $C(\abb)\rw\abb$ is natural in $\abb$.

\subsection{Application to the case of cat$^n$-groups}\label{apcat}
In this section we apply the construction of Section \ref{cons} to the case
where $\ccl=\cagpm$ to obtain, for each $1\leq i \leq n$, functors
$C_i:\cagp\rw\cagp$ and a weak equivalence $\alpha^{(i)}_\gcl:C_i(\gcl)\rw\gcl$
natural in $\gcl\in\cagp$ (Proposition \ref{spec.pro4}). These will be used in
Sections \ref{spelow} and \ref{specgen} to construct the functor $\spi$. We
first need the following lemma which shows that, under good conditions, the
construction of Section \ref{cons} applied to $\ccl=\cagpm$ does not change the
homotopy type. As explained in the proof of Proposition \ref{spec.pro4}, in the
construction of Section \ref{cons}, if $p_0$ is the counit of the adjunction
$\fl_\nm\dashv\ul_\nm$, then the hypothesis of Lemma \ref{spec.lem3} is
satisfied (thanks to Proposition \ref{twolem.pro2}). Together with the strong
contractibility of $\fl_n(X)$ (Theorem \ref{spec.the0}) this will imply the
existence of the functor $C_i$ and of $\alpha^{(i)}_\gcl:C_i(\gcl)\rw\gcl$ with
the required properties.
\begin{lemma}\label{spec.lem3}
    Let $1\leq i\leq n$ and $\gcl\in\cagp$. Let
    $\gcl_0^{(i)},\gcl_1^{(i)}\in\cagpm$ be as in Corollary \ref{simpdir.cor1} and let
    $p_0^{(i)}:\hcl_0\rw\gcl_0^{(i)}$ be a morphism in $\cagpm$ such that $\ncl
    p_0^{(i)}$ is levelwise surjective. Consider the pullback in $\cagpm$
\begin{equation}\label{apcat.eq1}
    \begin{diagram}[h=2.3em]
    \hcl_1 & \rTo^\nu & \hcl_0\times\hcl_0\\
    \dTo_{p_1^{(i)}} && \dTo_{p_0^{(i)}\times p_0^{(i)}}\\
    \gcl_1^{(i)} & \rTo_{(d_0^{(i)},d_1^{(i)})} & \gcl_0^{(i)}
    \times\gcl_0^{(i)}.
    \end{diagram}
\end{equation}
Then there is a unique $\cig\in \cagp$ which, as internal category in $\cagpm$
in the $i^{th}$ direction, has object of objects $\cig_{0}^{(i)}=\hcl_0$ and
object of arrows $\cig_1^{(i)}=\hcl_1$; also there is a weak equivalence in
$\cagp$ $\alpha_\gcl^{(i)}:\cig\rw\gcl$.
\end{lemma}
\prf Let\, $\xi_i,\;\xi'_i$ be\, as\, in\, Proposition \,\ref{simpdir.pro1}.
\,Recall\, that \,$(\xi_i\gcl)_0\,=\,\gcl^{(i)}_0$ and
$(\xi_i\gcl)_1=\gcl^{(i)}_1$. By Section \ref{cons}, given the pullback
(\ref{apcat.eq1}), there is a unique $\cg\in\cat(\cat^\nm(\gp))=\cagp$ with
$\cg_0=\hcl_0$, $\cg_1=\hcl_1$ and a functor
$(p_0^{(i)},p_1^{(i)}):\cg\rw\xi_i\gcl$. Let $\cig=\xi'_i\cg$; then
$\xi_i\cig=\cg$ so, by Proposition \ref{simpdir.pro1},
$\hcl_0=\cg_0=(\xi_i\cig)_0=\cig_0^{(i)}$ and
$\hcl_1=\cg_1=(\xi_i\cig)_1=\cig_1^{(i)}$. Also there is a map
$\alpha_\gcl^{(i)}=\xi'_i(p_0^{(i)},p_1^{(i)}):\cig=\xi'_i\cg\rw\xi'_i\xi_i\gcl=\gcl$.
To show that $\alpha_\gcl^{(i)}$ is a weak equivalence it is sufficient to show
that $(p_0^{(i)},p_1^{(i)})$ is a weak equivalence; that is,
$B(p_0^{(i)},p_1^{(i)})$ is a weak equivalence. We show this by induction on
$n$.

The case $n=1$ follows as a special case of \cite[Proposition 6.5]{ekv};
however, we provide a direct proof for completeness. We aim to show that the
map $(p_0,p_1):\cg\rw\gcl$ of $\cat^{1}(\gp)$ induces isomorphisms
$\pi_i(\cg)\cong\pi_i(\gcl)$ for $i=0,1$. Let $d'_0,d'_1:\hcl_1\rw \hcl_0$ be
the source and target maps. Recall from Lemma \ref{catn.lem1} that
$\pi_0(\cg)=\hcl_0/d'_1(\ker d'_0)$, $\pi_1(\cg)=\ker d'_{1_{|_{\ker
d'_0}}}=\ker(d'_0,d'_1)$ where $(d'_0,d'_1)=\nu:\hcl_1\rw\hcl_0\times\hcl_0$;
similarly for $\pi_0(\gcl),\pi_1(\gcl)$. It is easy to see that the pullback
(\ref{apcat.eq1}) in $\gp$  gives rise to a pullback
\begin{diagram}[h=2.3em]
    \ker(d'_0,d'_1)& \rTo & (1,1)\\
    \dTo && \dTo\\
    \ker(d_0^{(i)},d_1^{(i)})& \rTo & (1,1).
\end{diagram}
It follows that $\ker(d'_0,d'_1)\cong \ker(d_0^{(i)},d_1^{(i)})$ and therefore
$\pi_1(\cg)\cong\pi_1(\gcl)$. To show the isomorphism
$\pi_0(\cg)\!\cong\!\pi_0(\gcl)$, first notice that the map
$\pi_0(\cg)\!\rw\pi_0(\gcl)$ is surjective since by hypothesis $p_0^{(i)}$ is a
surjection in groups. To prove injectivity, let $x\in\hcl_0$ with
$p_0^{(i)}(x)=d_1^{(i)}(y)$, $y\in\ker d_0^{(i)}$. Recall from \ref{cons} that
$d'_0=\pr_1 \nu,\;d'_1=\pr_2\nu$. Hence
$(y,(1,x))\in\hcl_1=\gcl_1^{(i)}\tiund{\gcl_0^{(i)}\times\gcl_0^{(i)}}(\tens{\hcl_0}{})$
satisfies $d'_0(y,(1,x))=\pr_1(1,x)=1$ and $d'_1(y,(1,x))=\pr_2(1,x)=x$.
Therefore $x\in d'_1(\ker d'_0)$, which shows that $\pi_0(\cg)\rw\pi_0(\gcl)$
is also injective. This concludes the case $n=1$.

Inductively, suppose that $(p_0^{(i)},p_1^{(i)}):C(\gcl)\rw\xi_i\gcl$ is a weak
equivalence for $n-1$. The idea of the proof is similar to that used in the
proof of Lemma \ref{catn.lem3}. That is, we are going to ``slice" the diagram
(\ref{apcat.eq1}) along a simplicial direction, obtaining for each $[k]\in\dop$
a diagram (\ref{apcat.eq3a}) to which we can apply the inductive hypothesis.
This will yield a map of simplicial spaces whose geometric realizations are
$BC(\gcl)$ and $B\xi_i\gcl$ and which is a levelwise weak equivalence which
will imply $B\cg\simeq B\xi_i\gcl$. The formal details follow.

Since the multinerve functor, being a right adjoint, preserves pullbacks, from
(\ref{apcat.eq1}) we obtain a pullback in $[\Delta^{\nm^{op}},\gp]$
\begin{equation}\label{apcat.eq2}
    \begin{diagram}[h=2.3em]
    \ncl\hcl_1 & \rTo & \ncl\hcl_0\times \ncl\hcl_0\\
    \dTo && \dTo_{\ncl p_0^{(i)}\times \ncl p_0^{(i)}}\\
    \ncl\gcl_1^{(i)} & \rTo & \ncl \gcl_0^{(i)}
    \times \ncl\gcl_0^{(i)}.
    \end{diagram}
\end{equation}
Now regard each $(\nm)$-simplicial group in this pullback as a simplicial
object in $[\Delta^{{n-2}^{op}},\gp]$ in  a fixed direction. For each $k\geq
0$, put $(\ncl\gcl_1^{(i)})_k=(\ncl\gcl_1^{(i)})[k]$ and similarly for the
other terms. Since pullbacks in functor categories are computed pointwise, we
obtain from (\ref{apcat.eq2}) pullbacks in $[\Delta^{{n-2}^{op}},\gp]$
\begin{equation}\label{apcat.eq3}
    \begin{diagram}[h=2.3em]
    (\ncl\hcl_1)_k & \rTo & (\ncl\hcl_0)_k\times (\ncl\hcl_0)_k\\
    \dTo && \dTo\\
    (\ncl\gcl_1^{(i)})_k & \rTo & (\ncl \gcl_0^{(i)})_k
    \times (\ncl\gcl_0^{(i)})_k.
    \end{diagram}
\end{equation}
Since $\gcl_1^{(i)}\in\cagpm$, by Corollary \ref{simpdir.cor1},
$(\ncl\gcl_1^{(i)})_k$ is the multinerve of a cat$^{n-2}$-group, which we call
$\gcl_{1k}^{(i)}$; that is, $\ncl\gcl_{1k}^{(i)}=(\ncl\gcl_1^{(i)})_k$, and
similarly for the other terms. Since $\ncl$ is fully faithful, we obtain from
(\ref{apcat.eq3}) a pullback in $\cat^{n-2}(\gp)$
\begin{equation}\label{apcat.eq3a}
\begin{diagram}[h=2.3em]
    \hcl_{1k} & \rTo & \hcl_{0k}\times\hcl_{0k}\\
    \dTo && \dTo_{p_{0k}^{(i)}\times p_{0k}^{(i)}}\\
    \gcl_{1k}^{(i)} & \rTo_{(d_{0k}^{(i)},d_{1k}^{(i)})} & \gcl_{0k}^{(i)}
    \times\gcl_{0k}^{(i)}.
\end{diagram}
\end{equation}
We denote by $\gcl_{*k}^{(i)}$ the object of $\cat(\cat^{n-2}(\gp))=\cagpm$
given by
\begin{equation*}
    \tens{\gcl_{1k}^{(i)}}{\gcl_{0k}^{(i)}}\rTo\gcl_{1k}^{(i)}\;\pile{\rTo\\ \rTo\\
    \lTo}\;\gcl_{0k}^{(i)}
\end{equation*}
and similarly for $\hcl_{*k}^{(i)}$. Since by hypothesis $\ncl p_0^{(i)}$ is
levelwise surjective, the same is true for $\ncl p_{0k}^{(i)}$. We can
therefore apply the induction hypothesis and deduce that $(p_{0k}^{(i)},
p_{1k}^{(i)}): \hcl_{*k}^{(i)}\rw \gcl_{*k}^{(i)}$ is a weak equivalence in
$\cagpm$ for each $k\geq 0$.

Now let $\phi,\psi\in[\dop,\tp]$ be the simplicial spaces given by
$\phi_k=B\hcl_{*k}^{(i)}$ and $\psi_k=B\gcl_{*k}^{(i)}$ for each $[k]\in\dop$.
Then we have a map of simplicial spaces $\phi\rw\psi$ which is levelwise a weak
equivalence. Hence it induces a weak equivalence of geometric realizations
$|\phi|\simeq|\psi|$. On the other hand it is easy to see that $B\cg=|\phi|$
and $B\xi_i\gcl=|\psi|$; in fact $\ncl\xi_i\gcl$ can be regarded as the
simplicial object in $[\Delta^{{n-1}^{op}},\gp]$ taking each $[k]$ to
$\ncl\gcl_{*k}^{(i)}$; hence, as recalled in Section \ref{simtec}, the
classifying space $B\xi_i\gcl=B\ncl\xi_i\gcl$ can be computed by taking the
classifying space at each dimension and then taking the geometric realization
of the resulting simplicial space. Similarly for $\cg$. In conclusion the map
$(p^{(i)}_0,p^{(i)}_1):\cg\rw\xi_i\gcl$ induces a weak equivalence $B\cg\simeq
B\xi_i\gcl$, which proves the inductive step.\enpr
\begin{proposition}\label{spec.pro4}
    For each  $1\leq i\leq n$,  there is a functor $C_i:\cagp\rw\cagp$ such that,
    for each $\gcl\in\cagp$, the $\ctmen$ $C_i(\gcl)_0^{(i)}$
    is strongly contractible. Further, there is a weak equivalence
    $\alpha^{(i)}_\gcl:C_i(\gcl)\rw\gcl$ in
    $\cagp$, natural in $\gcl$.
\end{proposition}
\prf Let $\gcl\in\cagp$ and let $\gcl_0^{(i)}\vi\gcl_1^{(i)}\in\cagpm$ be as in
Corollary \ref{simpdir.cor1}. Consider the pullback in $\cagpm$
\begin{equation}\label{spec.dia1}
\begin{diagram}[h=2.8em]
    P_1 && \rTo& \fl_\nm\ul_\nm\gcl_0^{(i)}\times \fl_\nm\ul_\nm\gcl_0^{(i)}\\
    \dTo &&& \dTo_{\vep_{\gcl_0^{(i)}}\times\vep_{\gcl_0^{(i)}}}\\
    \gcl_1^{(i)} && \rTo_{(d_0^{(i)},d_1^{(i)})} & \gcl_0^{(i)}\times \gcl_0^{(i)}
\end{diagram}
\end{equation}
where $\vep$ is the counit of the adjunction $\fl_\nm\dashv \ul_\nm$. From the
triangle identities for the adjunction $\fl_\nm\dashv \ul_\nm$,
$\ul_\nm\vep_{\gcl_0}^{(i)}$ has a section, hence it is surjective. Thus by
Proposition \ref{twolem.pro2} $\ncl\vep_{\gcl_0}^{(i)}$ is levelwise
surjective. We can therefore apply Lemma \ref{spec.lem3}.  Hence there exists a
unique $\ctn$ $C_i(\gcl)$ with $C_i(\gcl)_1^{(i)}=P_1$,
$C_i(\gcl)_0^{(i)}=\fl_\nm\ul_\nm\gcl_0^{(i)}$ and a weak equivalence
$\alpha^{(i)}_\gcl:C_i(\gcl)\rw\gcl$ natural in $\gcl$. By Theorem
\ref{spec.the0}, $C_i(\gcl)_0^{(i)}$ is strongly contractible. \enpr
\begin{definition}\label{apcat.def1}
    We will call the functor $C_i$ of Proposition \ref{spec.pro4} the
    cofibrant replacement on cat$^n$-groups in the $i^{th}$ direction.
\end{definition}
\begin{remark}\rm\label{apcat.rem1}
    For the purpose of this paper it is sufficient to consider this just as  a
    terminology. However, it is motivated by the fact that there is indeed a
    model structure on $\ctn$s for which $C_i$ is a functorial cofibrant
    replacement. This can be obtained as an application of the results of
    \cite{ekv}. In \cite{ekv} they consider a model structure on $\cat\ccl$,
    under certain hypotheses, relative to a Grothendieck topology $\tau$ on $\ccl$. We
    take $\ccl=\cagpm$ and the Grothendieck topology to be the one whose basis
    $\kcl$ is given by $\kcl(\gcl)=\{F:\gcl'\rw\gcl\in\text{ Mor }\cagpm,\;
    \ncl F\text{ is levelwise surjective}\}$ for all $\gcl\in \cagpm$.
    Then a map $f$ of $\ctn$s, viewed as a map of internal categories in
    $\cagpm$ in the $i^{th}$ direction, is a weak equivalence if it is fully
    faithful and essentially $\tau$-surjective in the sense defined in
    \cite{ekv}; it is a fibration if it is a $\tau$-fibration as defined in
    \cite{ekv}. One can check that, for each $1\leq i\leq n$, the hypotheses of
    \cite[Theorem 5.5]{ekv} are satisfied, so this defines a model structure on
    $\ctn$s; from \cite{ekv}, the functor $C_i:\cagp\rw\cagp$ is a functorial
    cofibrant replacement. Since we
    do not need the details of this model structure elsewhere in the paper, we
    leave the proof of the above to the interested reader.
\end{remark}

\subsection{Constructing special cat$^n$-groups: low dimensional cases}\label{spelow}
In this section we discuss how to associate functorially to a $\ctn$ a special
one representing the same homotopy type in the cases $n=2$ and $n=3$. In
\ref{specgen} we will extend this construction to general $n$, providing a
functor $\spi:\cagp\rw\cagp_S$ and for each $\gcl\in \cagp$ a weak equivalence
$\alpha_\gcl:\spi\gcl\rw\gcl$, natural in $\gcl$. In principle it is not
necessary to study separately these low-dimensional cases, but we think it is
helpful to do so to gain intuition on the process for general $n$.

We recall the terminology introduced in Section \ref{defspec}: for a $\ctn$
$\gcl$, we say $\ncl\gcl$ is strongly contractible if $\gcl$ is strongly
contractible. Also, in this section and in the next one, we shall denote a
functor $F\in[\dnop,\gp]$ by $F(x_1\,x_2\ldots x_n)$, where the $x_i$'s are
``dummy variables". This notation will be convenient in keeping track of the
simplicial directions along which we will need to evaluate various pullbacks.

The case $n=2$ is simply an instance of Proposition \ref{spec.pro4} in the case
$n=2$ and $i=1$. In fact, the cofibrant replacement of $\gcl$ in direction 1
yields a cat$^2$-group $C_1(\gcl)$ such that $C_1(\gcl)^{(1)}_0$ is a strongly
contractible cat$^1$-group, and a weak equivalence $\alpha^{(1)}_\gcl:
C_1(\gcl)\rw\gcl$ natural in $\gcl$. Thus by definition $C_1(\gcl)$ is a
special cat$^2$-group and we take $\spi\gcl=C_1(\gcl)$.

Consider now the case $n=3$. Recall (see Remark \ref{defspec.rem1} b)) that a
special cat$^3$-group $\gcl$ is one for which $\gcl_0^{(1)}$ is a strongly
contractible cat$^2$-group and $\gcl_1^{(1)}$ is a special cat$^2$-group; in
turn, this is equivalent to saying that $\ncl\gcl(0,\mi,\mi)$ and
$\ncl\gcl(1,0,\mi)$ are strongly contractible.

Given $\gcl\in\cat^3(\gp)$, the cofibrant replacement of $\gcl$ in direction 1
gives a cat$^3$-group $C_1(\gcl)$ which, by Proposition \ref{spec.pro4}, is
such that $C_1(\gcl)^{(1)}_0$ is a strongly contractible cat$^2$-group or
equivalently that $\ncl C_1\gcl(0,\mi,\mi)$ is strongly contractible. This
gives one of the conditions for $C_1(\gcl)$ to be a special cat$^3$-group.
However, the remaining condition that $C_1(\gcl)_1^{(1)}$ is special, or
equivalently that $\ncl C_1\gcl(1,0,\mi)$ is strongly contractible is not in
general satisfied; thus $C_1(\gcl)$ is not in general a special cat$^3$-group.
However, we will see that if we apply $C_1$ not to $\gcl$ but to $C_2\gcl$,
then the resulting cat$^3$-group does satisfy the additional condition that
$\ncl C_1 C_2\gcl(1,0,\mi)$ is strongly contractible; hence $C_1 C_2 \gcl$ is
special. Thus we need to take two successive cofibrant replacements, in
direction 2 and then 1. We will then define $\spi\gcl=C_1 C_2\gcl$ and have a
weak equivalence given by the composite
\begin{equation*}
    \alpha_\gcl:\spi\gcl=C_1C_2\gcl\supar{\alpha^{(1)}_{C_2\gcl}}C_2\gcl\supar{\alpha^{(2)}_\gcl}\gcl.
\end{equation*}
The proof that $C_1 C_2\gcl$ is a special cat$^3$-group uses as fundamental
ingredients: the way in which $C_1$ and $C_2$ are obtained via the pullback
(\ref{spec.dia1}) as in Proposition \ref{spec.pro4} (and therefore a
corresponding pullback in $[\Delta^{2^{op}},\gp]$ upon application of the
multinerve); further, the fact that evaluating this pullback at $x_2=0$
exhibits $\ncl C_1 C_2\gcl(1,0,\mi)$ as a corresponding pullback in (nerves of)
cat$^1$-groups; and finally the fact that this pullback satisfies the
hypotheses of Lemma \ref{spec.lem1} b), which ensures that, under good
conditions, strong contractibility is preserved by pullbacks. This will yield
the required condition that $\ncl C_1 C_2\gcl(1,0,\mi)$ is strongly
contractible and hence $C_1C_2\gcl$ is special.

We now give the details of the proof. Let us put $\rcl=C_2\gcl$. By Proposition
\ref{spec.pro4}, $(C_2\gcl)_0^{(2)}$ is strongly contractible, that is
$\ncl\rcl(x_1\;0\;x_3)=\ncl C_2\gcl(x_1\;0\;x_3)=\ncl(C_2\gcl)_0^{(2)}$ is
strongly contractible. Recall that the notation here is as explained at the
beginning of the section; that is, $x_1$ and $x_3$ are ``dummy variables", so
 $\ncl\rcl(x_1\;0\;x_3)$ means the functor
$\ncl\rcl(\mi\;0\;\mi):\duop\rw\gp$ and similarly in what follows. Consider now
$C_1\rcl$; we know, again by Proposition \ref{spec.pro4}, that $(C_1 \rcl)_0^1$
is a strongly contractible cat$^2$-group; that is, $\ncl C_1\rcl(0\;x_2\;x_3)$
is strongly contractible. We need to show that $\ncl C_1\rcl(1\;0\;x_3)$ is
strongly contractible. Recall, from Proposition \ref{spec.pro4} that we have a
pullback in cat$^2$-groups
\begin{diagram}[h=2.2em]
    (C_1\rcl)^{(1)}_1 & \rTo & \tens{(C_1\rcl)^{(1)}_0}{}\\
    \dTo && \dTo_{\vep\times\vep}\\
    \rcl^{(1)}_1 & \rTo & \tens{\rcl^{(1)}_0}{}
\end{diagram}
where $(C_1\rcl)^{(1)}_0$ is strongly contractible and $\vep$ is the counit of
the adjunction $\fl_2\dashv \ul_2$. Since the multi-nerve functor $\ncl$, being
a right adjoint, preserves pullbacks, we obtain a pullback in
$[\Delta^{2^{op}},\gp]$:
\begin{diagram}
    \ncl C_1\rcl(1\;x_2\;x_3) & \rTo & \tens{\ncl C_1\rcl(0\;x_2\;x_3)}{}\\
    \dTo && \dTo_{\ncl\vep\times\ncl\vep} \\
    \ncl\rcl(1\;x_2\;x_3) & \rTo & \tens{\ncl\rcl(0\;x_2\;x_3)}{}.
\end{diagram}
Since limits in $[\Delta^{2^{op}},\gp]$ are computed pointwise, evaluating at
$x_2=0$, we obtain a pullback in $[\dop,\gp]$:
\begin{equation}\label{spelow.eq1}
 \begin{diagram}
    \ncl C_1\rcl(1\;0\;x_3) & \rTo & \tens{\ncl C_1\rcl(0\;0\;x_3)}{}\\
    \dTo && \dTo_{\ncl\vep_0\times\ncl\vep_0} \\
    \ncl\rcl(1\;0\;x_3) & \rTo & \tens{\ncl\rcl(0\;0\;x_3)}{}.
 \end{diagram}
\end{equation}
We are now going to show that this pullback satisfies the hypotheses of Lemma
\ref{spec.lem1} b); that is that,
 $$\ncl\rcl(1\;0\;x_3),\;\tens{\ncl\rcl(0\;0\;x_3)}{} \text{ and } \tens{\ncl
 C_1\rcl(0\;0\;x_3)}{}$$
are all strongly contractible and that $\ncl\vep_0\times\ncl\vep_0$ is
levelwise surjective. For the latter condition notice that from the triangle
identities for the adjunction $\fl_2\dashv \ul_2$, $\ul_2\vep$ has a section,
so it is surjective. Thus, by Proposition \ref{twolem.pro2}, $\ncl\vep$ is
levelwise surjective. Hence $\ncl\vep_0$, and therefore
$\ncl\vep_0\times\ncl\vep_0$, are levelwise surjective. By Lemma
\ref{spec.lem1} the product of strongly contractible cat$^1$-groups is strongly
contractible. Hence we are reduced to prove the strong contractibility of
\begin{equation*}
    \ncl\rcl(1\;0\;x_3),\quad \ncl\rcl(0\;0\;x_3),\quad \ncl
    C_1\rcl(0\;0\;x_3).
\end{equation*}
Here is where we use the fact, explained at the beginning of the proof, that
$\ncl\rcl(x_1\;0\;x_3)$ is strongly contractible. By definition of strong
contractibility this implies that $ \ncl\rcl(1\;0\;x_3),\; \ncl\rcl(0\;0\;x_3)$
are strongly contractible. As for the remaining one, recall $\ncl
C_1\rcl(0\;x_2\;x_3)$ is strongly contractible. By definition of strong
contractibility, this implies that $\ncl C_1\rcl(0\;0\;x_3)$ is strongly
contractible. Thus the hypotheses of Lemma \ref{spec.lem1} b) are satisfied for
the pullback (\ref{spelow.eq1}) and we conclude that $\ncl
C_1\rcl(1\;0\;x_3)=\ncl C_1 C_2(1\;0\;x_3)$ is strongly contractible, as
required.\enpr

\subsection{Constructing special cat$^n$-groups: general case}\label{specgen}
In this section we construct the functor $\spi:\cagp\rw\cagp_S$ for general $n$
and we show that there is a weak equivalence $\alpha_\gcl:\spi\gcl\rw\gcl$
natural in $\gcl$.

We start with an overview of the method we are going to use. We saw in Section
\ref{spelow} that in the case $n=3$ we had to take two cofibrant replacements
in directions 2 and 1 to obtain a special cat$^3$-group; that is,
$\spi\gcl=C_1C_2\gcl$, $\gcl\in\cat^3(\gp)$. Likewise, for each $n\geq 2$ we
need $(n-1)$ cofibrant replacements; that is, $\spi\gcl=C_1C_2\ldots
C_\nm\gcl$, $\;\gcl\in\cagp$. We need to understand the effect that each of the
$C_i$ produces at each step. At the first step we know from Proposition
\ref{spec.pro4} that $C_\nm\gcl$ is such that $(C_\nm\gcl)^{(\nm)}_0$ is
strongly contractible, that is $\ncl C_\nm \gcl(x_1\ldots x_{n-2}\;0\;x_n)$ is
strongly contractible. If we now apply $C_{n-2}$, we will see that
$C_{n-2}C_\nm\gcl$ is such that
\begin{equation*}
    \begin{split}
        & \ncl C_{n-2}C_\nm\gcl(x_1\ldots x_{n-3}\;0\;x_\nm\;x_n), \\
        & \ncl C_{n-2}C_\nm\gcl(x_1\ldots x_{n-3}\;1\;0\;x_n)
     \end{split}
\end{equation*}
are strongly contractible. Applying further $C_{n-3}$ will yield strong
contractibility for
\begin{equation*}
    \begin{split}
        & \ncl C_{n-3}C_{n-2}C_\nm\gcl(x_1\ldots x_{n-4}\;0\;x_{n-2}\;x_\nm\;x_n), \\
         & \ncl C_{n-3}C_{n-2}C_\nm\gcl(x_1\ldots x_{n-4}\;1\;0\;x_\nm\;x_n), \\
         & \ncl C_{n-3}C_{n-2}C_\nm\gcl(x_1\ldots x_{n-4}\;1\;1\;0\;x_n).
     \end{split}
\end{equation*}
It is easy to guess what is the emerging pattern. To formalize it, we use a
counting device, which is the following notion of \emph{$i$-special
cat$^n$-group}, for $1\leq i\leq \nm$.
\begin{definition}\label{specgen.def1}
    Let $\gcl\in\cagp$, $n\geq 2$. We say $\gcl$ is $(n-1)$-special if
    $\ncl\gcl(x_1\ldots x_{n-2}\;0\;x_n)$ is the multinerve of a strongly
    contractible cat$^\nm$-group. For each $1\leq i< \nm$, we say $\gcl$ is
    $i$-special if
    \begin{itemize}
      \item[(i)] $\ncl\gcl(x_1\ldots x_{i-1}\:0\;x_{i+1}\ldots x_n)$ is the
      multinerve of a strongly contractible cat$^\nm$-group.
      \item[(ii)] For each $1\leq k\leq \nm-i$,
      $\;\ncl\gcl(x_1\ldots x_{i-1}\:\underset{k}{\underbrace{1\;1\ldots 1}}\;0\ldots x_n)$ is
      the multinerve of a strongly contractible cat$^{n-k-1}$-group.
    \end{itemize}
\end{definition}
\begin{remark}\rm\label{specgen.rem1}
    It is clear that, from the definitions, that  1-special = special.
\end{remark}
We will show that, for each $1\leq k\leq \nm$, $C_{n-k}C_{n-k+1}\ldots
C_\nm\gcl$ is $(n-k)$-special. In particular for $k=n-1$ this will yield $C_1
C_2\ldots C_\nm\gcl$ special, as required. To prove the above it is sufficient
to show that
\begin{itemize}
  \item[a)] $C_\nm\gcl$ is $(n-1)$-special.
  \item[b)] For each $1< i\leq \nm$, if $\gcl$ is $i$-special then
  $C_{i-1}\gcl$ is $(i-1)$-special.
\end{itemize}
We know a) is true by Proposition \ref{spec.pro4}. The proof of b) uses a
method similar to the one for $n=3$. Namely, by Proposition \ref{spec.pro4}, we
have a pullback in $\cagpm$
\begin{diagram}[h=2.2em]
    (C_{i-1}\gcl)^{(i-1)}_1 & \rTo & \tens{(C_{i-1}\gcl)^{(i-1)}_0}{}\\
    \dTo && \dTo_{\vep_{\gcl^{(i-1)}_0}\times\vep_{\gcl^{(i-1)}_0}} \\
    \gcl^{(i-1)}_1 & \rTo & \tens{\gcl^{(i-1)}_0}{}.
\end{diagram}
Upon application of the multinerve functor one obtains a corresponding pullback
in $[\Delta^{\nm^{op}},\gp]$:
\begin{equation}\label{specgen.eq1}
\begin{diagram}[labelstyle=\scriptstyle]
    \parbox{53mm}{$\ncl C_{i-1}\gcl(x_1\ldots x_{i-2}\; 1\;x_i\ldots x_n)$} & \rTo &&
    \parbox{58mm}{$\ncl C_{i-1}\gcl(x_1\ldots x_{i-2}\; 0\;x_i\ldots x_n)
    \times\\ \ncl C_{i-1}\gcl(x_1\ldots x_{i-2}\; 0\;x_i\ldots x_n)$}\\
    \dTo &&& \dTo_{\ncl\vep_{\gcl^{(i-1)}_0}\times\ncl\vep_{\gcl^{(i-1)}_0}}\\
    \parbox{48mm}{$\ncl\gcl(x_1\ldots x_{i-2}\; 1\;x_i\ldots x_n)$} & \rTo &&
    \parbox{50mm}{$\ncl \gcl(x_1\ldots x_{i-2}\; 0\;x_i\ldots x_n)
    \times\\ \ncl \gcl(x_1\ldots x_{i-2}\; 0\;x_i\ldots x_n)$}.
\end{diagram}
\end{equation}
We know by Proposition \ref{spec.pro4} that $(C_{i-1}\gcl)^{(i-1)}_0$ is
strongly contractible; that is, $\ncl C_{i-1}\gcl(x_1\ldots x_{i-2}\;
0\;\;x_i\ldots x_n)$ is strongly contractible. To show that $C_{i-1}\gcl$ is
$(i-1)$-special it remains to check that for each $1\leq k\leq n-i$
$\;\ncl\gcl(x_1\ldots x_{i-2}\:$ $\underset{k}{\underbrace{1\;1\ldots
1}}\;0\ldots x_n)$ is strongly contractible. For this we first evaluate the
pullback (\ref{specgen.eq1}) at $x_i=0$ and, using the hypothesis that $\gcl$
is $i$-special, we verify that the resulting pullback satisfies the hypotheses
of Lemma \ref{spec.lem1}. This will yield the strong contractibility of $\ncl
C_{i-1}\gcl(x_1\ldots x_{i-2}\; 1\;0\;x_{i+1}\ldots x_n)$. Further, we will
evaluate the pullback (\ref{specgen.eq1}) at $x_{i+j}=1$ and $x_{i+k-1}=0$ for
each $0\leq j\leq k-2$, $\;1\leq k \leq n-i$ and again verify that to the
resulting pullback we can apply Lemma \ref{spec.lem1}. This will give the
strong contractibility of $\;\ncl C_{i-1}\gcl(x_1\ldots
x_{i-2}\:\underset{k}{\underbrace{1\;1\ldots 1}}\;0\ldots x_n)$ for each $2\leq
k \leq n-i$. In conclusion, from the above $\;\ncl C_{i-1}\gcl(x_1\ldots
x_{i-2}\:\underset{k}{\underbrace{1\;1\ldots 1}}\;0\ldots x_n)$ will be the
multinerve of a strongly contractible object for each $1\leq k \leq n-i$,
proving that $C_{i-1}\gcl$ is $(i-1)$-special. We now give the formal statement
and proofs.
\begin{lemma}\label{specgen.lem1}
    Let $1\leq i \leq n-1$, $n\geq 2$ and let $\gcl$ be an $i$-special cat$^n$-group.
    Then $C_{i-1}\gcl$ is $(i-1)$-special.
\end{lemma}
\prf From the definition of $C_{i-1}\gcl$ (see proof of Proposition
\ref{spec.pro4}), there is a pullback in $\cagpm$:
\begin{diagram}[h=2.2em]
    (C_{i-1}\gcl)^{(i-1)}_1 & \rTo & \tens{(C_{i-1}\gcl)^{(i-1)}_0}{}\\
    \dTo && \dTo_{\vep_{\gcl^{(i-1)}_0}\times\vep_{\gcl^{(i-1)}_0}} \\
    \gcl^{(i-1)}_1 & \rTo & \tens{\gcl^{(i-1)}_0}{}.
\end{diagram}
Since the multinerve functor, being a right adjoint, preserves limits, we
obtain the pullback (\ref{specgen.eq1}) in $[\Delta^{\nm^{op}},\gp]$. We know
from Proposition \ref{spec.pro4} that $(C_{i-1}\gcl)^{(i-1)}_0$ is strongly
contractible; that is, $\ncl C_{i-1}\gcl(x_1\ldots x_{i-2}\; 0\;x_i\ldots x_n)$
is strongly contractible. Hence to show that $C_{i-1}\gcl$ is $(i-1)$-special
it remains to show that for each $1\leq k\leq n-i$, $\;\ncl
C_{i-1}\gcl(x_1\ldots x_{i-2}\:\underset{k}{\underbrace{1\;1\ldots 1}}\;0\ldots
x_n)$ is strongly contractible. We shall treat separately the cases $k=1$ and
$2\leq k\leq n-i$. For the case $k=1$, notice that, as limits in
$[\Delta^{\nm^{op}},\gp]$ are computed pointwise, taking $x_i=0$ in
(\ref{specgen.eq1}) yields a pullback in $[\Delta^{{n-2}^{op}},\gp]$:
\begin{equation}\label{specgen.eq2}
\begin{diagram}[h=2.5em,labelstyle=\scriptstyle]
    \parbox{52mm}{$\ncl C_{i-1}\gcl(x_1... x_{i-2}\; 1\;0\;x_{i+1}... x_n)$} & \rTo &&
    \parbox{55mm}{$\ncl C_{i-1}\gcl(x_1... x_{i-2}\; 0\;0\;x_{i+1}... x_n)
    \times\\ \ncl C_{i-1}\gcl(x_1... x_{i-2}\; 0\;0\;x_{i+1}... x_n)$}\\
    \dTo &&& \dTo_{\ncl\vep_{\gcl^{(i-1)}_0}\times\ncl\vep_{\gcl^{(i-1)}_0}}\\
    \parbox{50mm}{$\ncl\gcl(x_1... x_{i-2}\; 1\;0\;x_{i+1}... x_n)$} & \rTo &&
    \parbox{55mm}{$\ncl \gcl(x_1... x_{i-2}\; 0\;0\;x_{i+1}... x_n)
    \times\\ \ncl \gcl(x_1... x_{i-2}\; 0\;0\;x_{i+1}... x_n)$}.
\end{diagram}
\end{equation}
We claim that (\ref{specgen.eq2}) satisfies the hypotheses of Lemma
\ref{spec.lem1}. In fact, from the triangle identities for the adjunction
$\fl_\nm\dashv\ul_\nm$, $\ul_\nm\vep$ has a section, so it is surjective. Thus
by Proposition \ref{twolem.pro2}, $\ncl\vep$ is levelwise surjective. Hence
$\ncl\vep_{\gcl^{(i-1)}_0}$ and therefore
$\ncl\vep_{\gcl^{(i-1)}_0}\times\ncl\vep_{\gcl^{(i-1)}_0}$, are levelwise
surjective. Also since, by hypothesis, $\gcl$ is $i$-special, $\ncl
\gcl(x_1\ldots x_{i-1}\;0\;x_{i+1}\ldots x_n)$ is strongly contractible; hence,
by definition of strongly contractibility and Remark \ref{contra.rem1},
\begin{equation*}
    \ncl \gcl(x_1\ldots x_{i-2}\; 1\;0\;x_{i+1}\ldots x_n)\quad \text{and}\quad
     \ncl \gcl(x_1\ldots x_{i-2}\; 0\;0\;x_{i+1}\ldots x_n)
\end{equation*}
are strongly contractible. As explained above, $\ncl C_{i-1}\gcl(x_1\ldots
x_{i-2}\; 0$ $\;x_{i}\ldots x_n)$ is strongly contractible, thus, by definition
of strong contractibility and Remark \ref{contra.rem1}, the same is true for
$\ncl C_{i-1}\gcl(x_1\ldots x_{i-2}\; 0\;0\;x_{i+1}\ldots x_n)$. Since, by
Lemma \ref{spec.lem1}, the product of strongly contractible objects is strongly
contractible, we conclude that the hypotheses of Lemma \ref{spec.lem1} b) are
satisfied for the pullback (\ref{specgen.eq2}), and therefore $\ncl
C_{i-1}\gcl(x_1\ldots x_{i-2}\; 1\;0\;x_{i+1}\ldots x_n)$ is strongly
contractible. This deals with the case $k=1$.

We are now going to treat the case $2\leq k\leq n-i$. For this, we evaluate the
pullback (\ref{specgen.eq1}) at $x_{i+j}=1$ and $x_{i+k-1}=0$ for each $0\leq
j\leq k-2$, $\;1\leq k\leq n-i$, obtaining a pullback in
$[\Delta^{{n-k-1}^{op}},\gp]$:
\begin{equation}\label{specgen.eq3}
\begin{diagram}[h=3.5em]
    \parbox{52mm}{$\ncl C_{i-1}\gcl(x_1.. x_{i-2}\;\underset{k}
    {\underbrace{1\,1.. 1}}\;0.. x_n)$} & \rTo &&
    \parbox{60mm}{$\ncl C_{i-1}\gcl(x_1.. x_{i-2}\;0\;\overset{k-1}
    {\overbrace{1.. 1}}\,0.. x_n)
    \times\\ \ncl C_{i-1}\gcl(x_1.. x_{i-2}\;0\;\underset{k-1}
    {\underbrace{1.. 1}}\;0.. x_n)$}\\
    \dTo &&& \dTo\\
    \parbox{52mm}{$\ncl\gcl(x_1.. x_{i-2}\;\underset{k}
    {\underbrace{1\;1.. 1}}\;0\ldots x_n)$} & \rTo &&
    \parbox{50mm}{$\ncl \gcl(x_1.. x_{i-2}\;0\;\overset{k-1}
    {\overbrace{1.. 1}}\;0.. x_n)
    \times\\ \ncl\gcl(x_1.. x_{i-2}\;0\;\underset{k-1}
    {\underbrace{1.. 1}}\;0.. x_n)$}.
\end{diagram}
\end{equation}
We claim that the  pullback (\ref{specgen.eq3}) also satisfies the hypotheses
of Lemma \ref{spec.lem1} b). In fact, since $\gcl$ is $i$-special, $\ncl
\gcl(x_1... x_{i-1}\;\overset{r}{\overbrace{1...1}}\;0... x_n)$ is strongly
contractible for $1\leq r\leq n-1-i$; hence, by definition of strong
contractibility and Remark \ref{contra.rem1}, taking $x_{i-1}=1$ and then
$x_{i-1}=0$ and letting $k=r+1$ we obtain that
\begin{equation*}
    \begin{split}
    & \ncl \gcl(x_1\ldots
    x_{i-2}\:\overset{k}{\overbrace{1\;1\ldots 1}}\;0\ldots x_n)\quad\text{and} \\
    & \ncl \gcl(x_1\ldots x_{i-2}\;0\;\overset{k-1}{\overbrace{1\ldots 1}}\;0\ldots x_n)
    \end{split}
\end{equation*}
are strongly contractible for each $2\leq k\leq n-i$. Also, as explained above,
$\ncl C_{i-1}\gcl(x_1\ldots x_{i-2}\; 0\;x_i\ldots x_n)$ is strongly
contractible, hence by definition and by Remark \ref{contra.rem1} the same is
true for $\ncl C_{i-1} \gcl(x_1\ldots
x_{i-2}\;0\;\overset{k-1}{\overbrace{1\ldots 1}}\;0\ldots x_n)$. As in the case
of the pullback (\ref{specgen.eq2}), since $\ncl\vep$ is levelwise surjective
and the product of strongly contractible objects is strongly contractible, we
conclude that the hypotheses of Lemma \ref{spec.lem1} b) are satisfied for the
pullback (\ref{specgen.eq3}), and therefore $\ncl C_{i-1} \gcl(x_1\ldots
x_{i-2}\:\overset{k}{\overbrace{1\;1\ldots 1}}\;0\ldots x_n)$ is strongly
contractible for each $2\leq k\leq n-i$, as required.\enpr
\begin{theorem}\label{specgen.the1}
    For each $n\geq 2$, there is a functor $\spi:\cagp\rw\cagp_S$ and, for each
    $\gcl\in\cagp$, a weak equivalence $\alpha_\gcl:\spi\gcl\rw\gcl$ natural in
    $\gcl$.
\end{theorem}
\prf Let $\spi\gcl=C_1C_2\ldots C_\nm\gcl$. From Proposition \ref{spec.pro4}
$C_\nm\gcl$ is $(n-1)$-special. Hence, by Lemma \ref{specgen.lem1} for each
$1\leq k\leq n-1$, $C_{n-k}C_{n-k+1}\ldots C_\nm\gcl$ is $(n-k)$-special. In
particular, for $k=\nm$ we deduce that $\spi\gcl$ is special. Let $\alpha_\gcl$
be the composite
\begin{equation*}
    \begin{split}
       & \spi\gcl\rTo^{\alpha^{(1)}_{C_2\ldots C_\nm\gcl}}C_2\ldots C_\nm\gcl
        \rTo^{\alpha^{(2)}_{C_3\ldots C_\nm\gcl}}C_3\ldots C_\nm\gcl\rTo\cdots\\
        & \cdots\rTo^{\alpha^{(n-k-1)}_{C_{n-k}\ldots C_\nm\gcl}}C_{n-k}\ldots C_\nm\gcl
        \rTo\cdots\rTo^{\alpha^{(n-1)}_{\gcl}}\gcl
    \end{split}
\end{equation*}
where each $\alpha^{(i)}$ is as in Proposition \ref{spec.pro4}. By Proposition
\ref{spec.pro4} each $\alpha^{(n-k-1)}_{C_{n-k}\ldots C_\nm\gcl}$ is a weak
equivalence and is natural in $\gcl$, hence the same is true for
$\alpha_\gcl$.\enpr

\bigskip

The following corollary shows that special $\ctn$s model path-connected
$(\np)$-types.
\begin{corollary}\label{spec.cor1}
    \hspace{3mm}There\, is\, an\, equivalence\, of\, categories $\;\cagp_S\bsim\;\simeq \;\;\hcl o(\tp^{(\np)}_*)$
\end{corollary}
\prf We claim that there is an equivalence of categories
$\cagp\bsim\;\simeq\;\cagp_S\bsim$. To prove this, let $f:\gcl\rw\gcl'$ be a
weak equivalence in $\cagp$. From Theorem \ref{specgen.the1}
$f\alpha_\gcl=\alpha_{\gcl'}\spi f$ and $\alpha_\gcl,\alpha_{\gcl'}$ are weak
equivalences. By the 2-out-of-3 property, $Sp\,f$ is also a weak equivalence.
Thus the functor $Sp$ induces a functor $Sp\,:\cagp\bsim\rw\cagp_S\bsim$. Let
$i:\cagp_S\bsim\;\rw\cagp\bsim$ be induced by the inclusion. Then the pair
$(Sp\, , i)$ gives an equivalence of categories. In fact, there are
isomorphisms in $\cagp\bsim$ and $\cagp_S\bsim$, natural in $\gcl$,
$i\:Sp\:\gcl\supar{i\alpha_\gcl} i\gcl=\gcl$ and
$\;Sp\;i\,\hcl=Sp\:\hcl\supar{\alpha_\hcl}\hcl$ for all $\gcl\in\cagp\bsim$,
$\hcl\in\cagp_S\bsim$. Hence $i\,Sp\:\cong\id$, $Sp\:i\cong\id$, so $(Sp\: ,i)$
is an equivalence of categories. By Theorem \ref{catn.the1}, the result
follows.\enpr

\subsection{The nerve of a special cat$^n$-group}\label{nerspec}
The goal of this section is to show (Corollary \ref{nerspec.cor1}) that the
nerve of a special $\ctn$, regarded as an internal category in $\ctmen$s in
direction 1 is not only a simplicial object in $\ctmen$s but in fact a
simplicial object in special $\ctmen$s. This fact will be very important in
Section \ref{globgen} as it will be used in the proof of Lemma
\ref{globgen.lem0} which is needed to define the functor $\dcl_n$ when $n>2$.
Corollary \ref{nerspec.cor1} is a consequence of the following Lemma
\ref{nerspec.lem1}, which shows that special $\ctn$s are well behaved with
respect to certain pullbacks. Lemma \ref{nerspec.lem1} is a direct consequence
of Lemma \ref{spec.lem1}, and it will also be used in the proof of Lemma
\ref{globgen.lem1}.
\begin{lemma}\label{nerspec.lem1}
    Consider the pullback in $\cagp$:
    \begin{diagram}[h=2.2em]
        \pcl & \rTo & \gcl'\\
        \dTo && \dTo_h\\
        \gcl & \rTo_f & \hcl \; .
    \end{diagram}
    Suppose that $\gcl,\hcl,\gcl'$ are special and that either $\ncl f$ or $\ncl h$
    is levelwise surjective. Then $\pcl$ is also special.
\end{lemma}
\prf Since the multinerve functor $\ncl$, being a right adjoint, preserves
pullbacks, there is a pullback in $[\dnop,\gp]$:
\begin{diagram}[h=2.2em]
        \ncl\pcl & \rTo & \ncl\gcl'\\
        \dTo && \dTo_{\ncl h}\\
        \ncl\gcl & \rTo_{\ncl f} & \ncl\hcl \; .
\end{diagram}
Since pullbacks in $[\dnop,\gp]$ are computed pointwise, this gives rise to a
pullback in $[\Delta^{\nm^{op}},\gp]$:
\begin{equation}\label{nerspec.eq1}
\begin{diagram}[h=2.2em]
        \ncl\pcl(0,\mi) & \rTo & \ncl\gcl'(0,\mi)\\
        \dTo && \dTo_{\ncl h_0}\\
        \ncl\gcl(0,\mi) & \rTo_{\ncl f_0} & \ncl\hcl(0,\mi) \;
\end{diagram}
\end{equation}
and, for each $1\leq r\leq n-2$, to a pullback in $[\Delta^{n-r-1^{op}},\gp]$:
\begin{equation}\label{nerspec.eq2}
\begin{diagram}[h=2.2em]
        \ncl\pcl(\overset{r}{\overbrace{1\ldots 1}\;}0\;\mi) & \rTo & \ncl\gcl'(\overset{r}
        {\overbrace{1\ldots 1}\;}0\;\mi)\\
        \dTo && \dTo_{\ncl h_{1\ldots1\;0}}\\
        \ncl\gcl(\overset{r}{\overbrace{1\ldots 1}\;}0\;\mi) & \rTo_{\ncl f_{1\ldots 1\;0}} & \ncl\hcl(\overset{r}
        {\overbrace{1\ldots 1}\;}0\;\mi) \; .
\end{diagram}
\end{equation}
Each object in (\ref{nerspec.eq1}) is the multinerve of a $\ctmen$ and each
object in (\ref{nerspec.eq2}) is the multinerve of a cat$^{n-r-1}$-group. We
shall denote by $\gcl_0$ the $\ctmen$ whose multinerve is $\ncl\gcl(0\;\mi)$
and by $\gcl_{1\ldots 1\;0}$ the cat$^{n-r-1}$-group whose multinerve is
$\ncl\gcl(\overset{r}{\overbrace{1\ldots 1}\;}0\;\mi)$, and similarly for the
other objects. Since the functor $\ncl$ is fully faithful, it reflects
pullbacks, therefore from (\ref{nerspec.eq1}) and  (\ref{nerspec.eq2}) we
obtain a pullback in $\ctmen$s
\begin{equation}\label{nerspec.eq3}
\begin{diagram}[h=2.2em]
        \pcl_0 & \rTo & \gcl'_0\\
        \dTo && \dTo_{h_0}\\
        \gcl_0 & \rTo_{f_0} & \hcl_0 \; ,
\end{diagram}
\end{equation}
and a pullback in cat$^{n-r-1}$-groups
\begin{equation}\label{nerspec.eq4}
\begin{diagram}[h=2.2em]
        \pcl_{1\ldots 1\;0} & \rTo & \gcl'_{1\ldots 1\;0}\\
        \dTo && \dTo_{h_{1\ldots 1\;0}}\\
        \gcl_{1\ldots 1\;0} & \rTo_{f_{1\ldots 1\;0}} & \hcl_{1\ldots 1\;0} \; .
\end{diagram}
\end{equation}
We claim that the pullbacks (\ref{nerspec.eq3}) and (\ref{nerspec.eq4}) satisfy
the hypotheses of Lemma \ref{spec.lem1}. In fact, since by hypothesis
$\gcl,\gcl',\hcl$ are special, by definition all of $\gcl_0,\gcl'_0,\hcl_0$,
$\gcl_{1\ldots 1\;0}$, $\gcl'_{1\ldots 1\;0}$, $\hcl_{1\ldots 1\;0}$ are
strongly contractible. If $\ncl f$ is levelwise surjective, the same is true
for $\ncl f_0$ and for $\ncl f_{1\ldots 1\;0}$. Therefore the hypothesis that
either $\ncl f$ or $\ncl h$ is levelwise surjective implies that either $\ncl
f_0$ or $\ncl h_0$ is levelwise surjective and, similarly, that either $\ncl
f_{1\ldots 1\;0}$ or $\ncl h_{1\ldots 1\;0}$ is levelwise surjective. Hence we
can apply Lemma \ref{spec.lem1} to (\ref{nerspec.eq3}) and (\ref{nerspec.eq4})
and conclude that $\pcl_0$ and $\pcl_{1\ldots 1\,0}$ are strongly contractible;
that is, $\ncl\pcl(0\,\mi)$ and $\ncl\pcl(\overset{r}{\overbrace{1\ldots
1}\;}0\,\mi)$ are strongly contractible, $1\leq r\leq n-2$. By definition, this
means that $\pcl$ is special.\enpr
\begin{corollary}\label{nerspec.cor1}
Let $\ner:\cagp\rw[\dop,\cagpm]$ be the nerve functor constructed by regarding
cat$^{n}$-groups as internal categories in cat$^{n-1}$-groups in direction 1.
Then $\ner$ restricts to a functor
    \begin{equation*}
        \ner:\cagp_S\rw[\dop,\cagpm_S]
    \end{equation*}
\end{corollary}
\prf Let $\gcl\in \cagp_S$. We shall denote for brevity $\gcl_0=\gcl_0^{(1)}$
and $\gcl_1=\gcl_1^{(1)}$; hence $\ner\gcl$ is the simplicial object in
$\cagpm$
\begin{equation*}
    \cdots\pile{\rTo\\ \rTo\\ \rTo\\ \rTo}\tens{\gcl_1}{\gcl_0}\pile{\rTo\\ \rTo\\
    \rTo}\gcl_1\pile{\rTo\\ \rTo\\ \lTo}\gcl_0
\end{equation*}
As noted in Remark \ref{defspec.rem1} b), since $\gcl$ is special, $\gcl_0$ is
strongly contractible (hence special) and $\gcl_1$ is special. We show by
induction on $m$ that, for each $m\geq 2$, $\;\gistr{m}$ is special . Since
$\pt_0:\gcl_1\rw\gcl_0$ has a section, $\ncl\pt_0$ is levelwise surjective.
Hence by Lemma \ref{nerspec.lem1}, $\tens{\gcl_1}{\gcl_0}$ is special. Suppose,
inductively, that $\;\gistr{m-1}$ is special; notice that  $\gistr{m}$ is
isomorphic to the pullback of the diagram
$\gistr{m-1}\rw\gcl_0\suparle{\pt_0}\gcl_1$. By induction hypothesis, the
conditions of Lemma \ref{nerspec.lem1} are satisfied and we conclude that
$\gistr{m}$ is special. \enpr

\section{Internal weak n-groupoids}\label{weak}
In this section we introduce the category $\dbb_n$ of internal weak
$n$-groupoids and we show one of our main theorems (Theorem \ref{globgen.the1})
asserting the existence of a functor $\dcl_n:\cagp_S\rw\dbb_n$ with the
property that $\;B\dcl_n\gcl=B\gcl$ for each $\gcl\in\cagp_S$. We recall from
Section \ref{catn} that the classifying space $B\gcl$ of a $\ctn$ $\gcl$ is the
classifying space $B\ncl\gcl$ of its multinerve. The existence of $\dcl_n$ is
proved by induction on $n$. The first step of induction, $n=2$, is dealt with
in Section \ref{glob2}, while the general case is treated in Section
\ref{globgen}.

\subsection{The category $\dbb_n$ of internal weak $n$-groupoids}\label{intweak}
In this section we define the category $\dbb_n$. The idea behind this notion,
which is defined inductively on dimension, is as follows. Suppose $\vcl$ is a
category with finite limits equipped with a notion of discrete object and a
notion of weak equivalence. Then we could define a category $\wcl_2(\vcl)$ of
Tamsamani weak 2-categories internal to $\vcl$ as the full subcategory of
$[\dop,\vcl]$ consisting of those simplicial objects $\phi$ such that $\phi_0$
is discrete and the Segal maps are weak equivalences.

Let $\vcl=\cat(\gp)$ and consider the nerve functor (which we will denote in
this section by $j_1$), $j_1=\ner:\cat(\gp)\rw[\dop,\gp]$. We take for discrete
objects in $\vcl$ the discrete internal categories; that is, those
$\gcl\in\cat(\gp)$ such that $j_1\gcl$ is constant, and for weak equivalences
the usual weak equivalences in $\cat(\gp)$; that is, the maps $f$ such that $B
j_1 f$ is a weak equivalence. Then we define $\dbb_2$ to be
$\wcl_2(\cat(\gp))$. Inductively, suppose we have defined $\dbb_\nm$ with a
fully faithful embedding $j_\nm:\dbb_\nm\rw[\dnomen,\gp]$. This embedding
allows us to define discrete objects and weak equivalences in $\dbb_\nm$
similarly to the case $n=2$. Then we define $\dbb_n$ to be $\wcl_2(\dbb_\nm)$.

We stress the fact that the terminology ``internal weak $n$-groupoid" is for us
just the name that we give to each object of the category $\dbb_n$ defined in
this way. A priori, other notions of internal weak $n$-groupoid may be given.
\begin{definition}\label{intweak.def1}
    For each $n\geq 1$ we define a category $\dbb_n$ and a fully faithful
    embedding $j_n:\dbb_n\rw[\dnop,\gp]$, as follows.
    For $n=1$, $\dbb_1=\cat(\gp)$ and $j_1=\ner:\cat(\gp)\rw[\dop,\gp]$ is the nerve
    functor. Suppose, inductively, we have defined $\dbb_\nm$ together with a
    full and faithful embedding $j_\nm:\dbb_\nm\rw[\Delta^{\nm^{op}},\gp]$. We
    say an object $\gcl$ of $\dbb_\nm$ is discrete if $j_\nm\gcl$ is constant; we
    say that a map $f$ in $\dbb_\nm$ is a weak equivalence iff $B j_\nm f$ is a
    weak equivalence, where $B:[\Delta^{\nm^{op}},\gp]\rw\tp$ is the
    classifying space functor. We define $\dbb_n$ to be the full subcategory of functors
    $\phi\in[\dop,\dbb_\nm]$ such that
    \begin{itemize}
      \item[(i)] $\phi_0$ is discrete.
      \item[(ii)] For each $k\geq 2$, the Segal maps $\phi_k\rw\fistr{k}$ are
      weak equivalences. We define the embedding $j_n:\dbb_n\rw[\dnop,\gp]$ to be the composite
    \begin{equation*}
        \dbb_n\overset{i}{\hookrightarrow}[\dop,\dbb_\nm]\supar{\ovl{j}_\nm}
        [\dop,[\Delta^{\nm^{op}},\gp]]\cong[\dnop,\gp]
    \end{equation*}
where $i$ is the inclusion and $\ovl{j}_\nm$ is induced by $j_\nm$ (as in
Definition \ref{simtec.def1}). We call $\dbb_n$ the category of internal weak
$n$-groupoids.
    \end{itemize}
\end{definition}
\begin{remark}\rm\label{intweak.rem1}
By definition, $\dbb_n\subset[\dop,\dbb_\nm]$ and $\dbb_1=\cgp$. It follows
that $\dbb_n$ is a full subcategory of
$\underset{\nm}{\underbrace{[\dop,[\dop,\ldots,[\dop}},\cgp]\cdots]$. Since the
latter is isomorphic to $[\dnomen,\cgp]$, we also have a fully faithful
embedding $\dbb_n\rw[\dnomen,\cgp]$.
\end{remark}

\subsection{The basic construction: the discrete nerve}\label{disner}
In Sections \ref{glob2} and \ref{globgen} we will construct a functor $\dcl_n$
which performs the passage from special cat$^n$-group, to the category $\dbb_n$
of internal weak $n$-groupoids in a way that preserves the homotopy type. As we
will see, the idea of the functor $\dcl_n$ is to ``squeeze" the strongly
contractible faces in a special cat$^n$-group to discrete ones, to recover the
globularity condition. The definition of the functor $\dcl_n$ relies on the
iterated use of a basic construction, which we call the \emph{discrete nerve}
functor. In this section we describe this construction in a general context.

Let $\ccl$ be a category with finite limits, $\ccl'\subset\ccl$ a full
subcategory and suppose we are given a functor $D:\ccl'\rw\ccl'$ together with
two natural transformations $d:\id_{\ccl'}\Rightarrow D$,
$t:D\Rightarrow\id_{\ccl'}$ satisfying the condition $dt=\id$. For consistency
with the notation used a later on, we shall denote $D(\gcl)=\gcl^d$, $D(f)=f^d$
for each object $\gcl$ and morphism $f$ in $\ccl'$. For notational simplicity
we shall denote each component $d_{\gcl},\;t_\gcl$ of the natural
transformations simply by $d$ and $t$ (this is a slight abuse of notation, but
the omitted subscript is clear from the context). Let $(\cat\ccl)'$ be the full
subcategory of $\cat\ccl$ whose objects $\gcl$ are such that the object of
objects $\gcl_0$ is in ${\ccl'}$. We define a functor
\begin{equation*}
    \ds\ncl:(\cat\ccl)'\rw[\dop,\ccl]
\end{equation*}
which we call the discrete nerve, as follows. Given $\gcl\in(\cat\ccl)'$
consider the usual nerve of $\gcl$, which is a simplicial object in $\ccl$ of
the form
\begin{equation*}
 \ner\gcl:\quad \cdots\gcl_1\tiund{\gcl_0}\gcl_1\tiund{\gcl_0}\gcl_1\;\pile{\rTo\\ \rTo\\ \rTo\\ \rTo
  }\;\gcl_1\tiund{\gcl_0}\gcl_1 \;\pile{\rTo\\ \rTo\\ \rTo}\;\gcl_1\;\pile{\rTo^{\pt_0}\\
  \lTo^{\sigma_0}\\ \rTo_{\pt_1} }\;\gcl_0.
\end{equation*}
We now modify $\ner\gcl$ by replacing $\gcl_0$ with $\gcl^d_0$ and by modifying
the maps $\pt_0,\pt_1,\sigma_0$ to produce the following simplicial object
$\ds\ncl\gcl$:
\begin{equation*}
  \cdots\gcl_1\tiund{\gcl_0}\gcl_1\tiund{\gcl_0}\gcl_1\;\pile{\rTo\\ \rTo\\ \rTo\\ \rTo
  }\;\gcl_1\tiund{\gcl_0}\gcl_1 \;\pile{\rTo\\ \rTo\\ \rTo}\;\gcl_1\;\pile{\rTo^{d\pt_0}\\
  \rTo^{d\pt_1}\\ \lTo_{\sigma_0 t}}\;\gcl_0^d.
\end{equation*}
In other words, $(ds\:\ncl\gcl)_0=\gcl_0^d$, $(ds\: \ncl\gcl)_n=\gcl_n$ for
$n\geq 1$, we have face and degeneracies $d\pt_i:\gcl_1\rw\gcl_0^d$, $i=0,1$,
$\sigma_0t:\gcl_0^d\rw\gcl_1$ while all the other faces and degeneracies are as
in $Ner \gcl$. Notice that the simplicial identities are satisfied because
$dt=\id$.

Let $F:\gcl\rw\gcl'$ be a morphism in $\cat\ccl$. Then $\ds\ncl
F:\ds\ncl\gcl\rw\ds\ncl\gcl'$ is given by
\begin{equation*}
    (\ds\ncl F)_k=
  \begin{cases}
    F_0^d, & {k=0;} \\
    F_k, & {k>0.}
  \end{cases}
\end{equation*}
The only point to note in checking that $\ds\ncl F$ is a simplicial map are the
relations for $i=0,1$ $\;F_0^d d \pt_i=d\pt'_i F_1$ and $F_1\sigma_0
t=\sigma'_0 t F_0^d$. These follow easily from the naturality of $d$ and $t$
together with the fact that $F$ is a morphism in $\cat\ccl$. This completes the
definition of the functor $\ds\ncl$. The reason for the name ``discrete nerve"
is that in Sections \ref{glob2} and \ref{globgen} we will apply this to the
case where $\ccl=\cagpm$, ${\ccl'}=SC\,\cagpm$, $\,D,d,t$ are as in Definition
\ref{contra.def1}, so that $\gcl^d$ is discrete.

\subsection{A technical lemma}\label{lemtec}
The following lemma will be used in Sections \ref{glob2} and \ref{globgen} to
show that the functor $\dcl_n$ preserves the homotopy type. Its proof (which is
also contained in the proof of \cite[Lemma 4.3]{paol}) amounts to a simple
computation. Let $R:[\dop,\gp]\rw[\dop,\set]$ be the composite
$[\dop,\gp]\supar{\ovl{\ner}}[\dop,[\dop,\set]]\cong[\duop,\set]\supar{\diag}[\dop,\set]$
where $\ner:\gp\rw[\dop,\set]$ is the nerve functor and $\ovl{\ner}$ is induced
by $\ner$ (as in Definition \ref{simtec.def1}).
\begin{lemma}\rm\cite{paol}\it\label{weak.lem01}
    Let $\psi$ (resp. $\chi$) be a bisimplicial group and let
    $\delta_i^h\vi\delta_i^v\vi\sigma_i^h\vi\sigma_i^v$
    (resp. $\mu_i^h\vi\mu_i^v\vi\nu_i^h\vi\nu_i^v$) be the horizontal and vertical face and
    degeneracies operators. Suppose that, for all $p>0$ and $q\geq 0$,
    \begin{equation*}
    \begin{split}
    & \chi_{pq}=\psi_{pq}\\
    & \delta^{v}_i=\mu^{v}_i:\psi_{p,q+1}\rw\psi_{pq}\\
    & \delta^{h}_i=\mu^{h}_i:\psi_{p+1,q}\rw\psi_{pq}\\
    & \sigma^{v}_i=\nu^{v}_i:\psi_{p,q}\rw\psi_{p,q+1}\\
    & \sigma^{h}_i=\nu^{h}_i:\psi_{p,q}\rw\psi_{p+1,q}.
    \end{split}
    \end{equation*}
    Then  $R\,\diag \psi=R\,\diag\chi$.
\end{lemma}
\prf Denote by $N$ the composite
$[\dop,\,\gp]\,\supar{\ovl{\ner}}\,[\dop,\,[\dop,\set]]\cong[\duop,\set]$.
Hence $(N\diag\psi)_{p\mathit{o}}=\{\cdot\}$ and
$(N\diag\psi)_{pq}=\pro{U\psi_{pp}}{}{q}$ for $q>0$, where $U:\gp\rw\set$ is
the forgetful functor. Also we have
\begin{equation*}
    (R\diag\psi)_k=(\diag N\diag\psi)_k=
      \begin{cases}
        \{\cdot\}, & {k=0}, \\
        \pro{U\psi_{kk}}{}{k}, & {k>0,}
      \end{cases}
\end{equation*}
and similarly for $R\diag\chi$. By hypothesis $\psi_{kk}=\chi_{kk}$ for $k > 0$
so that $(R\diag\psi)_k=(R\diag\chi)_k$ for all $k\geq 0$. Further, by
hypothesis, when $n>0$ the face and degeneracy operators
$(\diag\psi)_\np\rw(\diag\psi)_n$ and $(\diag\psi)_n\rw(\diag\psi)_\np$
coincide with the respective ones for $\diag\chi$. This implies that the face
and degeneracy operators in positive dimension of $R\diag\psi$ coincide with
the respective ones for $R\diag\chi$. As for the remaining face and degeneracy
operators, notice that the face maps $U\psi_{11}\pile{\rw\\ \rw}\{\cdot\}$ are
unique as $\{\cdot\}$ is the terminal object and the degeneracy maps
$\{\cdot\}\rw U\psi_{11}$ and $\{\cdot\}\rw U\chi_{11}$ coincide as they both
send $\{\cdot\}$ to the unit of the group $\psi_{11}=\chi_{11}$. In conclusion
$R\diag\psi=R\diag\chi$.\enpr

\subsection{From cubical to globular: case $n=2$}\label{glob2}
In this section we construct the functor $\dcl_2:\cat^2(\gp)_S\rw\dbb_2$ from
special cat$^2$-groups to internal weak 2-groupoids, and show that it preserves
the homotopy type. This is the first step in the inductive argument we will
give in Section \ref{globgen} to construct the functor
$\dcl_n:\cagp_S\rw\dbb_n$ for any $n\geq 2$. The case $n=2$ was also treated in
detail in \cite{paol}. We recall it here for completeness.

In the general setting of Section \ref{disner} let us take $\ccl\!=$
$\cat^1(\gp)$, $\ccl'\!=\sci\cat^1(\gp)$, $\,d,t$ as in Definition
\ref{contra.def1}. By Section \ref{disner} there is a functor
\begin{equation*}
    \ds\ncl:\cat^2(\gp)_S\rw[\dop,\cgp].
\end{equation*}
The next theorem asserts that the image of this functor is in $\dbb_2$ and that
$\ds\ncl$ preserves the homotopy type. We also record the fact that $\ds\ncl$
acts as the identity on globular objects: this will be needed in the inductive
argument in Theorem \ref{globgen.the1}. In what follows the functor
$j_2:\dbb_2\rw[\duop,\gp]$ is as in Definition \ref{intweak.def1}.
\begin{theorem}\label{glob2.the1}
    Let $\dcl_2:\cat^2(\gp)_S\rw[\dop,\cgp]$ be given by $\dcl_2=\ds\ncl$,
    then
    \begin{itemize}
      \item[(i)] $\dcl_2\gcl\in\dbb_2$ for each $\gcl\in\cat^2(\gp)_S$.
      \item[(ii)] $R\diag j_2\dcl_2\gcl=R\diag\ncl\gcl$ and $B\dcl_2\gcl=B\gcl$ for each $\gcl\in\cat^2(\gp)_S$
      \item[(iii)] $\dcl_2\gcl=\gcl$ for $\gcl\in 2\mi\cgp$.
    \end{itemize}
\end{theorem}
\prf Let $\gcl\in\cat^2(\gp)_S$. Then
\begin{equation*}
    (\dcl_2\gcl)_k=
      \begin{cases}
        \gcl_0^d, & {k=0;} \\
        \gcl_1, & {k=1;} \\
        \gistr{k}, & {k>1.}
      \end{cases}
\end{equation*}
Clearly $(\dcl_2\gcl)_0$ is discrete. So, to show that $\dcl_2\gcl\in\dbb_2$,
it remains to prove that for each $k\geq 2$ the Segal maps
\begin{equation*}
\begin{split}
     \eta_k:& \;(\dcl_2\gcl)_k=\gistr{k}\rw\gistpar{k}{\dcl_2}=\\
            & =\pro{\gcl_1}{\gcl_0^d}{k}
\end{split}
\end{equation*}
are weak equivalences in $\cgp$. Consider the case $k=2$. If $\ner:\cgp$
$\rw[\dop,\gp]$ is the nerve functor, there is a commutative diagram of
simplicial groups
\begin{equation}\label{glob2.eq1}
    \begin{diagram}
        \ner\gcl_1\tiund{\ner\gcl_0}\ner\gcl_1 & \rTo^{\ner\eta_2} &
        \ner\gcl_1\tiund{\ner\gcl_0^d}\ner\gcl_1\\
        \dTo^\nu && \dTo_{\nu'}\\
        \ner\gcl_1{\overset{h}{\times}}_{\ner\gcl_0}\ner\gcl_1 & \rTo^{\til{\ner}{\eta_2}} &
        \ner\gcl_1{\overset{h}{\times}}_{\ner\gcl_0^d}\ner\gcl_1
    \end{diagram}
\end{equation}
where ${\overset{h}{\times}}$ denotes the homotopy pullback,
$\til{\ner}{\eta_2}$ is induced by $\ner\eta_2$, and $\nu$ and $\nu'$ are the
canonical maps. Since the source and target maps $\pt_0,\pt_1:\gcl_1\rw\gcl_0$
have a section, the corresponding maps of simplicial groups
$\ner\gcl_1\pile{\rw\\ \rw}\ner\gcl_0$ are surjective. Hence they are
fibrations in the Quillen model structure on simplicial groups \cite{qui1}
recalled at the end of Section \ref{simtec}. Since this model structure is
right proper (every simplicial group is fibrant), it follows from a general
fact \cite[Corollary 13.3.8]{hirs} that $\nu$ and $\nu'$ are weak equivalences.
Since $\ner\gcl_0 \rTo\ner\gcl_0^d$ is a weak equivalence (since $\gcl_0$ is
strongly contractible, as $\gcl$ is special), by the homotopy invariance
property of homotopy pullbacks \cite[Proposition 13.3.4]{hirs},
$\til{\ner}{\eta_2}$ is a weak equivalence. Hence the commutativity of
(\ref{glob2.eq1}) and the 2-out-of-3 property imply that $\ner\eta_2$ is also a
weak equivalence. Thus $\eta_2$ is a weak equivalence in $\cat^1(\gp)$. The
case of the Segal maps for $k>2$ is completely similar. This proves (i).

The proof of (ii) amounts to using the definition of classifying space and
applying Lemma \ref{weak.lem01}. Namely, recall that the classifying space of
$\dcl_2\gcl$ is the classifying space of the bisimplicial group $j_2\dcl_2\gcl$
where $j_2$ is as in Definition \ref{intweak.def1}. Let
$\ncl:\cat^2(\gp)\rw[\Delta^{2^{op}},\gp]$ be the usual multinerve functor. It
is easy to check that the bisimplicial groups $j_2\dcl_2\gcl$ and $\ncl\gcl$
satisfy the hypotheses of Lemma \ref{weak.lem01}. In fact, recall that $j_2$ is
the composite
\begin{equation*}
    \dbb_2\hookrightarrow[\dop,\cat(\gp)]\supar{\ovl{j_1}}[\dop,[\dop,\gp]]\cong[\Delta^{2^{op}},\gp]
\end{equation*}
where $j_1=\ner:\cat(\gp)\rw[\dop,\gp]$ is the nerve functor. Therefore, by
definition of $\dcl_2=\ds\ncl$, we obtain
\begin{equation*}
(j_2\dcl_2\gcl)_{pq}=
      \begin{cases}
        (\ner\gcl_0^d)_q, & {p=0;} \\
        (\ner\gcl_1)_q, & {p=1;} \\
        (\ner(\gistr{p}))_q, & {p>1.}
      \end{cases}
\end{equation*}
On the other hand, the multinerve $\ncl$ is the composite,
\begin{equation*}
    \cat^2(\gp)\supar{\ner}[\dop,\cat(\gp)]\supar{\ovl{\ner}}
[\dop,[\dop,\gp]]\cong[\Delta^{2^{op}},\gp]
\end{equation*}
so that
\begin{equation*}
(\ncl\gcl)_{pq}=
      \begin{cases}
        (\ner\gcl_0)_q, & {p=0;} \\
        (\ner\gcl_1)_q, & {p=1;} \\
        (\ner(\gistr{p}))_q, & {p>1.}
      \end{cases}
\end{equation*}
Hence $(j_2\dcl_2\gcl)_{pq}=(\ncl\gcl)_{pq}$ for $p>0,\;q\geq 0$. From the
definition of $\ds\ncl$, the face and degeneracies operators of $j_2\dcl_2\gcl$
and $\ncl\gcl$ differ at most in the case of the horizontal face maps from
$(1,q)$ to $(0,q)$ $q\geq 0$, horizontal degeneracy maps from $(0,q)$ to
$(1,q)$, $q\geq 0$, vertical face maps $(0,q+1)\rw(0,q)$, $q\geq 0$ and
vertical degeneracy maps $(0,q)\rw(0,q+1)$, $q\geq 0$. Hence the hypotheses of
Lemma \ref{weak.lem01} are satisfied. We conclude that $R\diag
j_2\dcl_2\gcl=R\diag\ncl\gcl$, and therefore $B\dcl_2\gcl=|R\diag
j_2\dcl_2\gcl|=|R\diag\ncl\gcl|=B\gcl$, which is (iii).

Finally, if $\gcl\in 2\mi\cat(\gp)$, $\gcl_0$ is discrete, so
$\gcl_0^d=\gcl_0$, which immediately implies $\dcl_2\gcl=\gcl$, which is
(iii).\enpr

\subsection{From cubical to globular: general case}\label{globgen}
In this section we construct the functor $\dcl_n:\cagp_S\rw\dbb_n$ for any
$n\geq 2$. We start with an overview of the method we are going to use.

We saw in Section \ref{glob2} that $\dcl_2$ is the discrete nerve $\ds\ncl$. In
higher dimensions, the definition of $\dcl_n$ amounts to appropriately
iterating the discrete nerve construction. An essential fact that allows us to
do this is Corollary \ref{nerspec.cor1} which guarantees that the nerve of a
special $\ctn$ is a simplicial object in special cat$^{\nm}$-groups. Namely, in
the general setting of Section \ref{disner} we take $\ccl=\cagpm$,
$\ccl'=\sci\cagpm$ and $D,d,t$ as in Definition \ref{contra.def1}. Then we
obtain a functor $\ds\ncl:(\cat\ccl)'\rw[\dop,\ccl]$ where, as in Section
\ref{glob2}, $(\cat \ccl)'$ is the full subcategory of $\cat\ccl$ whose objects
$\gcl$ are such that $\gcl_0\in\ccl'$. Recall that if $\gcl$ is a special
$\cagp$ regarded as an internal category in $\cagpm$ in direction 1, the object
of objects $\gcl_0$ is strongly contractible; hence there is an inclusion
$\cagp_S\subset (\cat\ccl)'$, where $\ccl$ and $\ccl'$ are as above. The key
point is that from Corollary \ref{nerspec.cor1} the functor
$\ds\ncl:(\cat\ccl)'\rw[\dop,\ccl]$, when restricted to the subcategory
$\cagp_S$ of $(\cat\ccl)'$, yields the following:
\begin{lemma}\label{globgen.lem0}
    The discrete nerve construction yields a functor
\begin{equation*}
    \ds\ncl:\cagp_S\rw[\dop,\cat^{\nm}(\gp)_S].
\end{equation*}
\end{lemma}
\prf Let $\gcl\in\cagp_S$. We shall denote for brevity $\gcl_0=\gcl^{(1)}_0$
and $\gcl_1=\gcl^{(1)}_1$. Then $\ds\ncl\gcl$ is given by
\begin{equation}\label{globgen.eq2}
  \cdots\gcl_1\tiund{\gcl_0}\gcl_1\tiund{\gcl_0}\gcl_1\;\pile{\rTo\\ \rTo\\ \rTo\\ \rTo
  }\;\gcl_1\tiund{\gcl_0}\gcl_1 \;\pile{\rTo\\ \rTo\\ \rTo}\;\gcl_1\;\pile{\rTo^{d\pt_0}\\
  \rTo^{d\pt_1}\\ \lTo_{\sigma_0 t}}\;\gcl_0^d.
\end{equation}
Since $\gcl_0^d$ is discrete, it is strongly contractible and hence special; by
Corollary \ref{nerspec.cor1}, $\gcl_1$ and $\gistr{k}$ are special for all
$k\geq 2$.\enpr\bigskip

Suppose, inductively, that we have constructed the functor
$\dcl_\nm\!:\!\cat^{\nm}(\gp)_S$ $\rw\dbb_\nm$. By Lemma \ref{globgen.lem0} we
can then apply this functor levelwise in (\ref{globgen.eq2}) to obtain a
simplicial object in $\dbb_\nm$. That is, we define
\begin{equation*}
    \dcl_n=\ovl{\dcl_\nm}\circ\ds\ncl:\cagp_S\rw[\dop,\dbb_\nm]
\end{equation*}
where $\ovl{\dcl_\nm}:[\dop,\cagpm_S]\rw[\dop,\dbb_\nm]$ is induced by
$\dcl_\nm$ as in Definition \ref{simtec.def1}; that is,
$(\ovl{\dcl_\nm}\gcl)_k=\dcl_\nm\gcl_k$ for all $\gcl\in[\dop,\cagpm_S]$. In
Theorem \ref{globgen.the1} we are going to show by induction that the image of
the functor $\dcl_n$ is in $\dbb_n$, that $B\dcl_n\gcl=B\gcl$ for all
$\gcl\in\cagp_S$ and that $\dcl_n$ acts as the identity on globular objects.
The first step in the induction $(n=2)$ is Theorem \ref{glob2.the1} while the
inductive hypothesis consists in assuming that these properties hold for
$\dcl_\nm$.

There are two main ingredients in the proof that $\dcl_n\gcl\in\dbb_n$. One is
the observation that, for each $n\geq 2$, $\dcl_n$ preserves pullbacks of
diagrams of the type $\gcl\supar{f}\hcl\suparle{h}\gcl'$ where $\gcl,\gcl'$ are
special $\ctn$s, $\hcl$ is discrete and either $f$ or $h$ has a section. Note
that since either $f$ or $h$ has a section, either $\ncl f$ or $\ncl h$ is
levelwise surjective, hence,  by Lemma \ref{nerspec.lem1} the pullback
$\gcl\tiund{\hcl}\gcl'$ is special, and it makes sense to consider
$\dcl_n(\gcl\tiund{\hcl}\gcl')$. We show in Lemma \ref{globgen.lem1} that
$\dcl_n(\gcl\tiund{\hcl}\gcl')=\dcl_n\gcl\tiund{\dcl_n\hcl}\dcl_n\gcl'$. The
proof  of Lemma \ref{globgen.lem1} relies on the definition of $\dcl_n$ and
Lemma \ref{spec.lem1} which is the one that asserts the good behavior of strong
contractibility with respect to certain pullbacks. The second ingredient in the
proof that $\dcl_n\gcl\in\dbb_n$ is the observation that, for each $k\geq 2$,
the maps
$\gistr{k}\rw\gcl_1\tiund{\gcl_0^d}\overset{k}{\cdots}\tiund{\gcl_0^d}\gcl_1$
are weak equivalences. This fact is proved easily by passing to the diagonals
of the corresponding multinerves and then arguing similarly to the case $n=2$.

These two facts together with the induction hypothesis that
$B\dcl_\nm\gcl=B\gcl$ and $\dcl_\nm\gcl=\gcl$ if $\gcl$ is globular, yield the
weak equivalences of the Segal maps in $\dcl_n\gcl$. For instance, in the case
$k=2$, the weak equivalence
$\tens{\gcl_1}{\gcl_0}\simeq\tens{\gcl_1}{\gcl_0^d}$ together with the fact
that $B\dcl_\nm\gcl=B\gcl$, give a weak equivalence
$\dcl_\nm(\tens{\gcl_1}{\gcl_0})\simeq\dcl_\nm(\tens{\gcl_1}{\gcl_0^d})$. On
the other hand, since $\gcl$ is special, $\gcl_1$ is special (recall Remark
\ref{defspec.rem1} b)); further the map $d\pt_0:\gcl_1\rw\gcl_0^d$ has a
section. Hence by Lemma \ref{globgen.lem1},
$\dcl_\nm(\tens{\gcl_1}{\gcl_0^d})=\tens{\dcl_\nm\gcl_1}{\dcl_\nm\gcl_0^d}$.
Using the induction hypothesis $\dcl_\nm\gcl^d_0=\gcl_0^d$, one finally obtains
the weak equivalence
\begin{equation*}
    \begin{split}
      (\dcl_n\gcl)_2 & =\dcl_\nm(\tens{\gcl_1}{\gcl_0})\simeq\dcl_\nm(\tens{\gcl_1}{\gcl_0^d})= \\
        &
        =\tens{\dcl_\nm\gcl_1}{\dcl_\nm\gcl_0^d}=\tens{(\dcl_n\gcl)_1}{\gcl_0^d}.
    \end{split}
\end{equation*}
Similarly when $k>2$; thus the Segal maps in $\dcl_n\gcl\subset[\dop,\dbb_\nm]$
are weak equivalences; since $(\dcl_n\gcl)_0=\dcl_\nm\gcl_0^d=\gcl^d_0$ is
discrete, we conclude that $\dcl_n\gcl\in\dbb_n$. Finally the fact that
$B\dcl_n\gcl=B\gcl$ will follow from a computation using Lemma
\ref{weak.lem01}; the fact that $\dcl_n\gcl=\gcl$ when $\gcl\in n\mi\cat(\gp)$
is straightforward from the definition and the inductive hypothesis.

We now come to the formal statements and proofs.
\begin{lemma}\label{globgen.lem1}
    Let $\ds\ncl:\cagp_S\rw[\dop,\cat^\nm(\gp)_S]$ be as in
    Lemma \ref{globgen.lem0}. Let $\dcl_n:\cagp_S\rw[\dop,\dbb_\nm]$ be defined
    inductively by $\dcl_2=\ds\ncl$, $\dcl_n=\ovl{\dcl}_\nm\circ\ds\ncl$ for
    $n>2$. Let $\gcl\supar{f}\hcl\suparle{h}\gcl'$ be morphisms in $\cagp$
    with $\gcl,\gcl'$ special and $\hcl$ discrete, and suppose that
    either $f$ or $h$ has a section. Then $\dcl_n(\gcl\tiund{\hcl}\gcl')=
    \dcl_n\gcl\tiund{\dcl_n\hcl}\dcl_n\gcl'$.
\end{lemma}
\prf Since $\hcl$ is discrete, it is in particular strongly contractible, hence
special; since either $f$ or $h$ has a section, either $\ncl f$ or $\ncl h$ is
levelwise surjective so, by Lemma \ref{nerspec.lem1}, $\gcl\tiund{\hcl}\gcl'$
is special; hence it makes sense to consider $\dcl_n(\gcl\tiund{\hcl}\gcl')$.
We claim that, given $f$ and $h$ as in the hypothesis, we have
\begin{equation}\label{globgen.eq3}
    \ds\ncl(\gcl\tiund{\hcl}\gcl')=\ds\ncl\gcl\tiund{\ds\ncl\hcl}\ds\ncl\gcl'.
\end{equation}
In showing (\ref{globgen.eq3}) we shall denote, for each special $\ctn$ $\gcl$
and $k\geq 0$, $\gcl_k=(\ner\gcl)_k$ where $\ner:\cagp_S\rw[\dop,\cagpm_S]$ is
the nerve functor constructed by regarding $\ctn$s as internal categories in
$\ctmen$s in direction 1. Also, for each special $\ctn$s $\gcl$, $\ds\ncl\gcl$
is a simplicial object in $\cagpm_S$ by Lemma \ref{globgen.lem0} and we put
$(\ds\ncl\gcl)_k=(\ds\ncl\gcl)[k]$. Since $\gcl$ and $\gcl'$ are special,
$\gcl_0$ and $\gcl'_0$ are strongly contractible; also $\hcl_0$ is discrete
since $\hcl$ is discrete, so $\hcl_0$ is strongly contractible. From the
hypothesis, either $f_0:\gcl_0\rw\hcl_0$ or $h_0:\gcl'_0\rw\hcl_0$ has a
section, so either $\ncl f_0$ or $\ncl h_0$ is levelwise surjective. Thus the
hypotheses of Lemma \ref{spec.lem1} are satisfied for the diagram
$\gcl_0\supar{f_0}\hcl_0\suparle{h_0}\gcl'_0$ so that
$\gcl_0\tiund{\hcl_0}\gcl'_0$ is strongly contractible and, as in Lemma
\ref{spec.lem1}
\begin{equation}\label{globgen.eq4}
    (\gcl_0\tiund{\hcl_0}\gcl'_0)^d=\gcl_0^d\tiund{\hcl_0^d}{\gcl'}_0^d.
\end{equation}
Using (\ref{globgen.eq4}) and the definition of $\ds\ncl$, we then compute
\begin{equation}\label{globgen.eq5}
    \begin{split}
       & (\ds\ncl(\gcl\tiund{\hcl}\gcl'))_0=((\gcl\tiund{\hcl}\gcl')_0)^d=
       (\gcl_0\tiund{\hcl_0}\gcl'_0)^d=\gcl_0^d\tiund{\hcl_0^d}{\gcl'}_0^d\\
        & =(\ds\ncl\gcl)_0\tiund{(\ds\ncl\hcl)_0}(\ds\ncl\gcl')_0=
        (\ds\ncl\gcl\tiund{\ds\ncl\hcl}\ds\ncl\gcl')_0.
    \end{split}
\end{equation}
For $k>0$ we have, from the definition of $\ds\ncl$,
\begin{equation}\label{globgen.eq6}
    \begin{split}
       & (\ds\ncl(\gcl\tiund{\hcl}\gcl'))_k=(\gcl\tiund{\hcl}\gcl')_k=\gcl_k\tiund{\hcl_k}\gcl'_k= \\
        &(\ds\ncl\gcl)_k\tiund{(\ds\ncl\hcl)_k}(\ds\ncl\gcl')_k=
        (\ds\ncl\gcl\tiund{\ds\ncl\hcl}\ds\ncl\gcl')_k.
    \end{split}
\end{equation}
From (\ref{globgen.eq5}) and (\ref{globgen.eq6}), now (\ref{globgen.eq3})
follows.

We are now going to use (\ref{globgen.eq3}) to prove the lemma by induction on
$n$. For $n=2$, since $\dcl_2=\ds\ncl$ the lemma holds by (\ref{globgen.eq3}).
Suppose, inductively, that the lemma holds for $\nm$ and let
$\gcl\supar{f}\hcl\suparle{h}\gcl'$ be morphisms in $\ctn$s as in the
hypothesis. We claim that we can apply the inductive hypothesis to the diagram
\begin{equation}\label{globgen.eq6a}
    (\ds\ncl\gcl)_k\rTo^{(\ds\ncl f)_k}(\ds\ncl\hcl)_k\lTo^{(\ds\ncl h)_k}(\ds\ncl\gcl')_k
\end{equation}
in $\cagpm$. In fact since $\gcl$ and $\gcl'$ are special, from Lemma
\ref{globgen.lem0} also $(\ds\ncl\gcl)_k$ and $(\ds\ncl\gcl')_k$ are special.
Since $\hcl$ is discrete, $(\ds\ncl\hcl)_k$ is discrete. Since either $f$ or
$h$ has a section, either $(\ds\ncl f)_k$ or $(\ds\ncl h)_k$ has a section;
this is because $(\ds\ncl f)_0=f_0^d$ and $(\ds\ncl f)_k=f_k$ for $k>0$ and
similarly for $\ds\ncl h$ so that if, for instance $fs=\id$, then $f_ks_k=\id$
for all $k>0$ and $f_0^ds_0^d=\id$.

The inductive hypothesis applied to (\ref{globgen.eq6a}) together with
(\ref{globgen.eq3})  yields, for each $k\geq 0$,
\begin{equation*}
    \begin{split}
       & (\dcl_n(\gcl\tiund{\hcl}\gcl'))_k=\dcl_\nm(\ds\ncl(\gcl\tiund{\hcl}\gcl'))_k=
        \dcl_\nm(\ds\ncl\gcl\tiund{\ds\ncl\hcl}\ds\ncl\gcl')_k=\\
        & =\dcl_\nm((\ds\ncl\gcl)_k\tiund{(\ds\ncl\hcl)_k}(\ds\ncl\gcl')_k)=\\
       &= \dcl_\nm(\ds\ncl\gcl)_k\tiund{\dcl_\nm(\ds\ncl\hcl)_k}\dcl_\nm(\ds\ncl\gcl')_k=\\
        &=(\dcl_n\gcl)_k\tiund{(\dcl_n\hcl)_k}(\dcl_n\gcl')_k=
        (\dcl_n\gcl\tiund{\dcl_n\hcl}\dcl_n\gcl')_k.
    \end{split}
\end{equation*}
Hence $\dcl_n(\gcl\tiund{\hcl}\gcl')=\dcl_n\gcl\tiund{\dcl_n\hcl}\dcl_n\gcl'$,
proving the inductive step.\enpr
\begin{theorem}\label{globgen.the1}
    Let $\dcl_n:\cagp_S\rw[\dop,\dbb_\nm]$ be as in Lemma \ref{globgen.lem1}.
    Then
    \begin{itemize}
      \item[(i)] $\dcl_n\gcl\in\dbb_n$ for each $\gcl\in\cagp_S$.
      \item[(ii)] $R\diag j_n\dcl_n\gcl=R\diag\ncl\gcl$ and $B\dcl_n\gcl=B\gcl$ for each $\gcl\in\cagp_S$.
      \item[(iii)] $\dcl_n\gcl=\gcl$ for each $\gcl\in n\mi\cat(\gp)$.
    \end{itemize}
\end{theorem}
\prf We use induction on $n$. The case $n=2$ is Theorem \ref{glob2.the1}.
Suppose the theorem holds for $n-1$ and let $\gcl\in\cagp_S$.\medskip

\emph{Proof of (i) at step $n$}. By induction hypothesis (iii),
$\dcl_\nm\gcl_0^d=\gcl_0^d$ hence
$(\dcl_n\gcl)_0=\dcl_\nm(\ds\ncl\gcl)_0=\dcl_\nm\gcl_0^d=\gcl_0^d$ is discrete;
to show that $\dcl_n\gcl\in\dbb_n$ it remains to check that the Segal maps
\begin{equation*}
\begin{split}
   & (\dcl_n\gcl)_k=\dcl_\nm(\gistr{k})\rw\pro{(\dcl_n\gcl)_1}{(\dcl_n\gcl)_0}{k}= \\
    & =\pro{\dcl_\nm\gcl_1}{\gcl_0^d}{k},
\end{split}
\end{equation*}
for each $k\geq 2$ are weak equivalences in $\dbb_\nm$. Consider the case
$k=2$. We notice that the map
$\eta_2:\tens{\gcl_1}{\gcl_0}\rw\tens{\gcl_1}{\gcl_0^d}$ is a weak equivalence
in $\cagpm$. In fact, passing to the corresponding $(\nm)$-simplicial groups
via the multinerve and taking their diagonals we obtain the commutative diagram
of simplicial groups
\begin{equation*}
        \begin{diagram}
        \diag \ncl\gcl_1\tiund{\diag \ncl\gcl_0}\diag \ncl\gcl_1 & \rTo^{\diag\ncl\eta_2} &
        \diag \ncl\gcl_1\tiund{\diag \ncl\gcl_0^d}\diag \ncl\gcl_1\\
        \dTo^{\nu} && \dTo_{\nu'}\\
        \diag \ncl\gcl_1{\overset{h}{\times}}_{\diag \ncl\gcl_0}\diag \ncl\gcl_1 & \rTo^{\widetilde{\diag\ncl}\eta_2}&
        \diag \ncl\gcl_1{\overset{h}{\times}}_{\diag \ncl\gcl_0^d}\diag \ncl\gcl_1
    \end{diagram}
\end{equation*}
where $\overset{h}{\times}$ denotes the homotopy pullback and $\nu,\nu'$ are
the canonical maps. Since the maps $\diag \ncl\gcl_1\pile{\rw\\ \rw}\diag
\ncl\gcl_0$, being surjective, are fibrations of simplicial groups and the map
$\diag\ncl\gcl_0\rw\diag\ncl\gcl_0^d$ is a weak equivalence (since $\gcl_0$ is
strongly contractible as $\gcl$ is special), we can argue as in the proof of
Theorem \ref{glob2.the1} and conclude that the top map in the above diagram
\begin{equation*}
    \begin{split}
       & \diag \ncl(\tens{\gcl_1}{\gcl_0})=\tens{\diag \ncl\gcl_1}{\diag \ncl\gcl_0}\rw \\
        & \rw\tens{\diag \ncl\gcl_1}{\diag
    \ncl\gcl_0^d}=\diag\ncl(\tens{\gcl_1}{\gcl_0^d})
    \end{split}
\end{equation*}
is a weak\! equivalence \!of simplicial\! groups. \!This means that the map
$\tens{\gcl_1}{\!\gcl_0}\!\rw\tens{\gcl_1}{\gcl_0^d}$ is a weak equivalence in
$\cagpm$ as claimed. Using induction hypothesis (ii) we see that this implies
that the map
$\dcl_\nm(\tens{\gcl_1}{\gcl_0})\rw\dcl_\nm(\tens{\gcl_1}{\gcl_0^d})$ is a weak
equivalence in $\dbb_\nm$. This is because we have weak equivalences
\begin{equation}\label{globgen.eq1ad}
    B\dcl_\nm(\tens{\gcl_1}{\gcl_0})\!=\!B(\tens{\gcl_1}{\gcl_0})\simeq
    B(\tens{\gcl_1}{\gcl_0^d})\!=\!B\dcl_\nm(\tens{\gcl_1}{\gcl_0^d}).
\end{equation}
On the other hand, since $\gcl$ is special, $\gcl_1$ is special (recall Remark
\ref{defspec.rem1} b)), so the diagram
$\gcl_1\supar{d\pt_0}\gcl_0^d\suparle{d\pt_1}\gcl_1$ satisfies the hypotheses
of Lemma \ref{globgen.lem1}. This lemma, together with inductive hypothesis
(iii), gives
\begin{equation}\label{globgen.eq2ad}
    \dcl_\nm(\tens{\gcl_1}{\gcl_0^d})=\tens{\dcl_\nm\gcl_1}{\dcl_\nm\gcl_0^d}=
    \tens{\dcl_\nm\gcl_1}{\gcl_0^d}.
\end{equation}
From (\ref{globgen.eq1ad}) and (\ref{globgen.eq2ad}) we conclude that the Segal
map
\begin{equation*}
    (\dcl_n\gcl)_2=\dcl_\nm(\tens{\gcl_1}{\gcl_0})\rw\tens{\dcl_\nm\gcl_1}{\gcl_0^d}=
    \tens{(\dcl_n\gcl)_1}{(\dcl_n\gcl)_0}
\end{equation*}
is a weak equivalence in $\dbb_\nm$, as required. The case $k>2$ is completely
similar.\medskip

\emph{Proof of (ii) at step $n$}. Let $R=[\dop,\gp]\rw[\dop,\set]$ be as in
Section \ref{lemtec}. An easy computation shows that the following composite
functors from $n$-simplicial groups to simplicial sets coincide (the ``bar"
notation is as in Definition \ref{simtec.def1}):
\begin{equation*}
    \begin{split}
    R_1:& [\dnop,\gp]\supar{\ovl{\ner}}[\dnop,[\dop,\set]]\cong[\dnopiu,\set]\supar{\diag}
    [\dop,\set] \\
     R_2: & [\dnop,\gp]\supar{\diag}[\dop,\gp]\supar{R}[\dop,\set]\\
    R_3: &
    [\dnop,\gp]\cong[\dop,[\dnomen,\gp]]\supar{\overline{\diag}}[\dop,[\dop,\gp]]\supar{\ovl{R}}\\
     & [\dop,[\dop,\set]]\cong[\duop,\set]\supar{\diag}[\dop,\set].
    \end{split}
\end{equation*}
Also note that for $n=2$, $R_3$ is simply given by the composite
\begin{equation*}
    [\duop,\gp]\cong[\dop,[\dop,\gp]]\supar{\ovl{R}}[\dop,[\dop,\set]]
    \cong[\duop,\set]\supar{\diag}[\dop,\set].
\end{equation*}
We want to show that $R\,\diag j_n\dcl_n\gcl=R\,\diag\ncl\gcl$; that is,
$R_2j_n\dcl_n\gcl=R_2\ncl\gcl$. From above, this is the same as showing that
$R_3j_n\dcl_n\gcl=R_3\ncl\gcl$. We are going to compute $R_3j_n\dcl_n\gcl$
using the inductive hypothesis. Recall from Definition \ref{intweak.def1} that
$j_n:\dbb_n\rw[\dnop,\gp]$ is given by the composite
\begin{equation*}
    \dbb_n\hookrightarrow[\dop,\dbb_\nm]\supar{\ovl{j_\nm}}[\dop,[\dnomen,\gp]]\cong[\dnop,\gp].
\end{equation*}
Hence, given $\gcl\in\cagp_S$, $j_n\dcl_n\gcl$ as a simplicial object in
$[\dnomen,\gp]$ takes $[k]\in\dop$ to
$j_\nm(\dcl_n\gcl)_k=j_\nm\dcl_\nm(\ds\ncl\gcl)_k$. It is immediate that
$j_\nm\gcl_0^d=\ncl\gcl_0^d$; hence, from the form of $\ds\ncl\gcl$,
$j_n\dcl_n\gcl$, as a simplicial object in $[\dnomen,\gp]$ is given by
\begin{equation}\label{globgen.eq7}
\begin{split}
   & \cdots j_\nm\dcl_\nm(\tens{\gcl_1}{\gcl_0}\tiund{\gcl_0}\gcl_1)\;\pile{\rTo\\ \rTo\\ \rTo\\
     \rTo}\;j_\nm\dcl_\nm(\tens{\gcl_1}{\gcl_0})\;\pile{\rTo\\ \rTo\\  \rTo} \\
    & \pile{\rTo\\ \rTo\\  \rTo} \;j_\nm\dcl_\nm\gcl_1\;\pile{\rTo\\ \rTo\\
     \lTo}\;\ncl\gcl_0^d.
\end{split}
\end{equation}
By the definition of $R_3$, in order to compute $R_3j_n\dcl_n\gcl$ we need to
apply the composite $R\circ\diag:[\dnomen,\gp]\rw[\dop,\set]$ levelwise in
(\ref{globgen.eq7}) and then take the diagonal of the resulting bisimplicial
set. Thus $R_3j_n\dcl_n\gcl$ is the diagonal of the bisimplicial set
\begin{equation}\label{globgen.eq8}
\begin{split}
   & \cdot\cdot \pile{\rTo\\ \rTo\\  \rTo}
     R\,\diag j_\nm\dcl_\nm(\tens{\gcl_1}{\gcl_0})\pile{\rTo\\ \rTo\\  \rTo}
     \\
    &R\, \diag j_\nm\dcl_\nm\gcl_1 \pile{\rTo\\ \lTo\\ \lTo}R\, \diag\ncl\gcl_0^d.
\end{split}
\end{equation}
By the induction hypothesis, for each $k\geq 2$
\begin{equation*}
    \begin{split}
       & R\, \diag j_\nm\dcl_\nm\gcl_1=R\, \diag \ncl\gcl_1
    \quad\text{and}  \\
        & R\, \diag j_\nm\dcl_\nm(\gistr{k})=R \diag\ncl(\gistr{k})
    \end{split}
\end{equation*}
Hence (\ref{globgen.eq8}) coincides with
\begin{equation}\label{globgen.eq9}
    \cdots R\diag\ncl(\tens{\gcl_1}{\gcl_0})\pile{\rTo\\ \rTo\\  \rTo}
    R\diag\ncl\gcl_1\pile{\rTo\\ \rTo\\  \lTo}R\diag\ncl\gcl_0^d
\end{equation}
so that, in conclusion, $R_3 j_n\dcl_n\gcl$ is the diagonal of
(\ref{globgen.eq9}). On the other hand, from the definition of $R_3$ we see
that $R_3\ncl\gcl$ is the diagonal of the bisimplicial set
\begin{equation}\label{globgen.eq9bis}
   \cdots R\diag\ncl(\tens{\gcl_1}{\gcl_0})\pile{\rTo\\ \rTo\\  \rTo}
    R\diag\ncl\gcl_1\pile{\rTo\\ \rTo\\  \lTo}R\diag\ncl\gcl_0.
\end{equation}
Hence to show $R_3 j_n\dcl_n\gcl=R_3\ncl\gcl$ we need to show that the
diagonals of (\ref{globgen.eq9}) and (\ref{globgen.eq9bis}) coincide. To prove
this, consider the two bisimplicial groups
\begin{equation}\label{globgen.eq10}
    \cdots \diag\ncl(\tens{\gcl_1}{\gcl_0})\pile{\rTo\\ \rTo\\  \rTo}
    \diag\ncl\gcl_1\pile{\rTo\\ \rTo\\ \lTo}\diag\ncl\gcl_0^d
\end{equation}
\begin{equation}\label{globgen.eq11}
    \cdots \diag\ncl(\tens{\gcl_1}{\gcl_0})\pile{\rTo\\ \rTo\\  \rTo}
    \diag\ncl\gcl_1\pile{\rTo\\ \rTo\\ \lTo}\diag\ncl\gcl_0.
\end{equation}
It is easy to see that (\ref{globgen.eq10}) and (\ref{globgen.eq11}) satisfy
the hypotheses of Lemma \ref{weak.lem01}, so that by Lemma \ref{weak.lem01}
applying $R_2$ to (\ref{globgen.eq10}) and (\ref{globgen.eq11}) yields the same
simplicial set. But $R_2=R_3$, hence applying $R_3$ to (\ref{globgen.eq10}) and
(\ref{globgen.eq11}) gives the same simplicial set. From the expression of
$R_3$ for $n=2$ we see that $R_3$ applied to (\ref{globgen.eq10}) and
(\ref{globgen.eq11}) gives respectively the diagonals of
 (\ref{globgen.eq9}) and (\ref{globgen.eq9bis}).
This shows that the diagonals of (\ref{globgen.eq9}) and (\ref{globgen.eq9bis})
coincide. Thus $R_3j_n\dcl_n\gcl=R_3\ncl\gcl$. Since $R_3=R_2$ this means that
$R_2j_n\dcl_n\gcl=R_2\ncl\gcl$ which proves the first part of (ii); since (as
recalled in Section \ref{simtec}) $B=|R_1|$ and $R_1=R_2$ we also have
$B\dcl_n\gcl=B j_n\dcl_n\gcl=|R_1
j_n\dcl_n\gcl|=|R_1\ncl\gcl|=B\ncl\gcl=B\gcl$. This concludes the inductive
step.\medskip

\emph{Proof \,of \,(iii)\, at \,step\, $n$}. Since $\;\gcl\in n\mi\cat(\gp)$,
$\;\gcl_0$\, is \,discrete\, so $\;\gcl_0^d\,=\,\gcl_0$, which implies
$\;\ds\ncl\gcl=\gcl$. Therefore\,
$\dcl_n\gcl\,=\,\ovl{\dcl_\nm}\circ\ds\ncl\gcl\,=\,\ovl{\dcl_\nm}\gcl$, so that
$\;(\dcl_n\gcl)_0\,=\,\dcl_\nm\gcl_0$, $\;(\dcl_n\gcl)_1=\dcl_\nm\gcl_1$,
$(\dcl_n\gcl)_k=\dcl_\nm(\gistr{k})$ for $k\geq 2$. On the other hand, since
$\;\gcl\in n\mi\cat(\gp)$, it follows that $\gcl_0$, $\gcl_1$, $\gistr{k}$ are
all in $(n-1)\mi\cat(\gp)$. Hence by the induction hypothesis
$\dcl_\nm\gcl_0=\gcl_0$, $\dcl_\nm\gcl_1=\gcl_1$,
$\dcl_\nm(\gistr{k})=\gistr{k}$. It follows that $\dcl_n\gcl=\gcl$. \enpr


\section{From internal weak $n$-groupoids to Tamsamani's weak
$(\np)$-groupoids}\label{tamwe} In this section we perform the third main step
in the construction of the comparison functor $F:\cagp\rw\hcl_\np$, which
consists of proving (Theorem \ref{semist.the01}) the existence of a functor
$V_n:\dbb_n\rw\hcl_\np$ from internal weak $n$-groupoids to semistrict
Tamsamani $(n+1)$-groupoids which preserves the homotopy type and sends weak
equivalences to $(\np)$-equivalences.

The construction of $V_n$ is by induction on $n$. The first step of the
induction, $n=2$, is dealt with in Section \ref{tam3}. The formal treatment of
the general case in Section \ref{semist} is preceded by a more informal
overview section (Section \ref{tamgen}) which is especially aimed at clarifying
the logical structure of the inductive argument in Section \ref{semist}. The
three main theorems of the paper (Theorems \ref{specgen.the1},
\ref{globgen.the1} and \ref{semist.the01}) immediately imply (Corollary
\ref{semist.cor01}) the existence of the comparison functor
$F:\cagp_S\rw\hcl_\np$ preserving the homotopy type. The main result of the
paper, which is the semistrictification result of Theorem \ref{result.the1},
follows easily.

\subsection{Geometric versus categorical equivalences}\label{geom}
The following lemma is an essential ingredient in the proof of Theorem
\ref{tam3.the1}, Lemma \ref{semist.lem01}, Proposition \ref{semist.pro01},
leading to one of our main results, Theorem \ref{semist.the01}.
\begin{lemma}\label{geom.lem1}
    A morphism $f:\phi\rw\psi$ in $\tcl_n$ is an $n$-equivalence if and only if $Bf$
    is a weak equivalence.
\end{lemma}
\prf If $f$ is an $n$-equivalence, then $Bf$ is a weak equivalence as proved in
\cite[Proposition 11.2 b]{tam}. We prove the converse by induction on $n$.

The statement is well known to hold for groupoids; that is, for $n=1$ (see for
instance \cite{ms}). In fact, if $\phi$ is a groupoid its nerve $\ner\phi$ is a
fibrant simplicial set. Hence from the expression for the homotopy groups of a
fibrant simplicial set (see for instance \cite{weib}) one can see that
$\pi_0\ner\phi=\phi_0/\!\!\sim$ where, for each $y\in\phi_1$, we have $\pt_0
y\sim\pt_1 y$ and for each $x\in\phi_0$ we have
$\pi_1(\ner\phi,x)=\mathrm{Hom}_\phi(x,x)$ and $\pi_i(\ner\phi,x)=0$ for $i>1$.
Since $\phi$ is a groupoid, for any $x,x'$ such that $\mathrm{Hom}_\phi(x,x')$
is not empty we have $\mathrm{Hom}_\phi(x,x)=\mathrm{Hom}_\phi(x,x')$. Let
$f:\phi\rw\psi$ be a map of groupoids such that $Bf$ is a weak equivalence;
then $\phi_0\bsim\,\cong\,\psi_0\bsim$ so that $f$ is essentially surjective on
objects; also for each $x,x'\in\phi_0$, we have
$\pi_1(\ner\phi,x)\cong\pi_1(\ner\psi,f(x))$ that is, from above,
$\mathrm{Hom}_\phi(x,x')\cong\mathrm{Hom}_\psi(fx,fx')$; so $f$ is fully
faithful. In conclusion $f$ is an equivalence of categories.

Suppose the lemma is true for $n$ and let $f:\phi\rw\psi$ be a morphism in
$\tcl_\np$. By \cite[Proposition 11.4]{tam} there is an $(n+1)$-equivalence
$\alpha:\phi\rw\pp_\np(B\phi)$ natural in $\phi$ and, since $Bf$ is a weak
equivalence, $\pp_\np(Bf)$ is a $(n+1)$-equivalence. We therefore have a
commutative diagram:
\begin{equation}\label{semi.eq01}
\begin{diagram}[h=2.2em]
    \phi & \rTo^{\alpha_\phi} && \pp_\np B\phi \\
    \dTo^{f} &&& \dTo_{\pp_\np Bf}\\
    \psi & \rTo_{\alpha_\psi} && \pp_\np B\psi.
\end{diagram}
\end{equation}
Thus for each $x,y\in\phi_0$ we have a commutative diagram in $\tcl_n$
\begin{equation}\label{semi.eq02}
\begin{diagram}[h=2.2em]
    \phi_{(x,y)} & \rTo^{\alpha_{(x,y)}} && (\pp_\np B\phi)_{(\alpha x,\alpha y)} \\
    \dTo^{f_{(x,y)}} &&& \dTo_{(\pp_\np Bf)_{(\alpha x,\alpha y)}}\\
    \psi_{(fx,fy)} & \rTo_{\alpha_{(fx,fy)}} && (\pp_\np B\psi)_{(\alpha f x,\alpha f y)}
\end{diagram}
\end{equation}
in which $\alpha_{(x,y)},\alpha_{(fx,fy)},(\pp_\np Bf)_{(\alpha x,\alpha y)}$
are $n$-equivalences. Thus in particular $B\alpha_{(x,y)}$,
$B\alpha_{(fx,fy)}$, $B(\pp_\np Bf)_{(\alpha x,\alpha y)}$ are weak
equivalences; the commutativity of (\ref{semi.eq02}) and the 2 out of 3
property of weak equivalences imply that $Bf_{(x,y)}$ is a weak equivalence. By
the induction hypothesis, we conclude that $f_{(x,y)}$ is an $n$-equivalence.
Finally, from (\ref{semi.eq01}) we obtain a commutative diagram in $\tcl_1$
\begin{diagram}
    \tau_1^{(n+1)}\phi && \rTo^{\tau_1^{(n+1)}\alpha_\phi} && \tau_1^{(\np)}\pp_\np(B\phi) \\
    \dTo^{\tau_1^{(n+1)}f} &&&& \dTo_{\tau_1^{(n+1)}\pp_\np Bf}\\
    \tau_1^{(n+1)}\psi && \rTo_{\tau_1^{(n+1)}\alpha_\psi} && \tau_1^{(n+1)}\pp_\np(B\psi) .
\end{diagram}
in which $\tau_1^{(n+1)}\alpha_\phi$, $\tau_1^{(n+1)}\alpha_\psi$,
$\tau_1^{(n+1)}\pp_\np Bf$ are equivalences of groupoids. By the 2 out of 3
property, so is $\tau_1^{(n+1)}f$. This concludes the proof that $f$ is an
$(n+1)$-equivalence.\enpr

\subsection{Semistrict Tamsamani $n$-groupoids}\label{semtam}
We define the category $\hcl_n$ of semi-strict Tamsamani $n$-groupoids. We will
prove in Sections \ref{tam3}, \ref{tamgen}, \ref{semist} that there is a
functor $V_n:\dbb_n\rw\hcl_\np$ preserving the homotopy type and sending weak
equivalences to $(n+1)$-equivalences.

Let $\phi\in\tcl_n$. Recall from Corollary \ref{teclem.cor1} that
$\tcl_n\subset[\dop,\tcl_\nm]$ and that $\tcl_1=\gpd$ is the category of
groupoids. For each $[k]\in\dop$ we denote as usual
$\phi_k=\phi[k]\in\tcl_\nm$. Recall from Section \ref{formal} that an object of
$\tcl_\nm$ is discrete if it is in the image of the functor
$\delta^{(\nm)}:\set\rw\tcl_\nm$. We denote by $\{\cdot\}$ the one-element set.
Since $\{\cdot\}$ is the terminal object in $\set$, $\delta^{(\nm)}\{\cdot\}$
is the terminal object in $\tcl_\nm$. Hence a fibre product in $\tcl_\nm$ over
$\delta^{(\nm)}\{\cdot\}$ coincides with the categorical product.
\begin{definition}\label{tam3.def1}
    The category $\hcl_n$ of semistrict Tamsamani $n$-groupoids is the full
    subcategory of $\tcl_n$ of those objects $\phi$ such that
    \begin{itemize}
      \item[i)] $\phi_0=\delta^{(\nm)}\{\cdot\}$
      \item[ii)] For each $k\geq 0$, the Segal maps
      $\phi_k\rw\fistrr{k}{}=\phi_1\times\overset{k}{\cdots}\times\phi_1$ are
      isomorphisms.
    \end{itemize}
    A map $f$ in $\hcl_n$ is an $n$-equivalence if and only if it is an
    $n$-equivalence in $\tcl_n$.
\end{definition}
The following are straightforward facts about $\hcl_n$.
\begin{remark}\rm\label{tam3.rem1}\text{ }
\begin{itemize}
  \item[a)] Recall from Remark \ref{teclem.rem1}  that there is a full and faithful embedding of
$\tcl_\nm$ in $[\Delta^{{n-2}^{op}},\gpd]$. Under this embedding the discrete
$(\nm)$-groupoid $\delta^{(\nm)}\{\cdot\}$ corresponds to a functor
$\Delta^{{n-2}^{op}}\rw\gpd$ which is constant with value in the trivial
groupoid.
  \item[b)] By definition a map $f:\phi\rw\psi$ in $\hcl_n$ is an
  $n$-equivalence if and only if for each $x,y\in\phi_0$,
  $f_{(x,y)}:\phi_{(x,y)}\rw\psi_{(fx,fy)}$ is a $(n-1)$-equivalence
  and $\tau_1^{(n)}f$ is an equivalence of groupoids.
  Since $\phi_0$ and $\psi_0$ are terminal objects, it is immediate that
  $\phi_{(x,y)}=\phi_1$ and $\psi_{(fx,fy)}=\psi_1$. Hence $f$ is an
  $n$-equivalence if and only if $f_1:\phi_1\rw\psi_1$ is an
  $(\nm)$-equivalence and $\tau^{(n)}_1 f$ is an equivalence of groupoids.
  \item[c)] Let $(\tcl_\nm,\!\times)$ be the category of Tamsamani weak
  $(n-1)$-groupoids equip-ped with the cartesian monoidal structure. We observe
  that $\hcl_n$ is isomorphic to a full subcategory of the category
$\mathrm{Mon}(\tcl_\nm,\times)$ of monoids in $(\tcl_\nm,\times)$. In fact, it
is easy to see that $\mathrm{Mon}(\tcl_\nm,\times)$ is isomorphic to the
category of simplicial objects $\phi$ in $\tcl_\nm$ such that
$\phi_0=\delta^{(\nm)}\{\cdot\}$ and the Segal maps are isomorphisms. On the
other hand, from its definition it is immediate that the category $\hcl_n$ is
the full subcategory of simplicial objects $\phi$ in $\tcl_\nm$ with
$\phi_0=\delta^{(\nm)}\{\cdot\}$ and Segal maps isomorphisms such that
$\tau^{(n)}_1\phi$ is a groupoid.
\end{itemize}
\end{remark}

\subsection{Delooping from a group-based structure to a set-based structure}\label{deloop}
Recall that, by Remark \ref{intweak.rem1}, $\dbb_n$ is a full subcategory
\begin{equation}\label{deloop.eq1}
    \dbb_n\subset\underset{\nm}{\underbrace{[\dop,[\dop,\ldots,[\dop}},\cgp]\cdots]
\end{equation}
and, by Remark \ref{teclem.rem1}, $\tcl_n$ is a full subcategory
\begin{equation}\label{deloop.eq2}
    \tcl_n\subset\underset{\nm}{\underbrace{[\dop,[\dop,\ldots,[\dop}},\gpd]\cdots].
\end{equation}
The\, right\, sides\, of\; (\ref{deloop.eq1})\; and\; (\ref{deloop.eq2})\;
are\, isomorphic\, to\; $[\dnomen,\cgp]$ and $[\dnomen,\gpd]$ respectively, but
for their use later in Sections \ref{totam} and \ref{semist} it is clearer to
keep them in the form (\ref{deloop.eq1}) and (\ref{deloop.eq2}). In this
section we define for each $n\geq 2$ a functor
\begin{equation*}
    V_n:\underset{\nm}{\underbrace{[\dop,[\dop,\ldots,[\dop}},\cgp]\cdots]
    \rw\underset{n}{\underbrace{[\dop,[\dop,\ldots,[\dop}},\gpd]\cdots].
\end{equation*}
This functor will be used in Sections \ref{tam3}, \ref{tamgen}, \ref{semist} to
perform the passage from the category $\dbb_n$ of internal weak $n$-groupoids
to the category $\hcl_\np$ of semistrict \tam $(n+1)$-groupoids. More
precisely, we will show in \ref{semist} that when restricted to the full
subcategory $\dbb_n$ we obtain a functor $V_n:\dbb_n\rw\hcl_\np$.

Informally, the functor $V_n$ is easily described as follows. Recall that a
group can be considered as a category with just one object, hence there is a
nerve functor from the category of groups to the category of simplicial sets.
Given an object $\phi$  of
$\underset{\nm}{\underbrace{[\dop,[\dop,\ldots,[\dop}},\cgp]\cdots]$, to
construct $V_n\phi$ we first take the nerve of each group occurring in the
diagram defining $\phi$. Next we regard the resulting multi-simplicial
structure as a simplicial object using as simplicial direction the ``delooping"
direction; that is, the direction along which we took the nerve of each group.

It is essential for its use in the subsequent sections that this definition of
the functor $V_n$ is made precise and that the appropriate notation is in
place. Since every internal category in groups is an internal groupoid, the
forgetful functor $U:\gp\rw\set$ gives rise to a functor $U_1:\cgp\rw\gpd$ as
follows. If $\gcl\in\cgp$ is given by
\begin{equation*}
    G_1\tiund{G_0} G_1\rTo G_1 \;\pile{\rTo\\ \rTo\\ \lTo}\; G_0
\end{equation*}
then $U_1 \gcl$ is the groupoid obtained by taking the underlying diagram of
sets; that is,
\begin{equation*}
    UG_1\tiund{UG_0}UG_1\rTo UG_1\;\pile{\rTo\\ \rTo\\ \lTo}\; UG_0.
\end{equation*}
We call $U_1$ the \emph{underlying groupoid functor}. For each $n\geq 2$ we
define a functor
\begin{equation*}
    U_n:\underset{\nm}{\underbrace{[\dop,[\dop,\ldots,[\dop}},\cgp]\cdots]\rw
    \underset{\nm}{\underbrace{[\dop,[\dop,\ldots,[\dop}},\gpd]\cdots].
\end{equation*}
Informally, this functor takes the underlying groupoid of every cat$^1$-group
occurring in each object of its domain. Formally $U_n$ is defined inductively
as follows. For $n=2$ the functor $U_2:[\dop,\cgp]\rw[\dop,\gpd]$ is obtained
by applying $U_1$ levelwise; that is, with our usual notation as in Definition
\ref{simtec.def1}, $U_2=\ovl{U_1}$. Thus if $\phi\in[\dop,\cgp]$ and we denote,
as usual, $\phi_k=\phi([k])\in\cgp$, we have $(U_2\phi)_k=U_1\phi_k\in\gpd$.
Suppose, inductively, that we have defined $U_\nm$. Then we set
$U_n=\ovl{U_\nm}$.

We now define $V_n$. In what follows the functor
$\delta^{(n)}:\set\!\rw\!\underset{\nm}{\underbrace{[\dop,[\dop,\ldots,[\dop}},$
$\gpd]\cdots]$ is as in Section \ref{formal}. That is,
$\delta^{(1)}:\set\rw\gpd$ takes a set $X$ to the discrete groupoid on $X$.
Inductively, $\delta^{(n)}=d\delta^{(\nm)}$ where
$d:\tcl_\nm\rw[\dop,\tcl_\nm]$ takes $\psi\in\tcl_\nm$ to the constant
simplicial object in $\tcl_\nm$, $(d\psi)_k=\psi_k$ for all $[k]\in\dop$. Given
$\phi\in\underset{\nm}{\underbrace{[\dop,[\dop,\ldots,[\dop}},\cgp]\cdots]$
then
$V_n\phi\in\underset{n}{\underbrace{[\dop,[\dop,\ldots,[\dop}},\gpd]\cdots]$
takes each $[r]\in\dop$ to
$(V_n\phi)_r=(V_n\phi)[r]\in\underset{\nm}{\underbrace{[\dop,[\dop,\ldots,[\dop}},\gpd]\cdots]$
given by
\begin{equation}\label{deloop.eq3bis}
    (V_n\phi)_r=
  \begin{cases}
    \delta^{(n)}\{\cdot\}, & {r=0;} \\
    U_n\phi, & {r=1;} \\
     U_n\phi\times\overset{r}{\cdots}\times U_n\phi, & { r > 1}.
  \end{cases}
\end{equation}
The faces and degeneracy operators in the simplicial object $V_n\phi$ are those
given by the nerve construction. More precisely, the nerve functor
$\ner:\gp\rw[\dop,\set]$ induces a functor $N:\cgp\rw[\dop,\gpd]$ obtained by
taking the nerves of the group of objects, groups of arrows, groups of
composable arrows of the internal category. Explicitly, for each $p\geq 0$ and
$\gcl\in\cgp$, denoting by $*$ the trivial groupoid, we have
\begin{equation}\label{deloop.eq4}
    (N\gcl)_p=
      \begin{cases}
        *, & { p=0;} \\
        U_1\gcl, & { p=1;} \\
        U_1\gcl\times\overset{p}{\cdots}\times U_1\gcl, & { p > 1.}
      \end{cases}
\end{equation}
Here is a sketch of $N\gcl\in[\dop,\gpd]$.
\begin{diagram}[small,h=2em]
    \cdots & \pile{\rTo\\ \rTo\\ \rTo\\ \rTo}&
    \scriptstyle(UG_1\tiund{UG_0}UG_1)\times(UG_1\tiund{UG_0}UG_1)& \pile{\rTo\\ \rTo\\ \rTo} &
    \scriptstyle UG_1\tiund{UG_0}UG_1 & \pile{\rTo\\ \rTo\\ \lTo} & \{\cdot\}\\
    \dTo && \dTo && \dTo && \dTo\\
    \cdots \scriptstyle UG_1\times UG_1\times UG_1 & \pile{\rTo\\ \rTo\\ \rTo\\ \rTo} & \scriptstyle UG_1\times UG_1
    & \pile{\rTo\\ \rTo\\ \rTo} & \scriptstyle UG_1 & \pile{\rTo\\ \rTo\\ \lTo} & \{\cdot\}\\
    \dTo \dTo \uTo && \dTo \dTo \uTo && \dTo \dTo \uTo && \dTo \dTo \uTo\\
    \scriptstyle UG_0\times UG_0\times UG_0 & \pile{\rTo\\ \rTo\\ \rTo\\ \rTo} &\scriptstyle UG_0\times
     UG_0 & \pile{\rTo\\ \rTo\\ \rTo} & \scriptstyle UG_0 & \pile{\rTo\\ \rTo\\ \lTo} &
     \{\cdot\}\\
     \cdots p=3 && p=2 && p=1 && p=0.
\end{diagram}
On the other hand, from the definition of $U_n$ we know that, for each
$[p_1]...[p_n]\in\dop$,
$(U_n\phi)([p_1])([p_2])\ldots([p_\nm])=U_1(\phi([p_1])([p_2])\ldots([p_\nm]))$.
Hence the face and degeneracy operators in (\ref{deloop.eq4}) give the
corresponding ones in $V_n\phi$. We finally illustrate the case $n=2$ with a
picture. Given $\phi\in[\dop,\cgp]$ and denoting $\phi_k=\phi([k])\in\cgp$ as
above, $V_2\phi\in[\dop,[\dop,\gpd]]$ is given by
\begin{diagram}[h=1.8em]
    \vdots && \vdots && \vdots && \vdots & \vdots\\
    \cdots \scriptstyle U_1\phi_2\times U_1\phi_2\times U_1\phi_2& \pile{\rTo\\ \rTo\\ \rTo}&
    \scriptstyle U_1\phi_2\times U_1\phi_2 & \pile{\rTo\\ \rTo\\ \rTo} &
    \scriptstyle U_1\phi_2 & \pile{\rTo\\ \rTo\\ \lTo} & * & q=2 \\
    \dTo\dTo\dTo && \dTo\dTo\dTo && \dTo\dTo\dTo && \dTo\dTo\dTo\\
    \cdots \scriptstyle U_1\phi_1\times U_1\phi_1\times U_1\phi_1& \pile{\rTo\\ \rTo\\ \rTo}&
    \scriptstyle U_1\phi_1\times U_1\phi_1 & \pile{\rTo\\ \rTo\\ \rTo} &
    \scriptstyle U_1\phi_1 & \pile{\rTo\\ \rTo\\ \lTo} & * & q=1 \\
    \dTo\dTo\uTo && \dTo\dTo\uTo && \dTo\dTo\uTo && \dTo\dTo\uTo\\
    \cdots \scriptstyle U_1\phi_0\times U_1\phi_0\times U_1\phi_0& \pile{\rTo\\ \rTo\\ \rTo}&
    \scriptstyle U_1\phi_0\times U_1\phi_0 & \pile{\rTo\\ \rTo\\ \rTo} &
    \scriptstyle U_1\phi_0 & \pile{\rTo\\ \rTo\\ \lTo} & * & q=0 \\
    \cdots \quad r=3 & & r=2 && r=1 && r=0 & \uTo_q\\
    &&&&&\text{ }&\lTo_r & \text{ }
\end{diagram}
where $V_2\phi$ is a simplicial object in $[\dop,\gpd]$ along the delooping
direction $r$; that is, the direction along which we take the nerve of each
group.

\subsection{From internal weak 2-groupoids to semistrict Tamsamani 3-group-oids}\label{tam3}
In this section we are going to perform the passage from the category $\dbb_2$
of internal weak 2-groupoids to the category $\hcl_3$ of semistrict Tamsamani
3-groupoids. This is the first step in the inductive argument we will give in
Section \ref{tamgen} to construct a functor $V_n:\dbb_n\rw\hcl_{\np}$ for
general $n$ which preserves the homotopy type and sends weak equivalences to
$(n+1)$-equivalences. This low dimensional case has also been treated in detail
in \cite{paol}. We include it here for completeness and also because we use a
slightly different notation from \cite{paol}. We first give an overview of the
main ideas. Throughout this section, the ``bar" notation is as in Definition
\ref{simtec.def1}.

Let $V_2:[\dop,\cgp]\rw[\dop,[\dop,\gpd]]$ be as in (\ref{deloop.eq3bis}). As
noticed in Remarks \ref{teclem.rem1} and \ref{intweak.rem1} respectively,
$\dbb_2$ is a full subcategory of $[\dop,\cgp]$ and $\tcl_3$ is a full
subcategory of $[\dop,[\dop,\gpd]]$. In the next theorem we show that if
$\phi\in\dbb_2$, then $V_2\phi\in\hcl_3$. In other words $V_2$ restricts to a
functor $V_2:\dbb_2\rw\hcl_3$. We will also see that this functor preserves the
homotopy type and sends weak equivalences to 3-equivalences.

Notice that the semistrictness of $V_2\phi$ arises immediately by the way in
which $V_2$ has been constructed: the Segal maps in the ``delooping" direction
along which we took the nerve of each group are automatically isomorphisms;
that is, for each $k\geq 2$
$(V_2\phi)_k\cong(V_2\phi)_1\tiund{(V_2\phi)_0}\cdots\tiund{(V_2\phi)_0}(V_2\phi)_1$
and $(V_2\phi)_0=\delta^{(2)}\{\cdot\}$. These are precisely the two conditions
defining the full subcategory $\hcl_3$ of $\tcl_3$. What remains to be proved
to show that $V_2\phi\in\hcl_3$ is that $U_2\phi\in\tcl_2$ and $\tau^{(3)}_1
V_2\phi$ is a groupoid.

Recall from Section \ref{deloop} that $U_2=\ovl{U_1}$, so that given
$\phi\in\dbb_2\subset[\dop,\cagp]$, $U_2\phi\in[\dop,\gpd]$ is the simplicial
object
\begin{equation*}
    \cdots U_1\phi_2\pile{\rTo\\ \rTo\\ \rTo} U_1\phi_1\pile{\rTo\\ \rTo\\
    \lTo} U_1\phi_0.
\end{equation*}
Here we have denoted, as usual, $\phi_k=\phi[k]\in\cgp$, for each $[k]\in\dop$.
There are two main steps involved in the proof that $U_2\phi\in\tcl_2$. One is
that the weak equivalences of $\cgp$ $\phi_k\rw\fistr{k}$, which hold as
$\phi\in\dbb_2$, give rise to weak equivalences of the underlying groupoids
$U_1\phi_k\rw U_1\phi_1\tiund{U_1\phi_0}\overset{k}{\cdots}$
$\tiund{U_1\phi_0}U_1\phi_1$. As observed in Lemma \ref{geom.lem1} (case
$n=1$), in the category of groupoids weak equivalences coincide with
categorical equivalences. Hence the Segal maps in $U_2\phi\in[\dop,\gpd]$ are
equivalences of categories.

The second main step in the proof that $U_2\phi\in\tcl_2$ is the computation of
$\tau^{(1)}_{1}U_2\phi$, which has to be a groupoid for $U_2\phi$ to be an
object of $\tcl_2$. The strategy is to show that $\tau^{(2)}_1U_1\phi$ is the
underlying groupoid of an object of $\cgp$. This computation involves the
functor $\ovl{\pi_0}:[\dop,\cgp]\rw[\dop,\gp]$ induced by $\pi_0:\cgp\rw\gp$
where the latter is as in Lemma \ref{catn.lem1}. One observes that since, by
Lemma \ref{catn.lem1}, $\pi_0$ preserves fibre products over discrete objects
and sends weak equivalences to isomorphisms, $\ovl{\pi_0}$ restricts to a
functor $\ovl{\pi_0}:\dbb_2\rw\ner(\cgp)$ from $\dbb_2$ into the full
subcategory $\ner(\cgp)$ of those simplicial groups which are nerves of objects
of $\cgp$. Composing $\ovl{\pi_0}$ with the isomorphism
$\eta:\ner\cgp\cong\cgp$, this yields a functor
\begin{equation*}
    T^{(1)}:\dbb_2\rw\cgp.
\end{equation*}
The way to think about $\tiu$ is that it plays for $\dbb_2$ the same role that
the functor $\tau^{(2)}_1:\tcl_2\rw\gpd$ has for $\tcl_2$. Recall in fact that
$\tau^{(2)}_1$ is the composite of the restriction to $\tcl_2$ of the functor
$\ovl{\tau^{(1)}_0}:[\dop,\cat]\rw[\dop,\set]$ induced by
$\tau^{(1)}_0:\cat\rw\set$ with the isomorphism $\nu:\ner(\cat)\cong\cat$. The
functor $\pi_0:\cgp\rw\gp$ is an internal analogue of
$\tau^{(1)}_0:\cat\rw\set$. It turns out that $\tau^{(2)}_1$ and $\tiu$ are
related as follows:
\begin{equation}\label{tam3.eq1}
    \tau^{(2)}_1 U_2\phi=U_1\tiu\phi
\end{equation}
for every $\phi\in\dbb_2$. Thus, for each $\phi\in\dbb_2$, (\ref{tam3.eq1})
exhibits $\tau^{(2)}_1U_2\phi$ as the underlying  groupoid of an object of
$\cgp$.

To prove that $V_2\phi\in\tcl_3$ it will remain to show that
$\tau^{(3)}_{1}V_2\phi$ is the nerve of a groupoid. This will be an easy
computation using the above functors and the definition of $\tau^{(3)}_{1}$.
Finally, the fact that $V_2$ preserves the homotopy type will be a
straightforward consequence of the notion of classifying space of a
multi-simplicial group. With the use of Lemma \ref{geom.lem1} this also will
imply easily that $V_2$ sends weak equivalences to 3-equivalences. We now give
the formal statement and proof of the theorem.
\begin{theorem}\label{tam3.the1}
    Let $V_2:[\dop,\cgp]\rw[\dop,[\dop,\gpd]]$ be as in (\ref{deloop.eq3bis}).
    Then $V_2$ restricts to a
    functor $V_2:\dbb_2\rw\hcl_3$. Further, for each $\phi\in\dbb_2$, $B\phi=B V_2
    \phi$ and $V_2$ sends weak equivalences in $\dbb_2$ to 3-equivalences in
$   \hcl_3$.
\end{theorem}
\prf Given $\phi\in\dbb_2\subset[\dop,\cgp]$  we  denote as usual
$\phi_k=\phi([k])\in\cgp$ for all $[k]\in\dop$. By definition of $\hcl_3$, to
show that $V_2\phi\in\hcl_3$ we need to show that
\begin{itemize}
  \item[(i)] $(V_2\phi)_0=\delta^{(2)}\{\cdot\}$.
  \item[(ii)] For each $k\geq 2$ the Segal maps
  $(V_2\phi)_k\rw(V_2\phi)_1\tiund{(V_2\phi_0)}\cdots\tiund{(V_2\phi)_0}(V_2\phi)_1$
  are isomorphisms.
  \item[(iii)] $U_2\phi\in\tcl_2$.
  \item[(iv)] $\tau^{(3)}_{1}V_2\phi$ is a groupoid.
\end{itemize}\medskip

\emph{Proof of (i) and (ii)}. These hold automatically from the definition of
$V_2$.

\medskip

\emph{Proof of (iii)}. Recall from Section \ref{deloop} that $U_2=\ovl{U_1}$,
so $U_2\phi\in[\dop,\gpd]$ is the simplicial object taking $[k]\in\dop$ to
$(U_2\phi)_k=U_1\phi_k$. By definition of $\tcl_2$, to show that
$U_2\phi\in\tcl_2$ we need to check that
\begin{itemize}
  \item[a)] $U_1\phi_0$ is discrete.
  \item[b)] For all $k\geq 2$, the Segal maps $U_1\phi_k\rw\fistrr{}{U_1}$
are equivalences of categories.
  \item[c)] For all $x,y\in U_1\phi_0$, $(U_2\phi)_{(x,y)}$ is a groupoid.
  \item[d)] $\tau^{(2)}_{1}U_2\phi$ is the nerve of a groupoid.
\end{itemize}
Notice that a) is immediate since, as $\phi\in\dbb_2$, $\phi_0$ is a discrete
internal category in groups. Also c) is immediate since $U_1\phi_1$ is a
groupoid. To prove b), recall that since $\phi\in\dbb_2$, for all $k\geq 0$ the
Segal maps $\phi_k\rw\fistr{k}$ are\, weak\, equivalences\, in\, $\cgp$. \,That
is,\, by\, Lemma\; \ref{catn.lem1},\; if $\ner:\cgp\rw[\dop,\gp]$ is the nerve
functor, the maps of simplicial groups $\ner\phi_k\rw\ner(\fistr{k})$ induce
isomorphisms of homotopy groups for all $i\geq 0$
\begin{equation}\label{tam3.eq1bis}
    \pi_i\ner\phi_k\cong\pi_i\ner(\fistr{k}).
\end{equation}

But the homotopy groups of a simplicial group are the homotopy groups of the
underlying (fibrant) simplicial set (see for instance \cite{weib}). The
underlying simplicial set of the simplicial group $\ner\phi_k$ is $N
U_1\phi_k$, where $N:\gpd\rw[\dop,\set]$ is the nerve functor and similarly for
$\ner(\fistr{k})$. It follows from (\ref{tam3.eq1bis}) that for each $i\geq 0$
there are isomorphisms
\begin{equation*}
    \pi_i N U_1\phi_k\cong\pi_i N(\fistrr{k}{U_1}).
\end{equation*}
Since the homotopy groups of a fibrant simplicial set coincide with the
homotopy groups of its geometric realization \cite{weib}, this means that the
maps of groupoids $U_1\phi_k\!\rw\!\fistrr{k}{U_1}$ are weak equivalences; that
is, they induce weak equivalences of classifying spaces. By Lemma
\ref{geom.lem1} (case $n=1$) these maps are also categorical equivalences. This
completes the proof of b).

To prove d), we are going to show that $\tau^{(2)}_{1}U_2\phi$ is the
underlying groupoid of an object of $\cgp$. Let $\pi_0:\cgp\rw\gp$ be as in
Lemma \ref{catn.lem1} and let $\ovl{\pi_0}:[\dop,\cgp]\rw[\dop,\gp]$ be induced
by $\pi_0$ (as in Definition \ref{simtec.def1}). As recalled in Lemma
\ref{catn.lem1} $\pi_0$ preserves fibre products over discrete objects and
sends weak equivalences to isomorphisms. Hence, if $\phi\in\dbb_2$, for all
$k\geq 2$ we have
\begin{equation*}
    \begin{split}
       &(\ovl{\pi_0}\phi)_k=\pi_0\phi_k\cong\pi_0(\fistr{k})\cong   \\
    &\cong\fistrr{k}{\pi_0}=(\ovl{\pi_0}\phi)_1\tiund{(\ovl{\pi_0}\phi)_0}
    \overset{k}{\cdots}\tiund{(\ovl{\pi_0}\phi)_0}
    (\ovl{\pi_0}\phi)_1.
     \end{split}
\end{equation*}
Hence $\ovl{\pi_0}\phi$ is the nerve of an object of $\cgp$. In other words,
$\ovl{\pi_0}$ restricts to a functor $\ovl{\pi_0}:\dbb_2\rw\ner(\cgp)$ from
$\dbb_2$ to the full subcategory $\ner(\cgp)$ of those simplicial groups which
are nerves of objects of $\cgp$. Composing $\ovl{\pi_0}$ with the isomorphism
$\eta:\ner\cgp\cong\cgp$ we obtain a functor
\begin{equation*}
    \tiu:\dbb_2\rw\cgp
\end{equation*}
with
\begin{equation*}
    \tiu=\eta\ovl{\pi_0}_{|_{\dbb_2}},\quad
    \ner\circ\tiu=\ovl{\pi_0}_{|_{\dbb_2}}.
\end{equation*}
We claim that, for all $\phi\in\dbb_2$
\begin{equation}\label{tam3.eq2}
    \tau^{(2)}_1 U_2\phi=U_1 \tiu\phi.
\end{equation}
To prove this, observe that if $U:\gp\rw\set$ is the forgetful functor
\begin{equation}\label{tam3.eq3}
    U\pi_0=\tau^{(1)}_0 U_1
\end{equation}
and there is an obvious commutative diagram
\begin{equation}\label{tam3.eq3bis}
    \xymatrix{
    [\dop,\gp]  \ar[rr]^{\ovl{U}} &&  [\dop,\set]\\
    \ner\cgp \ar@{^{(}->}[u] \ar[rr]^{\ovl{U}} \ar[d]^\eta && \ner(\gpd)\ar@{^{(}->}[u]\ar[d]^\nu\\
    \cgp\ar[rr]^{U_1} && \gpd}
\end{equation}
where $\ovl{U}$ is induced by $U$. Recalling that, by definition (see
(\ref{formal.eq1}) in Section \ref{formal}),
$\tau^{(2)}_1=\nu\ovl{\tau^{(1)}_0}$, $\tiu=\eta\ovl{\pi_0}_{|_{\dbb_2}}$ and
$U_2=\ovl{U_1}$, we obtain from (\ref{tam3.eq3}) and (\ref{tam3.eq3bis}):
\begin{equation*}
     \tau^{(2)}_{1}U_2\phi=\nu\ovl{\tau^{(1)}_0}\ovl{U_1}\phi=\nu\ovl{\tau^{(1)}_0U_1}\phi
    =\nu\ovl{U\pi_0}\phi=\nu\ovl{U}\ovl{\pi_0}\phi=U_1\eta\ovl{\pi_0}\phi=U_1\tiu\phi
\end{equation*}
which is (\ref{tam3.eq2}). Hence for any $\phi\in\dbb_2$, $\tau^{(2)}_1U_2\phi$
is the underlying groupoid of the internal category in groups $\tiu\phi$. This
concludes the proof of d).

\emph{Proof of iv)}. Recalling from (\ref{formal.eq1}) in Section \ref{formal}
that $\ner\tau^{(3)}_1=\ovl{\tau^{(2)}_0}$ and
$\tau^{(2)}_0=\tau^{(1)}_0\tau^{(2)}_1$ we obtain, for all $k\geq 0$,
\begin{equation*}
    (\ner\tau^{(3)}_1V_2\phi)_k=(\ovl{\tau^{(2)}_0}V_2\phi)_k=
    \tau^{(2)}_0(V_2\phi)_k=\tau^{(1)}_0\tau^{(2)}_1(V_2\phi)_k.
\end{equation*}
Using (\ref{tam3.eq2}), (\ref{tam3.eq3}), the expression of $(V_2\phi)_k$ and
the fact, noted in Section \ref{formal}, that $\tau^{(2)}_1$ preserves fibre
products over discrete objects we therefore calculate, for all $k>1$,
\begin{equation*}
    \begin{split}
       (\ner\tau^{(3)}_1V_2\phi)_0 & =\{\cdot\} \\
       (\ner\tau^{(3)}_1V_2\phi)_1 & =\tau^{(1)}_0\tau^{(2)}_1 U_2\phi
       =\tau^{(1)}_0 U_1 \tiu\phi=U\pi_0 \tiu\phi\\
       (\ner\tau^{(3)}_1V_2\phi)_k & =\tau^{(1)}_0\tau^{(2)}_1
       (\pro{U_2\phi}{}{k})=\pro{\tau^{(1)}_0\tau^{(2)}_1U_2\phi}{}{k}=\\
          &=\pro{U\pi_0 \tiu\phi}{}{k}.
     \end{split}
\end{equation*}
This shows that $\ner\tau^{(3)}_1V_2\phi$ is the nerve of the group
$\pi_0\tiu\phi$. Thus $\tau^{(3)}_1V_2\phi$ is a group, so in particular a
groupoid, as required.

To complete the proof of the theorem we need to show that $V_2$ preserves the
homotopy type and that it sends weak equivalences to 3-equivalences. The
classifying space of an object $\phi\in\dbb_2$ is by definition the classifying
space of the bisimplicial group $j_2\phi$, where the embedding
$j_2:\dbb_2\rw[\Delta^{2^{op}},\gp]$ is as in Definition \ref{intweak.def1}. As
recalled in Section \ref{simtec}, to compute the classifying space of a
multi-simplicial group we take the nerve of the group in each dimension and
then compute the classifying space of the resulting multi-simplicial set.
Formally, we have the composite functor $B:\dbb_2\rw\mathrm{Top}$ given by
\begin{equation*}
    \dbb_2\supar{j_2}[\Delta^{2^{op}},\gp]\supar{\ovl{\ner}}[\Delta^{2^{op}},[\dop,\set]]\cong
    [\Delta^{3^{op}},\set]\supar{diag}[\dop,\set]\supar{|\,\cdot\,|}\mathrm{Top}
\end{equation*}
where $\ovl{\ner}$ is induced by the nerve functor $\ner:\gp\rw[\dop,\set]$,
$diag$ is the multi-diagonal and $|\,\cdot\,|$ is the geometric realization. On
the other hand, recall  that the classifying space of any object of $\tcl_3$
(hence in particular of $\hcl_3$) is obtained by taking the nerve of each
groupoid in the embedding $\tcl_3\subset[\dop,\tcl_2]\subset[\dop,[\dop,\gpd]]$
and then taking the classifying space of the resulting 3-simplicial set.
Formally, we have the composite functor $B:\tcl_3\rw\mathrm{Top}$ given by
\begin{equation*}
\begin{split}
  &\tcl_3\hookrightarrow[\dop,\tcl_2]\hookrightarrow[\dop,[\dop,\gpd]]\supar{\ovl{\ner}}
  [\dop,[\dop,[\dop,\set]]]\cong    \\
    & \cong [\Delta^{3^{op}},\set]\supar{diag}[\dop,\set]\supar{|\,\cdot\,|}\mathrm{Top}.
\end{split}
\end{equation*}
Since the functor $V_2:\dbb_2\rw\hcl_3$ is obtained by taking the nerve of the
group at each dimension it is clear that the two functors
\begin{equation*}
    \dbb_2\supar{j_2}[\Delta^{2^{op}},\gp]\supar{\ovl{\ner}}[\Delta^{2^{op}},
    [\dop,\set]]\cong[\Delta^{3^{op}},\set]
\end{equation*}
and
\begin{equation*}
    \dbb_2\supar{V_2}\hcl_3\hookrightarrow[\dop,[\dop,\gpd]]\supar{\ovl{\ner}}[\dop,[\dop,[\dop,\set]]]\cong
    [\Delta^{3^{op}},\set]
\end{equation*}
when applied to the same $\phi\in\dbb_2$ yield 3-simplicial sets which can
differ at most by a permutation of the simplicial coordinates, and which
therefore have the same diagonal. Hence $B\phi=BV\phi$ for all $\phi\in\dbb_2$.

Finally, suppose that $f:\phi\rw\psi$ is a weak equivalence in $\dbb_2$; then
$BV_2\phi=B\phi\simeq B\psi=BV_2\psi$ so that $BV_2f$ is a weak equivalence. By
Lemma \ref{geom.lem1}, $V_2f$ is therefore a 3-equivalence.\enpr\bigskip

In Section \ref{semist} we will construct by induction the functor
$V_n:\dbb_n\rw\hcl_\np$ for each $n\geq 2$. In proving the inductive step we
will need the following corollary. A discussion of why this will be needed in
the construction of $V_3$ is contained in the overview section \ref{tamgen}.
\begin{corollary}\label{tam3.cor1}
    The functor $\tiu:\dbb_2\rw\cgp$ preserves weak equivalences and
    fibre products over discrete objects and sends discrete objects to discrete objects.
\end{corollary}
\prf Recall that $\ner\circ \tiu=\ovl{\pi_0}_{|_{\dbb_2}}$ and that, by Lemma
\ref{catn.lem1}, $\pi_0:\cgp\rw\gp$ preserves fibre products over discrete
objects. For each $\phi\in\dbb_2$ denote as usual $\phi_k=\phi([k])\in\cgp$,
$[k]\in\dop$. Consider the fibre product in $\dbb_2$ $\phi\tiund{\psi}\phi$,
where $\psi$ is discrete. Then, for each $k\geq 0$, $\psi_k$ is a discrete
object of $\cgp$ so we have
\begin{equation*}
    \begin{split}
      &(\ner \tiu(\phi\tiund{\psi}\phi))_k=(\ovl{\pi_0}(\phi\tiund{\psi}\phi))_k=
      \pi_0(\phi\tiund{\psi}\phi)_k=\pi_0(\phi_k\tiund{\psi_k}\phi_k)= \\
        &=\pi_0\phi_k\tiund{\pi_0\psi_k}\pi_0\phi_k=
        (\ovl{\pi_0}\phi)_k\tiund{(\ovl{\pi_0}\psi)_k}(\ovl{\pi_0}\phi)_k=\\
        &= (\ner \tiu\phi)_k\tiund{(\ner \tiu\psi)_k}(\ner \tiu\phi)_k.
    \end{split}
\end{equation*}
We conclude that
\begin{equation*}
    \ner \tiu(\phi\tiund{\psi}\phi)=\ner \tiu\phi\tiund{\ner \tiu\psi}\ner \tiu\phi
    =\ner (\tiu\phi\tiund{\tiu\psi} \tiu\phi),
\end{equation*}
and therefore $\tiu(\phi\tiund{\psi}\phi)=\tiu\phi\tiund{\tiu\psi} \tiu\phi$.
That is, $\tiu$ preserves fibre products over discrete objects.

It is clear that if $\psi$ is a discrete object of $\dbb_2$ then $\tiu\psi$ is
a discrete object of $\cgp$. It remains to prove that $\tiu$ preserves weak
equivalences. Let $f:\phi\rw\psi$ be a weak equivalence in $\dbb_2$. By Theorem
\ref{tam3.the1} $V_2 f$  is a 3-equivalence in $\hcl_3$. Hence (see Remark
\ref{tam3.rem1} b)) $(V_2f)_1=U_2 f$ is a 2-equivalence in $\tcl_2$. By
definition this implies that $\tau^{(2)}_1U_2f$ is an equivalence of groupoids.
But recall from (\ref{tam3.eq2}) that $\tau^{(2)}_1U_2\phi=U_1 \tiu\phi$, for
each $\phi\in\dbb_2$; therefore we conclude that
\begin{equation}\label{tam3.eq4}
    \tau^{(2)}_1U_2f=U_1 \tiu f:U_1 \tiu\phi\rw U_1 \tiu\psi
\end{equation}
is a weak equivalence. This easily implies that the map of $\cgp$ $\tiu
f:\tiu\phi\rw\tiu\psi$ is a weak equivalence. In fact, if
$N:\gpd\rw[\dop,\set]$ is the nerve functor, the weak equivalence
(\ref{tam3.eq4}) gives isomorphisms of homotopy groups
\begin{equation}\label{tam3.eq5}
    \pi_iNU_1\tiu\phi\cong\pi_iNU_1\tiu\psi
\end{equation}
for all $i\geq 0$.  If $\ner:\cgp\rw[\dop,\gp]$ is the nerve functor, arguing
as in the proof of Theorem \ref{tam3.the1} (proof of (iii)), the underlying
simplicial set of the simplicial group $\ner \tiu\phi$ is $NU_1\tiu\phi$, and
similarly for $\psi$. But the homotopy groups of a simplicial group are the
homotopy groups of the underlying fibrant simplicial set. Hence
(\ref{tam3.eq5}) yields isomorphisms for all $i\geq 0$, $\pi_i\ner
\tiu\phi\cong\pi_i\ner \tiu\psi$. By Lemma \ref{catn.lem1} this shows that
$\tiu f$ is a weak equivalence in $\cgp$.\enpr

\subsection{From internal weak $n$-groupoids to Tamsamani model:\label{totam}
overview of the general case}\label{tamgen} In this section we provide an
overview of the most important steps involved in the main theorem of this
chapter, Theorem \ref{semist.the01}, which establishes the passage, for any
$n$, from the category $\dbb_n$ of internal weak $n$-groupoids to the category
$\hcl_\np$ of semistrict Tamsamani $(\np)$-groupoids. Let
\begin{equation*}
    V_n:\underset{\nm}{\underbrace{[\dop,[\dop,\ldots,[\dop}},\cgp]\cdots]
    \rw\underset{n}{\underbrace{[\dop,[\dop,\ldots,[\dop}},\gpd]\cdots].
\end{equation*}
be as in (\ref{deloop.eq3bis}). Recall from Remark \ref{teclem.rem1} and Remark
\ref{intweak.rem1} that $\dbb_n$ and $\tcl_\np$ are full subcategories
respectively of the domain and codomain of $V_n$. Our aim is to show that $V_n$
restricts to a functor $V_n:\dbb_n\rw\hcl_\np$. We will also see that this
functor preserves the homotopy type and sends weak equivalences to
$(\np)$-equivalences.

By definition of $\hcl_\np$ to show that $V_n\phi\in\hcl_\np$ when
$\phi\in\dbb_n$ we need to show that:
\begin{itemize}
  \item[(i)] $(V_n\phi)_0=\delta^{(n)}\{\cdot\}$.
  \item[(ii)] For each $k\geq 0$, the Segal maps $(V_n\phi)_k\rw\pro{(V_n\phi)_1}
  {(V_n\phi)_0}{}$ are isomorphisms.
  \item[(iii)] $U_n\phi\in\tcl_n$.
  \item[(iv)] $\tau_1^{(\np)}V_n\phi$ is a groupoid.
\end{itemize}
Notice that (i) and (ii) hold automatically from the definition of $V_n$: the
Segal maps in the delooping direction along which we took the nerve of each
group are automatically isomorphisms.

We shall now illustrate the main points in the proof. We are going to give an
inductive argument to show that $V_n\phi\in\hcl_\np$ when $\phi\in\dbb_n$ and
that $V_n:\dbb_n\rw\hcl_\np$ has the desired properties. However, this will not
be a simple induction starting from the hypothesis that there is
$V_\nm:\dbb_\nm\rw\hcl_n$. The argument is more complex and we are going to
explain below why this needs to be the case. We know by Section \ref{tam3} that
$V_2\phi\in\hcl_3$ when $\phi\in\dbb_2$, so in particular (iii) above holds for
$n=2$. We could start inductively by supposing that (iii) holds for $(n-1)$
(that is, $U_\nm\psi\in\tcl_\nm$ for $\psi\in\dbb_\nm$), and try to show that
it holds for $n$. As we are going to see, this will help us to identify what is
the correct inductive hypothesis we need to assume. In what follows the ``bar"
notation is as in Definition \ref{simtec.def1}.

Suppose that for each $\psi\in\dbb_\nm$, $U_\nm\psi\in\tcl_\nm$. Recall from
Section \ref{deloop} that $U_n=\ovl{U_\nm}$. Thus if
$\phi\in\dbb_n\subset[\dop,\dbb_\nm]$ and for each $[k]\in\dop$ we denote as
usual $\phi_k=\phi([k])\in\dbb_\nm$ we have $(U_n\phi)_k=U_\nm\phi_k$; hence,
by the induction hypothesis, $U_n\phi\in[\dop,\tcl_\nm]$. Notice also that
$U_\nm\phi_0$ is discrete, as $\phi_0$ is discrete (since $\phi\in\dbb_n$). So
what remains to be proved for $U_n\phi$ to be an object of $\tcl_n$ is that:
\label{p1}
\begin{itemize}
  \item[a)] For all $k\geq 0$, the Segal maps $U_\nm\phi_k\rw\fistrr{k}{U_\nm}$
  are $(\nm)$-equivalences.
  \item[b)] $\tau^{(n)}_1 U_n\phi$ is a groupoid.
\end{itemize}
The proof of a) is very similar to the case $n=2$ and relies on Lemma
\ref{geom.lem1} which guarantees that, in the category $\tcl_\nm$, the
$(n-1)$-equivalences and the weak equivalences coincides. More precisely, since
$\phi\in\dbb_n$ there is a weak equivalence $\phi_k\rw\fistr{k}$ in $\dbb_\nm$,
for each $k\geq 2$. It is not hard to show that this gives rise to a weak
equivalence in $\tcl_\nm$:
\begin{equation*}
    U_\nm\phi_k\rw U_\nm(\fistr{k})\cong\fistrr{k}{U_\nm}.
\end{equation*}
By Lemma \ref{geom.lem1}, these are also $(n-1)$-equivalences, proving a).

The strategy to prove b) is, as in the case $n=2$, to show that
$\tau^{(n)}_1U_n\phi$ is the underlying groupoid of an object of $\cgp$. The
argument for general  $n$ however involves an issue which does not appear in
the case $n=2$. To understand why this is the case, let us consider for
instance the case $n=3$ and try to compute $\tau^{(3)}_1U_3\phi$, for each
$\phi\in\dbb_3$. Using Lemma \ref{teclem.lem1} and the fact that
$U_3=\ovl{U}_2$ we compute
\begin{equation}\label{tamgen.eq1}
    \tau^{(3)}_1U_3\phi=\tau^{(2)}_1\ovl{\tau^{(2)}_1}\ovl{U}_2\phi=
    \tau^{(2)}_1\ovl{\tau^{(2)}_1U_2}\phi.
\end{equation}
Let $\tiu:\dbb_2\rw\cgp$ be as in the proof of Theorem \ref{tam3.the1} and let
$\ovl{\tiu}:[\,\dop,\;\dbb_2\,]\;\rw\;[\;\dop,\;\cgp\;]$ be induced by $\tiu$.
Recalling from (\ref{tam3.eq1}) in Section \ref{tam3} that
$\tau^{(2)}_1U_2\psi=U_1\tiu\psi$ for each $\psi\in\dbb_2$ and that
$U_2=\ovl{U_1}$, we calculate, for each $k\geq 0$
\begin{equation}\label{tamgen.eq2}
    (\ovl{\tau^{(2)}_1U_2}\phi)_k=\tau^{(2)}_1U_2\phi_k=U_1
    \tiu\phi_k=U_1(\ovl{\tiu}\phi)_k=(U_2\ovl{\tiu}\phi)_k.
\end{equation}
Hence $\ovl{\tau^{(2)}_1U_2}\phi=U_2\ovl{\tiu}\phi$ and therefore from
(\ref{tamgen.eq1})
\begin{equation}\label{tamgen.eq3}
    \tau^{(3)}_1U_3\phi=\tau^{(2)}_1U_2\overline{\tiu}\phi.
\end{equation}
Recall (see (\ref{tam3.eq2}) in proof of Theorem \ref{tam3.the1}) that
$\tau^{(2)}_1 U_2\psi=U_1\tiu\psi$ for each $\psi\in\dbb_2$. We would like to
use this fact in (\ref{tamgen.eq3}) taking $\psi=\ovl{\tiu}\phi$ to rewrite the
right hand side as $U_1\tiu\ovl{\tiu}\phi$; this would exhibit $\tau^{(3)}_1
U_3\phi$ as the underlying groupoid of the cat$^1$-group $\tiu\ovl{\tiu}\phi$.
However, the relation $\tau^{(2)}_1 U_2\psi=U_1\tiu\psi$ holds when
$\psi\in\dbb_2$ so, in order to apply this relation to the case where
$\psi=\ovl{\tiu}\phi$, we first need to show that $\ovl{\tiu}\phi$ is in
$\dbb_2$. This is the new issue appearing when $n>2$. The fact that
$\ovl{\tiu}\phi\in\dbb_2$ is a consequence of the fact, proved in Corollary
\ref{tam3.cor1}, that $\tiu:\dbb_2\rw\cgp$ preserves weak equivalences, fibre
products over discrete objects, and sends a discrete object to a discrete
object. In fact, if $\phi\in\dbb_3$, then $\phi_0$ is discrete, so
$(\ovl{\tiu}\phi)_0=\tiu\phi_0$ is discrete. Also, the Segal maps
$\phi_k\rw\fistr{k}$ are weak equivalences, therefore the maps
\begin{equation*}
    \begin{split}
      (\ovl{\tiu}\phi)_k & =\tiu\phi_k\rw\;  \tiu(\fistr{k})\cong\fistrr{k}{\tiu}= \\
        & =\pro{(\ovl{\tiu}\phi)_1}{(\ovl{\tiu}\phi)_0}{k}
    \end{split}
\end{equation*}
are weak equivalences in $\cgp$. By definition this shows that
$\ovl{\tiu}\phi\in\dbb_2$. We can then apply the relation
$\tau^{(2)}_1U_2\psi=U_1\tiu\psi$ for $\psi\in\dbb_2$ in (\ref{tamgen.eq3})
taking for $\psi$ the object $\ovl{\tiu}\phi\in\dbb_2$ to obtain
\begin{equation}\label{tamgen.eq4}
    \tau^{(3)}_1U_3\phi=U_1\tiu\ovl{\tiu}\phi.
\end{equation}
This exhibits $\tau^{(3)}_1U_3\phi$ as the underlying groupoid of the
cat$^1$-group $\tiu\ovl{\tiu}\phi$.

Let us now denote by $\tid:\dbb_3\rw\cgp$ the composite functor
$\tid=\tiu\ovl{\tiu}$. Then (\ref{tamgen.eq4}) immediately gives
$\tau^{(3)}_1U_3\phi=U_1\tid\phi$ for all $\phi\in\dbb_3$. One can ask if
$\tid$ satisfies the same property of $\tiu$ of preserving weak equivalences
and fibre products over discrete objects and of sending discrete objects to
discrete objects. This is indeed the case and, as is easily guessed, this
property of $\tid$ is needed for the next step up, $n=4$, in the proof that
$\tau^{(4)}_1U_3\phi$ is a groupoid when $\phi\in\dbb_4$.

Recall from the proof of Corollary \ref{tam3.cor1} that the similar properties
of $\tiu$ arise from having constructed the functor $V_2:\dbb_2\rw\hcl_3$
sending weak equivalences to 3-equivalences. Likewise, the analogous properties
of $\tid$ arise from having the functor $V_3:\dbb_3\rw\hcl_4$ sending weak
equivalences to 4-equivalences. We have already given the most important steps
in the construction of $V_3$. Given $\phi\in\dbb_3$ what remains to be proved
to show that $V_3\phi\in\hcl_4$ is the fact that $\tau^{(4)}_1V_3\phi$ is a
groupoid: this is an easy computation formally analogous to the case $n=2$,
which shows that $\tau^{(4)}_1V_3\phi$ is in fact the group $\pi_0\tid\phi$.
The fact that $V_3$ sends weak equivalences to 4-equivalences is also
straightforward with the use of Lemma \ref{geom.lem1}, and formally analogous
to the case $n=2$.

This discussion indicates how the inductive step works for general $n$. Namely,
at step $(\nm)$ we will have the following situation (inductive hypothesis):
for all $\psi \in\dbb_\nm$, \label{ref1}
\begin{itemize}
  \item[(i)] $U_\nm\psi\in\tcl_\nm$.
  \item[(ii)]There exists $T^{(n-2)}\!:\!\dbb_\nm\!\rw\!\cgp$ such that
    $\tau^{(\nm)}U_\nm\psi\!=\!U_1T^{(n-2)}\psi$.
  \item[(iii)] $T^{(n-2)}$ preserves weak equivalences and fibre products over
  discrete objects and sends discrete objects to discrete objects.
\end{itemize}
Let us now consider $\phi\in\dbb_n$. We want to show that given the inductive
hypothesis, we have $U_n\phi\in\tcl_n$; that is, a) and b) are satisfied. As
explained at the beginning of the section, using the inductive hypothesis (i),
condition a) holds, relative to the Segal maps for $U_n\phi$. As for condition
b), using inductive hypothesis (ii), a computation formally analogous to the
one we made above for $n=3$ yields
\begin{equation}\label{tamgen.eq5}
    \tau_1^{(n)}U_n\phi=\tau^{(2)}_1U_2\ovl{T^{(n-2)}}\phi
\end{equation}
where $\ovl{T^{(n-2)}}:[\dop,\dbb_\nm]\rw[\dop,\cgp]$ is induced by
$T^{(n-2)}$. We will then use induction hypothesis (iii) to show that if
$\phi\in\dbb_n$ then $\ovl{T^{(n-2)}}\phi$ is an object of $\dbb_2$: the
argument is formally analogous to the one we gave for $n=3$. Hence using the
relation $\tau^{(2)}_1U_2\psi=U_1\tiu\psi$ for all $\psi\in\dbb_2$, which holds
as it is the first step of the induction, and taking
$\psi=\ovl{T^{(n-2)}}\phi\in\dbb_2$ we rewrite (\ref{tamgen.eq5}) as
\begin{equation}\label{tamgen.eq6}
    \tau_1^{(n)}U_n\phi=U_1 \tiu\ovl{T^{(n-2)}}\phi.
\end{equation}
This exhibits $\tau^{(n)}_1\! U_n \phi$ as the underlying groupoid of the
cat$^1$-group $\tiu\ovl{T^{(n-2)}}\phi$, proving b); in conclusion a) and b)
are proved at step $n$, so $U_n\phi\in\tcl_n$, which is (i) at step $n$. Now
let $T^{(\nm)}$ be the composite functor
\begin{equation*}
    T^{(\nm)}=\tiu\ovl{T^{(n-2)}}:\dbb_n\rw\cgp.
\end{equation*}
Then (\ref{tamgen.eq6}) gives $\tau_1^{(n)}U_n\phi=U_1T^{(\nm)}\phi$ for all
$\phi\in\dbb_3$. This proves (ii) at step $n$. It remains to prove (iii) at
step $n$. This will be a consequence of the fact that, having proved (i) and
(ii) at step $n$, one can show there is a functor $V_n:\dbb_n\rw\hcl_\np$
sending weak equivalences to $(n+1)$-equivalences. The argument is formally
analogous to the one we sketched for $n=3$. We conclude with a schematic
summary of the inductive argument (in what follows (i), (ii), (iii) refer to
the properties listed on page \pageref{ref1}):\medskip

\noindent\emph{\textbf{First step of the induction: $n=2$}}
\begin{itemize}
  \item[a)] (i) and (ii) hold for $n=2$ (proved in Section \ref{tam3})\; $\Rightarrow$
  \item[b)] $V_2:\dbb_2\rw\hcl_3$ preserving the homotopy type $\Rightarrow$
  \item[c)] (iii) holds for $n=2$.
\end{itemize}
\emph{\textbf{Inductive step}}
\begin{itemize}
  \item[a)] (inductive hypothesis) (i), (ii), (iii) hold for $(n-1)$ $\Rightarrow$
  \item[b)] (i) and (ii) hold for $n$ $\Rightarrow$
  \item[c)] $V_n:\dbb_n\rw\hcl_\np$ preserving the homotopy type $\Rightarrow$
  \item[d)] (iii) holds for $n$.
\end{itemize}
We reiterate the following point, which is important to correctly understand
the argument. Notice that the information we are interested in is that for each
$n$ there is a functor $V_n:\dbb_n\rw\hcl_\np$ which preserves the homotopy
type and sends weak equivalences to $(n+1)$-equivalences. For each $n$ this is
a consequence of having (i) and (ii) alone at step $n$. However, we also need
condition (iii) in order to prove inductively that (i) and (ii) hold for each
$n$. In fact, except at the first step of the induction (where we proved
directly that (i) and (ii) hold for $n=2$), having (i) and (ii) at step $(n-1)$
is not sufficient to prove (i) and (ii) at step $n$, and inductive hypothesis
(iii) is needed. This is because, in proving (i) at step $n$, we need to show
that $\tau_1^{(n)}U_n\phi$ is a groupoid. The key point is that in the
computation of $\tau_1^{(n)}U_n\phi$ for $\phi\in\dbb_n$ we find
$\tau_1^{(n)}U_n\phi=\tau^{(2)}_1U_2\ovl{T^{(n-2)}}\phi$. We need to know that
$\ovl{T^{(n-1)}}\phi\in\dbb_2$ in order to apply the relation
$\tau_1^{(2)}U_2\psi=U_1\tiu\psi$ taking $\psi=\ovl{T^{(n-2)}}\phi$ and obtain
$\tau^{(2)}_1U_2\ovl{T^{(n-1)}}\phi=U_1\tiu\ovl{T^{(n-1)}}\phi$, which proves
that $\tau_1^{(n)}U_n\phi$ is a groupoid. The fact that $\ovl{T^{(n-1)}}\phi$
is an object of $\dbb_2$ is a consequence of the inductive hypothesis (iii).

In the next section we will first show (Lemma \ref{semist.lem01}) the
implication that for each $n\geq 2$ conditions (i) and (ii) at step $n$ give
the functor $V_n:\dbb_n\rw\hcl_\np$ with the desired properties. We will then
give an inductive argument in Proposition \ref{semist.pro01} to verify that
(i), (ii) and (iii) hold for each $n$. We will then conclude (Theorem
\ref{semist.the01}) the existence for each $n\geq 2$ of the functor
$V_n:\dbb_n\rw\hcl_\np$ with the desired properties.


\subsection{From internal weak $n$-groupoids to Tamsamani model: general case}\label{semist}
This section is devoted to the formal statement and proof of one of our main
theorems, asserting the existence for every $n\geq 2$ of a functor
$V_n:\dbb_n\rw\hcl_\np$ from internal weak $n$-groupoids to semistrict
Tamsamani $(n+1)$-groupoids preserving the homotopy type. Throughout this
section, for each $n\geq 2$, let
\begin{equation*}
    V_n:\underset{\nm}{\underbrace{[\dop,[\dop,\ldots,[\dop}},\cgp]\cdots]
    \rw\underset{n}{\underbrace{[\dop,[\dop,\ldots,[\dop}},\gpd]\cdots].
\end{equation*}
be as in (\ref{deloop.eq3bis}). Throughout this section, the ``bar" notation is
as in Definition \ref{simtec.def1}.
\begin{lemma}\label{semist.lem01}
    Let $V_n$ be as in (\ref{deloop.eq3bis}) and suppose that the following
conditions are satisfied:
\begin{itemize}
  \item[i)] $U_n\phi\in\tcl_n$ for each $\phi\in\dbb_n$.
  \item[ii)] There exists a functor $T^{(\nm)}:\dbb_n\rw\cgp$ such that
$\tau^{(n)}_1U_n\phi=U_1T^{(\nm)}\phi$ for each $\phi\in\dbb_n$.
\end{itemize}
Then $V_n$ restricts to a functor $V_n:\dbb_n\rw\hcl_\np$ which preserves the
homotopy type and which sends weak equivalences to $(n+1)$-equivalences.
\end{lemma}
\prf By hypothesis (i), and the definition of $V_n$,
$\,V_n\phi\in[\dop,\tcl_n]$. Also, by definition of $V_n$, for each $k\geq 0$,
$(V_n\phi)_k=\fistpar{k}{V_n}=(V_n\phi)_1\overset{k}{\times\cdots\times}(V_n\phi)_1$
where the last equality holds as $(V_n\phi)_0=\delta^{(n)}\{\cdot\}$ is the
terminal object. Hence $V_n\phi\in\hcl_\np$ provided we show that
$\ta{\np}{1}V_n\phi$ is a groupoid. Recall from Section \ref{formal} that
$\ner\ta{\np}{1}=\ovl{\ta{n}{0}}$ and $\ta{n}{0}=\ta{1}{0}\ta{n}{1}$. Hence,
for all $k\geq 0$,
\begin{equation*}
    (\ner\ta{\np}{1}V_n\phi)_k=(\ovl{\ta{n}{0}}V_n\phi)_k=\ta{n}{0}
    (V_n\phi)_k=\ta{1}{0}\ta{n}{1}(V_n\phi)_k.
\end{equation*}
Thus, using hypothesis (ii), the expression of $V_n$, the fact (see
(\ref{tam3.eq3})) that $U\pi_0=\ta{1}{0}U_1$ and the fact, noticed in Section
\ref{formal}, that $\tau^{(n)}_1$ preserves fibre products over discrete
objects, we calculate for each $k\geq 1$,
\begin{equation*}
    \begin{split}
      (\ner\ta{\np}{1}V_n\phi)_0 & = \{\cdot\} \\
       (\ner\ta{\np}{1}V_n\phi)_1  &
    =\ta{1}{0}\ta{n}{1}U_n\phi=\ta{1}{0}U_1T^{(\nm)}\phi=U\pi_0T^{(\nm)}\phi\\
    (\ner\ta{\np}{1}V_n\phi)_k &
    =\ta{1}{0}\ta{n}{1}(U_n\phi\overset{k}{\times...\times}U_n\phi)=
    \ta{1}{0}\ta{n}{1}U_n\phi\overset{k}{\times...\times}\ta{1}{0}\ta{n}{1}U_n\phi\\
    &=U\pi_0T^{(\nm)}\phi\overset{k}{\times...\times}U\pi_0T^{(\nm)}\phi.
    \end{split}
\end{equation*}
This shows that $\,\ner\,\ta{\np}{1}\,V_n\,\phi$ is the nerve of the group
$\pi_0T^{(\nm)}\phi$. Thus $\ta{\np}{1}V_n\phi$ is a group, so in particular is
a groupoid, as required.

We are now going to show that $V_n$ preserves the homotopy type. There are full
and faithful maps (denoted with the same symbol for convenience):
\begin{equation}\label{semist.eq12}
    \begin{split}
        & \ncl:\dbb_n\rw[\Delta^{\np^{op}},\set] \\
        & \ncl:\tcl_\np\rw[\Delta^{\np^{op}},\set].
    \end{split}
\end{equation}
Each of these is obtained by composing the fully faithful embeddings
$\dbb_n\rw[\dnomen,$ $\cgp]$ and $\tcl_\np\rw[\dnop,\gpd]$ of Remarks
\ref{intweak.rem1} and \ref{teclem.rem1} with the functors
$[\dnomen,\cgp]\rw[\Delta^{\np^{op}},\set]$ and
$[\dnop,\gpd]\rw[\Delta^{\np^{op}},\set]$ induced by the nerve functors. Recall
from Section \ref{deloop} that the functor $V_n:\dbb_n\rw\hcl_\np$ was obtained
by taking the nerve of each group in the diagram defining $\dbb_n$. It follows
that, for any $\phi\in\dbb_n$, the $(n+1)$-simplicial sets $\ncl\phi$ and $\ncl
V_n\phi$ differ at most by a permutation of the multi-simplicial coordinates.
Hence $\ncl\phi$ and $\ncl V_n\phi$ have the same diagonal:
\begin{equation}\label{semist.eq02}
    diag\,\ncl\phi=diag\,\ncl V_n\phi.
\end{equation}
On the other hand the classifying space functors are given by the composites
\begin{equation*}
    \begin{split}
      & B:  \dbb_n\supar{\ncl}[\Delta^{\np^{op}},\set]\supar{diag}[\dop,\set]\supar{|\cdot|}\tp \\
      & B:
      \tcl_\np\supar{\ncl}[\Delta^{\np^{op}},\set]\supar{diag}[\dop,\set]\supar{|\cdot|}\tp.
    \end{split}
\end{equation*}
It follows from (\ref{semist.eq02}) that $B\phi=BV_n\phi$. Finally we are going
to show that $V_n$ sends weak equivalences to $(\np)$-equivalences. Suppose
that $f:\phi\rw\psi$ is a weak equivalence in $\dbb_n$. Then, from above
$BV_n\phi=B\phi\simeq B\psi=BV_n\psi$ so that $BV_n f$ is a weak equivalence.
By Lemma \ref{geom.lem1}, $V_n f$ is also a $(\np)$-equivalence.\enpr
\begin{proposition}\label{semist.pro01}
    Let $V_n$ be as in (\ref{deloop.eq3bis}), $n\geq 2$. For each $\phi\in\dbb_n$, the
    following hold:
    \begin{itemize}
      \item[(i)] $U_n\phi\in\tcl_n$.
      \item[(ii)] There exists a functor $T^{(\nm)}:\dbb_n\rw\cgp$ such
      that $\ta{n}{1}U_n\phi=U_1T^{(\nm)}\phi$.
      \item[(iii)] $T^{(\nm)}$ preserves weak equivalences and fibre products
      over discrete objects and sends each discrete object to a discrete object.
    \end{itemize}
\end{proposition}
\prf We use induction on $n$. For $n=2$ properties (i), (ii), (iii) have been
proved in \ref{tam3} (see Theorem \ref{tam3.the1}, (\ref{tam3.eq2}) and
Corollary \ref{tam3.cor1}). Suppose, inductively, that (i), (ii), (iii) hold
for $n-1$ and let $\phi\in\dbb_n$. Denote, as usual,
$\phi_k=\phi([k])\in\dbb_\nm$ for all $k\geq 0$.\medskip

\emph{Proof of (i) at step $n$}. Recall from Section \ref{deloop} that
$U_n=\ovl{U_\nm}$; that is, $(U_n\phi)_k=U_\nm\phi_k$ for all $k\geq 0$. By
induction hypothesis (i) , $U_\nm\phi_k\in\tcl_\nm$, so that
$U_n\phi\in[\dop,\tcl_\nm]$. Also, $U_\nm\phi_0$ is discrete as $\phi_0$ is
discrete (since $\phi\in\dbb_n$). What remains to be proved for $U_n\phi$ to be
an object of $\tcl_n$ is that:
\begin{itemize}
  \item[a)] For all $k\geq 2$ the Segal maps $U_\nm\phi_k\rw\fistrr{k}{U_\nm}$
  are $(\nm)$-equivalences
  \item[b)] $\ta{n}{1}U_n\phi$ is a groupoid.
\end{itemize}
To prove a) recall that, since $\phi\in\dbb_n$, for each $k\geq 2$ the Segal
maps $\phi_k\rw\fistr{k}$ are equivalences in $\dbb_\nm$; that is, they induce
weak equivalences of classifying spaces. Our aim is to show that this implies
that the maps of $(\nm)$-groupoids $U_\nm\phi_k\rw\fistrr{k}{U_\nm}$ are weak
equivalences, from which a) will follows by Lemma \ref{geom.lem1}. Recall that
the classifying space of an object $\psi$ of $\dbb_\nm$ is the classifying
space of the $(\nm)$-simplicial group $j_\nm\psi$ where $j_\nm$ is as in
Definition \ref{intweak.def1}. Hence the map of $(\nm)$-simplicial groups
\begin{equation}\label{semist.eq04}
    j_\nm\phi_k\rw j_\nm(\fistr{k})
\end{equation}
is a weak equivalence. We notice that this implies that the map of the
underlying $(\nm)$-simplicial sets
\begin{equation}\label{semist.eq05}
    \ovu j_\nm\phi_k\rw\ovu j_\nm(\fistr{k})
\end{equation}
is a weak equivalence, where
$\ovu:[\Delta^{{n-1}^{op}},\gp]\rw[\Delta^{{n-1}^{op}},\set]$ is induced by the
forgetful functor $U:\gp\rw\set$. To see this, we just pass to the diagonal
simplicial groups and argue similarly to the case $n=2$. That is, applying the
diagonal functor $diag:[\dnop,\gp]\rw[\dop,\gp]$ to (\ref{semist.eq04}), we
find a weak equivalence of simplicial groups $diag\,j_\nm\phi_k\rw
diag\,j_\nm(\fistr{k})$. Arguing as in the proof of the case $n=2$, since the
homotopy groups of a simplicial group are the homotopy groups of the underlying
fibrant simplicial set, this implies that there is a weak equivalence of the
underlying simplicial sets $\ovu\, diag\,j_\nm\phi_k\rw\ovu\,
diag\,j_\nm(\fistr{k})$; here we have denoted with the same symbol $\ovu$ the
functor $\ovu:[\dop,\gp]\rw[\dop,\set]$. The obvious commutative diagram
\begin{diagram}[w=5em,h=2.4em]
    {[\dnop,\gp]}& \rTo^{\ovu} & [\dnop,\set]\\
    \dTo_{diag} && \dTo_{diag}\\
    {[\dnop,\gp]}& \rTo^{\ovu} & [\dnop,\set]
\end{diagram}
therefore implies that the maps of simplicial sets
\begin{equation}\label{semist.eq06}
    diag\,\ovu j_\nm\phi_k\rw diag\,\ovu j_\nm(\fistr{k})
\end{equation}
are weak equivalences. On the other hand, from the definition of $U_\nm$ (see
Section \ref{deloop}) it is easy to see that there is a commutative diagram
\begin{equation*}
    \xy
    0;/r.18pc/:
    (-50,25)*+{\scriptstyle\dbb_\nm}="1";
    (30,25)*+{\scriptstyle[\Delta^{\nm^{op}},\gp]}="2";
    (-50,-25)*+{\scriptstyle\underset{n-2}{\underbrace{\scriptstyle[\dop,[\dop,\ldots,[\dop}},\gpd]\cdots]}="3";
    (30,-25)*+{\scriptstyle\underset{n-1}{\underbrace{\scriptstyle[\dop,[\dop,\ldots,[\dop}},\set]\cdots]\quad
    \overset{\xi}{\cong}}="4";
    (-50,0)*+{\scriptstyle\underset{n-2}{\underbrace{\scriptstyle[\dop,[\dop,\ldots,[\dop}},\cgp]\cdots]}="5";
    (70,-25)*+{\underset{}{\scriptstyle[\Delta^{\nm^{op}},\set]}}="6";
    {\ar^{j_\nm}"1";"2"};
    {\ar^{\ovu}"2";"6"};
    {\ar_{U_\nm}"5";"3"};
    {\ar^{\scriptstyle\ovl{\ner}}"3";"4"};
    {\ar@{_{(}->}^{}"1";"5"};
\endxy
\end{equation*}
where $\ovl{\ner}$ is induced by the nerve functor $\ner:\gpd\rw[\dop,\set]$
and $\xi$ is a canonical bijection. Hence (\ref{semist.eq06}) gives a weak
equivalence of simplicial sets
\begin{equation}\label{semist.eq07}
    diag\,\xi\ovl{\ner}U_\nm\phi_k\rw
    diag\,\xi\ovl{\ner}U_\nm\phi_k(\fistr{k}).
\end{equation}
But notice that
\begin{equation}\label{semist.eq13}
    \diag\xi\ovl{\ner}U_\nm\phi_k=\diag\ncl U_\nm\phi_k
\end{equation}
where $\ncl:\tcl_\nm\rw[\Delta^{\nm^{op}},\set]$ is the full and faithful map
(\ref{semist.eq12}), and similarly for the right hand side of
(\ref{semist.eq07}). Hence from (\ref{semist.eq07}) and (\ref{semist.eq13}) we
obtain weak equivalences
\begin{equation*}
    \begin{split}
      BU_\nm\phi_k & =|\,diag\,\ncl U_\nm\phi_k|\simeq |\,diag\,\ncl U_\nm(\fistr{k})|= \\
        & =B U_\nm(\fistr{k})=B(\fistrr{k}{U_\nm}).
    \end{split}
\end{equation*}
By Lemma \ref{geom.lem1}, the maps of $(n-1)$-groupoids
\begin{equation*}
    U_\nm\phi_k\rw\fistrr{}{U_\nm},
\end{equation*}
being weak equivalences, are also $(n-1)$-equivalences. This proves a).

We are now going to prove b). Using Lemma \ref{teclem.lem1} and the fact (see
Section \ref{deloop}) that $U_n=\ovu_\nm$, we compute
\begin{equation}\label{semist.eq08}
    \ta{n}{1}U_n\phi=\ta{2}{1}\ovl{\ta{\nm}{1}}\ovu_\nm\phi=\ta{2}{1}\ovl{\ta{\nm}{1}U_\nm}\,\phi.
\end{equation}
By induction hypothesis (ii), $\ta{\nm}{1}U_\nm=U_1T^{(n-2)}$. Hence for each
$k\geq 0$ we have
\begin{equation}\label{semist.eq09}
    (\ovl{\ta{\nm}{1}U_\nm}\phi)_k=\ta{\nm}{1}U_\nm\phi_k=U_1T^{(n-2)}\phi_k=
    U_1(\ovl{T^{(n-2)}}\phi)_k=(U_2\ovl{T^{(n-2)}}\phi)_k
\end{equation}
where $\ovl{T^{(n-2)}}:[\dop,\dbb_\nm]\rw[\dop,\cgp]$ is induced by
$T^{(n-2)}$. Hence from (\ref{semist.eq09}) $\ovl{\tau^{(n-1)}_1
U_\nm}\,\phi=U_2\ovl{T^{(n-2)}}\phi$ so that from (\ref{semist.eq08}) we
conclude
\begin{equation}\label{semist.eq10}
    \ta{n}{1}U_n\phi=\ta{2}{1}U_2\ovl{T^{(n-2)}}\phi.
\end{equation}
We now claim that $\ovl{T^{(n-2)}}\phi\in\dbb_2$. This is a consequence of the
induction hypothesis (iii) that $T^{(n-2)}$ preserves weak equivalences and
fibre products over discrete objects and sends discrete objects to discrete
objects. In fact, since $\phi\in\dbb_n$, $\phi_0$ is discrete so
$(\ovl{T^{(n-2)}}\phi)_0=T^{(n-2)}\phi_0$ is discrete. Further, for each $k\geq
0$, the Segal maps
\begin{equation*}
\begin{split}
    (\ovl{T^{(n-2)}}\phi)_k & =T^{(n-2)}\phi_k\rw T^{(n-2)}(\fistr{k})\cong \\
    & \cong\fistrr{}{T^{(n-2)}}=\\
    & =\fistpar{}{\ovl{T^{(n-2)}}}
\end{split}
\end{equation*}
are weak equivalences. This shows that $\ovl{T^{(n-2)}}\phi\in\dbb_2$, as
claimed.

From the first step of the induction $(n=2)$ we know that
$\ta{2}{1}U_2\psi=U_1\tiu\psi$ for all $\psi\in\dbb_2$. Applying this to
(\ref{semist.eq10}), taking $\psi=\ovl{T^{(n-2)}}\phi\in\dbb_2$, we deduce that
\begin{equation}\label{semist.eq11}
    \ta{n}{1}U_n\phi=U_1\tiu\ovl{T^{(n-2)}}\phi.
\end{equation}
This shows that $\ta{n}{1}U_n\phi$ is a groupoid, proving b).\medskip

\emph{Proof of (ii) at step $n$}. Let $T^{(\nm)}:\dbb_n\rw\cgp$ be defined by
the composite $T^{(\nm)}=\tiu \ovl{T^{(n-2)}}$. By (\ref{semist.eq11}),
$\ta{n}{1}U_n\phi=U_1 T^{(n-1)}\phi$ for every $\phi\in\dbb_n$, which shows
(ii) at step $n$.\medskip

\emph{Proof of (iii) at step $n$}. Consider the fibre product
$\phi\tiund{\psi}\phi$ in $\dbb_n$ with $\psi$ discrete. We want to show that
$T^{(n-1)}$ preserves fibre products over discrete objects. By the induction
hypothesis, $T^{(n-2)}$ preserves fibre products over discrete objects. This
easily implies that $\ovl{T^{(n-2)}}$ preserves fibre products over discrete
objects. In fact, since $\psi$ is discrete, for each $k\geq 0$,
$\psi_k\in\dbb_\nm$ is discrete and we have
\begin{equation*}
    \begin{split}
      &(\ovl{T^{(n-2)}}(\phi\tiund{\psi}\phi))_k =T^{(n-2)}(\phi\tiund{\psi}\phi)_k
       =T^{(n-2)}(\phi_k\tiund{\psi_k}\phi_k)=\\
        & = T^{(n-2)}\phi_k\tiund{T^{(n-2)}\psi_k}T^{(n-2)}\phi_k =
        (\ovl{T^{(n-2)}}\phi)_k\tiund{(\ovl{T^{(n-2)}}\psi)_k}(\ovl{T^{(n-2)}}\phi)_k
    \end{split}
\end{equation*}
so that
$\ovl{T^{(n-2)}}(\phi\tiund{\psi}\phi)=\ovl{T^{(n-2)}}\phi\tiund{\ovl{T^{(n-2)}}\psi}
\ovl{T^{(n-2)}}\phi$. Since $\psi$ is discrete, by induction hypothesis (iii),
so is $T^{(n-2)}\psi_k$ for each $k$, and therefore also
$\ovl{T^{(n-2)}}\psi\in\dbb_2$ is discrete. From step $n=2$ of the induction,
$\tiu$ preserves fibre products over discrete objects. Together with the
analogous property of $\ovl{T^{(n-2)}}$ explained above, this implies
\begin{equation*}
    \begin{split}
       & T^{(n-1)}(\phi\tiund{\psi}\phi)=\tiu\ovl{T^{(n-2)}}(\phi\tiund{\psi}\phi)=
       \tiu(\ovl{T^{(n-2)}}\phi\tiund{\ovl{T^{(n-2)}}\psi}\ovl{T^{(n-2)}}\phi)= \\
        &=\tiu\ovl{T^{(n-2)}}\phi\tiund{\tiu\ovl{T^{(n-2)}}\psi}\tiu\ovl{T^{(n-2)}}\phi)=
        T^{(n-1)}\phi\tiund{T^{(n-1)}\psi}T^{(n-1)}\phi.
    \end{split}
\end{equation*}
Thus $T^{(n-1)}$ preserves fibre products over discrete objects, as required.
It is clear that, when $\psi$ is discrete, so
$T^{(n-1)}\psi=\tiu\ovl{T^{(n-2)}}\psi$.

It remains to show that $T^{(n-1)}$ preserves weak equivalences. Since (i) and
(ii) have been proved at step $n$, by Lemma \ref{semist.lem01} we have a
functor $V_n:\dbb_n\rw\hcl_\np$ sending weak equivalences to
$(n+1)$-equivalences. Let $f:\phi\rw\psi$ be a weak equivalence in $\dbb_n$;
then $V_n f$ is an $(n+1)$-equivalence in $\hcl_\np$. In particular (see Remark
\ref{tam3.rem1} b)), $(V_n f)_1=U_n f:U_n\phi\rw U_n\psi$ is an $n$-equivalence
in $\tcl_n$ and therefore $\ta{n}{1}U_n f:\ta{n}{1}U_n \phi\rw\ta{n}{1}U_n
\psi$ is an equivalence of groupoids. But from having proved (ii) at step $n$,
we have $\ta{n}{1}U_n f=U_1 T^{(\nm)}f$ so $U_1 T^{(\nm)}f$ is a weak
equivalence. Hence, arguing as in the proof of Corollary \ref{tam3.cor1}, we
conclude that the map $T^{(\nm)}f:T^{(\nm)}\phi\rw T^{(\nm)}\psi$ of $\cgp$ is
a weak equivalence. This completes the proof of (iii) at step $n$.\enpr
\begin{theorem}\label{semist.the01}
    Let $V_n$ be as in (\ref{deloop.eq3bis}). For each $n\geq 2$, $V_n$ restricts to
    a functor $V_n:\dbb_n\rw\hcl_\np$
    which preserves the homotopy type and sends weak equivalences to
    $(n+1)$-equivalences.
\end{theorem}
\prf From Proposition \ref{semist.pro01} for each $n\geq 2$ conditions (i) and
(ii) hold. Hence by Lemma \ref{semist.lem01} the theorem follows.\enpr\bigskip

Using the main theorems of Sections \ref{specia} and \ref{weak} we immediately
deduce the existence of a comparison functor from $\ctn$s to the Tamsamani
model.
\begin{corollary}\label{semist.cor01}
    There is a functor $F:\cagp\rw\hcl_\np$ such that
    $B F \gcl\cong B\gcl$ in $\hcl o(\tp^{(\np)}_{*})$ for every $\gcl\in\cagp$.
\end{corollary}
\prf Let $F=V_n\dcl_n\spi$. By Theorems \ref{specgen.the1}, \ref{globgen.the1}
and \ref{semist.the01}, for each $\gcl\in\cagp$ we have $B F \gcl=B
V_n\dcl_n\spi\gcl=B\spi\gcl\simeq B\gcl$. Hence there is an isomorphism
$BF\gcl\cong B\gcl$ in $\hcl o(\tp^{(\np)}_{*})$.\enpr

\subsection{Semistrictification result}\label{result}
In this section we prove our main semistrictification result. Recall from
Section \ref{catn} that there is a functor $\pcl_n:\tp\rw\cagp$ which gives an
equivalence of categories (Theorem \ref{catn.the1}) $\cagp/\!\!\sim\,\simeq\hcl
o(\tp_*^{(\np)})$ where $\tp_*^{(\np)}$ is the category of connected
$(n+1)$-types. Consider now the composite functor
\begin{equation*}
    S:\tcl_\np\supar{B}\tp\supar{\pcl_n}\cagp\supar{Sp}\cagp_S\supar{\dcl_n}\dbb_n
    \supar{V_n}\hcl_\np.
\end{equation*}
Since the image of this functor lies in the category $\hcl_\np$ of semistrict
Tamsamani $(n+1)$-groupoids, we call $S$ the \emph{semistrictification
functor}.

We saw in Corollary \ref{semist.cor01} that the composite $F=V_n\dcl_nSp$
preserves the homotopy type. The functor $\pcl_n:\tp\rw\cagp$ does not preserve
the homotopy type but, from Theorem \ref{catn.the1}, its restriction to the
subcategory $\tp_*^{(\np)}$ of connected $(n+1)$-types does preserve the
homotopy type. The image of the functor $B:\tcl_\np\rw\tp$ is a $(n+1)$-type
but not, in general, a connected $(n+1)$-type. Hence the functor $S$ does not
in general preserve the homotopy type. If, however, $\gcl\in\tcl_\np$ is such
that $B\gcl$ is a connected $(n+1)$-type, then from the above discussion
$B\gcl$ and $BS\gcl$ do have the same homotopy type; that is, $B\gcl\cong
BS\gcl$ in $\hcl o (\tp_*^{(\np)})$. Using Proposition \ref{group.pro1} this
gives immediately the following theorem, which is our semistrictification
result.
\begin{theorem}\label{result.the1}
    Let $\gcl\in\tcl_\np$ be such that $B\gcl$ is a connected $(n+1)$-type.
    There exists a zig-zag of $(n+1)$-equivalences in $\tcl_\np$ connecting
    $\gcl$ and $S(\gcl)\in\hcl_\np$.
\end{theorem}
\prf Since $B\gcl$ is a connected $(n+1)$-type, from Theorem \ref{catn.the1},
$B\gcl\cong B\pcl_n B\gcl$ in $\hcl o (\tp_*^{(\np)})$. From Corollary
\ref{semist.cor01}, $B\pcl_n B\gcl\cong B F\pcl_n B\gcl=BS\gcl$ in $\hcl o
(\tp_*^{(\np)})$. Therefore $B\gcl\cong BS\gcl$ in $\hcl o (\tp_*^{(\np)})$.
Hence there is a zig-zag of weak equivalences in $\tp_*^{(\np)}$ connecting
$B\gcl$ and $BS\gcl$. From Proposition \ref{group.pro1} it follows that there
is a zig-zag of $(n+1)$-equivalences in $\tcl_\np$ connecting $\gcl$ and
$S(\gcl)$.\enpr\bigskip

The next corollary illustrates the relevance of our semistrictification result
to the modelling of connected $(\np)$-types. A further discussion about this
can be found in Section \ref{furth} c).

Let $\,\tcl_\np^*\,$ denote\, the \,full\, subcategory\, of $\tcl_\np$ of\,
those\, Tamsamani weak $(\np)$-groupoids whose classifying space is path
connected. Clearly $\hcl_\np$ is a full subcategory of $\tcl_\np^*$. It is
immediate that the equivalence of categories  $\tcl_\np/\!\!\sim^\np\simeq\hcl
o(\tp^{(\np)})$ of \cite{tam} restricts to $\tcl^*_\np/\!\!\sim^\np\simeq\hcl
o(\tp_*^{(\np)})$. In the next corollary we will see that Theorem
\ref{result.the1} implies that $\hcl o(\tp_*^{(\np)})$ is equivalent to a
category which is smaller than $\tcl^*_\np/\!\!\sim^\np$. More precisely, let
$\hcl o_{\hcl_\np}(\tcl_\np)$ be the full subcategory of
$\tcl^*_\np/\!\!\sim^\np$ whose objects are in $\hcl_\np$. Then we have
\begin{corollary}\label{result.cor1}
    There is an equivalence of categories
    \begin{equation}\label{result.eq1}
        \hcl o_{\hcl_\np}(\tcl^*_\np)\simeq\hcl o(\tp_*^{(\np)})
    \end{equation}
\end{corollary}
\prf Let $Q_n:\tp_*^{(\np)}\rw\hcl_\np$ be the composite
$Q_n=V_n\dcl_n\spi\pcl_n$ so that $Q_n B=S$. Since
$\pcl_n:\tp_*^{(\np)}\rw\cagp$ as well as $\spi$, $\dcl_n$, $V_n$ preserve the
homotopy type and $V_n$ sends weak equivalences to $(\np)$-equivalences, $Q_n$
induces a functor $Q_n:\hcl o (\tp_*^{(\np)})\rw\hcl o_{\hcl_\np}(\tcl^*_\np)$.
The classifying space also induces a functor $B:\hcl
o_{\hcl_\np}(\tcl^*_\np)\rw\hcl o(\tp^{(\np)}_{*})$. We claim that $(B,Q_n)$
gives the equivalence of categories (\ref{result.eq1}). In fact, given
$\gcl\in\hcl o_{\hcl_\np}(\tcl^*_\np)$, $B\gcl$ is a connected $(\np)$-type,
hence by Theorem \ref{result.the1} there is a zig-zag of $(\np)$-equivalences
in $\tcl_\np$ connecting $\gcl$ and $S\,\gcl=Q_n B\gcl$, say
\begin{equation*}
    \gcl\rTo X_1\lTo\cdots\rTo X_n\lTo Q_n B\gcl.
\end{equation*}
Then all $BX_i$ are connected $(\np)$-types, so this is a zig-zag in
$\tcl_\np^*$. Hence $\gcl$ and $Q_n B\gcl$ are isomorphic in
$\tcl^*_\np/\!\!\sim^\np$. Since both $\gcl$ and $Q_n B\gcl$ are $\hcl_\np$, by
definition of $\hcl o_{\hcl_\np}(\tcl_\np)$, it follows that $\gcl\cong Q_n
B\gcl$ in $\hcl o_{\hcl_\np}(\tcl_\np)$. Let $X\in \tp_*^{(\np)}$, then by
Theorem \ref{catn.the1} and Corollary \ref{semist.cor01} there are isomorphisms
 $B Q_n X=B F \pcl_n X\cong B \pcl_n X\cong X$ in $\hcl o(\tp_*^{(\np)})$. This
concludes the proof that $(B,Q_n)$ is an equivalence of categories. \enpr

\section{Themes for further investigations}\label{furth}
We would like to conclude this paper with a list of open questions and
conjectures that arise from our results. Although extremely interesting, the
following questions are substantial projects and we will tackle them in
subsequent papers.\medskip

\textbf{a)} In \cite{paol} we considered two distinct semistrictifications of
Tamsamani weak 3-groupoids for the path-connected case. One is the same
category $\hcl_3$ considered in this paper. The other is the subcategory
$\kcl_3$ of $\tcl_3$ consisting of those objects
$\phi\in\tcl_3\subset[\dop,\tcl_2]$ such that $\phi_0=\delta^{(2)}\{\cdot\}$
and for each $[k]\in\dop$ $\phi_k$ is a strict 2-groupoid. Objects of $\kcl_3$
are not strict 3-groupoids because the Segal maps $\phi_k\rw\fistr{k}$ are just
2-equivalences. The categories $\hcl_3$ and $\kcl_3$ both have objects which
are semistrict, but they do not coincide because if $\phi\in\hcl_3$ each
$\phi_k$ is in general a weak, not a strict, 3-groupoid while the  Segal maps
$\phi_k\rw\fistr{k}$ are isomorphisms. We constructed in \cite{paol} a functor
$S':\tcl_3\rw\kcl_3$ and showed that, given $\phi\in\tcl_3$ such that $B\phi$
is path-connected there is a zig-zag of 3-equivalences in $\tcl_3$ connecting
$\phi$ and $S'\phi$. We also associated to each object of $\hcl_3$ and of
$\kcl_3$ a Gray groupoid with one object having the same homotopy type.

Gray groupoids are well known to model 3-types \cite{jt}, \cite{le}. More
recently, a different type of semistrict structure has been proposed for
modelling homotopy types. Carlos Simpson conjectured \cite{simp} that every
homotopy $n$-type arises as the geometric realization of a strict $n$-groupoid
equipped with a suitable notion of weak identity arrow. An elegant
formalization of this conjecture was given by Joachim Kock \cite{kock} in terms
of so called ``fair $n$-groupoids". Thus the \emph{Kock-Simpson conjecture} for
homotopy types states that fair $n$-groupoids model $n$-types. At present a
proof of this conjecture exists for 1-connected 3-types \cite{joko}.

We conjecture that the two semistrictifications of Tamsamani weak 3-group-oids
in terms of $\kcl_3$ and $\hcl_3$ correspond respectively to Gray groupoids and
fair 3-groupoids for the path connected case. Intuitively it is clear that
$\kcl_3$ rather than $\hcl_3$ corresponds to Gray groupoids because in a Gray
groupoid the $\hom$ objects between two 0-cells form a strict 2-groupoid. More
generally, for any $n\geq 3$, we conjecture that $\hcl_n$ corresponds to fair
$n$-groupoids (with one object).

We think that a suitable approach to this question could be to establish a
direct comparison between special $\ctmen$s and fair $n$-groupoids. Namely, we
conjecture that there is a functor relating directly these two categories and
preserving the homotopy type which exhibits the strongly contractible faces of
a special $\ctn$ as ``spaces of weak units". Since by Corollary \ref{spec.cor1}
special $\ctmen$s model connected $n$-types, this approach would provide a
potential tool for proving the Kock-Simpson conjecture for connected
$n$-types.\medskip

\textbf{b)} Our semistrictification result establishes a zig-zag of
$(\np)$-equivalences between $\gcl\in\tcl_\np$ and $S\gcl\in\hcl_\np$ when
$B\gcl$ is path-connected. One could ask if this zig-zag could be made into a
single map using some higher order information. We know for instance that
biequivalences of bicategories have pseudo-inverses so one could speculate that
in a theory of weak $n$-categories $n$-equivalences should also have
pseudo-inverses when regarding $n$-categories not just as a category but as
some suitable $(\np)$-category.

However, the results of \cite{lp} illustrate that in the Tamsamani model, at
least for $n=2$, this principle is not satisfied precisely in this form:
instead one needs to modify the morphisms of $\wcl_2$ (viewed as a 2-category)
to include pseudo-natural transformations between functors from $\dop$ to Cat,
in order for 2-equivalences to have pseudo-inverses. The modified 2-category
$\wcl_2^{ps}$ is then proved in \cite[Theorem 7.3]{lp} to be biequivalent to
the 2-category of bicategories, normal homomorphisms and oplax natural
transformations with identity components. These considerations lead us to
speculate that in order to invert $n$-equivalences in $\tcl_n$ we may need to
enlarge the class of morphisms in the Tamsamani model to include some kind of
pseudo-natural maps.\medskip

\textbf{c)} In Corollary \ref{result.cor1} we showed that ${\hcl
o(\tp^{(\np)}_*)}$ and $\hcl o_{\hcl_\np}(\tcl^*_\np)$  are equivalent. This is
homotopically significant because the category $\hcl o_{\hcl_\np}(\tcl^*_\np)$
is smaller that $\tcl^*_\np/\!\!\sim^\np$, being the full subcategory whose
objects are in $\hcl_\np$. On the other hand, we can also consider the
localization $\hcl_\np/\!\!\sim^\np$. Clearly $\hcl_\np/\!\!\sim^\np$ is a
subcategory of $\hcl o_{\hcl_\np}(\tcl^*_\np)$ with the same set of objects,
but a priori it is not necessarily a full subcategory. An interesting question
we point out is whether $\hcl_\np/\!\!\sim^\np$ is a full subcategory of $\hcl
o_{\hcl_\np}(\tcl^*_\np)$. If this was the case, then clearly
$\hcl_\np/\!\!\sim^\np$ and $\hcl o_{\hcl_\np}(\tcl^*_\np)$ would coincide
hence, by Corollary \ref{result.cor1}, $\hcl_\np/\!\!\sim^\np$ would be
equivalent to $\hcl o(\tp^{(\np)}_*)$.\medskip

\textbf{d)} One could ask if the ideas of this paper could be extended to the
non-path connected case as well as to the non-groupoidal case, to establish
corresponding semistrictification results for a general Tamsamani weak
$n$-groupoid and for Tamsamani weak $n$-categories. We believe that in the
first case the ideas of this paper could be implemented;  the latter case
presents more challenging issues. In fact, there is in the literature
\cite{cere} a model of 3-types called in \cite{cere} 2-cat groupoids; this
reduces to cat$^2$-groups in the path connected case. Using this model it is
possible to adapt the ideas of this paper for $n=2$ (as well as the other
results in \cite{paol}) to obtain corresponding semistrictification results for
general Tamsamani weak 3-groupoids. A proof of this was presented by the author
in \cite{pata}. When $n>2$, a generalization of the model in \cite{cere} does
not seem to have yet appeared in the literature, but it is very reasonable to
conjecture that it exists. We think it would be possible to use it to adapt the
ideas and techniques of this paper to the general non path-connected case.

The non-groupoidal case is more challenging. Here is, for instance, one of the
issues that one encounters when trying to extend the structures studied in this
paper to a non-groupoidal setting. One could look for a suitable notion of
``special $n$-fold category" to be compared to Tamsamani weak $n$-categories,
in the same way as special $\ctn$s relate to the Tamsamani model. Such a notion
would require a notion of strong contractibility and then the requirement that
the faces, which in a globular object are discrete, should be strongly
contractible. This could be defined formally as in the case of special $\ctn$s.

We observe, however, that in order to extend the functor $\dcl_n$ to this more
general setting, we would need some extra conditions, which are automatically
satisfied in the case of special $\ctn$s. Namely, remember that an important
point that made possible the passage from $\cagp_S$ to $\dbb_n$ via the functor
$\dcl_n$ is the fact that all pullbacks occurring in the nerve of a special
$\ctn$ are weakly equivalent to homotopy pullbacks. More precisely, this was
the reason why the Segal maps in $\dcl_2\gcl\in[\dop,\cgp]$ are weak
equivalences (see proof of Theorem \ref{glob2.the1}, diagram
(\ref{glob2.eq1})). Similarly when $n>2$, after passing diagonals (see proof of
Theorem \ref{globgen.the1}). We would need to require a similar condition to
hold in a ``special $n$-fold category". This could be achieved, for instance,
by requiring the source and target maps to be fibrations in an appropriate
model structure on $(\nm)$-fold categories. Model structures on $n$-fold
categories are only now beginning to be investigated. Even for $n=2$,
interesting model structures on double categories are non trivial to establish;
progress in this direction has been made in joint work of the author
\cite{fpp}.

\section{Appendix}\label{app}
\emph{\textbf{Proof of Lemma \ref{oper.lem3}}}\medskip

It is known from general theory (see for instance \cite{orz}) that $\dcl$ is
cocomplete, so the coproduct $A\amalg B$ exists. Let $V$ be the subobject of
$A\amalg B$ consisting of elements of the form
\begin{equation}\label{app.eq1}
    u_1(a_1)u_2(b_1)u_1(a_2)u_2(b_2)...\;\text{ or }\;
    u_2(b_1)u_1(a_1)u_2(b_2)u_1(a_2)...
\end{equation}
where $a_i\in A,\;\; b_i\in B$. Since $\omega(a,b)=\omega a\; \omega b$ for
each $\omega\in\Omega'_1$, $V$ is closed under the operations in $\dcl$, so it
is a subobject of $A\amalg B$. By the universal property of coproducts, we have
the following commutative diagram
\begin{equation*}
     \xy
    0;/r.25pc/:
    (0,0)*{A\amalg B}="1";
    (0,15)*{V}="2";
    (-20,-15)*{\overset{}{A}}="3";
    (20,-15)*{\overset{}{B}}="4";
    (17,-14)*{}="5";
    (-17,-14)*{}="6";
    {\ar@/^1.2pc/^{u_1}"3";"2"};
    {\ar@/_1.2pc/_{u_2}"4";"2"};
    {\ar@{^(->}_{u_1}"6";"1"};
    {\ar@{_(->}^{u_2}"5";"1"};
    {\ar_\id"1";"2"};
    \endxy
\end{equation*}
Therefore $A\amalg B=V$, so every element of $A\amalg B$ has the form
(\ref{app.eq1}), as required.\enpr\bigskip

\emph{\textbf{Proof of Theorem \ref{oper.the1}}}\medskip

We need to show that there is an isomorphism
\begin{equation}\label{app.eq2}
    \hund{\cat\dcl}(\lcl H,\abb)\cong\hund{\dcl}(H,\ker\pt_0\times A_0)
\end{equation}
where $\lcl H=(H\amalg H)_p$ as in the statement of the theorem. Let
$\pt_0,\pt_1,\sigma_0$ be the source, target and identity maps of the internal
category $\abb$ and denote by $i$ the inclusion of $\ker \pt_0$ in $A_1$.
Suppose we are given a morphism $(h,g):H\rw\ker\pt_0\times A_0$ in $\dcl$. We
are going to build an internal functor $\alpha(h,g):\lcl H\rw\abb$. Let
$\gamma: H\amalg H\rw A_1$ and $f: H\amalg H\rw A_0$ be defined by
\begin{equation}\label{app.eq3}
    \gamma u_1=ih,\quad\gamma u_2=\sigma_0 g,\quad f u_1=\pt_1 i h,\quad f
    u_2=g.
\end{equation}
Notice that we have
\begin{equation}\label{app.eq4}
    \pt_1\gamma=f,\quad gp=\pt_0 \gamma
\end{equation}
where $p:H\amalg H\rw H$ is determined by $pu_1=1$, $pu_2=\id$. In fact,
composing with the coproduct injections and using (\ref{app.eq3}) we find
\begin{equation*}
    \begin{split}
        & \pt_1\gamma u_1=\pt_1 i h=fu_1,\quad \pt_1\gamma u_2=\pt_1\sigma_0 g=g=fu_2 \\
         & \pt_0\gamma u_1=\pt_0 i h=1=gpu_1,\quad \pt_0\gamma u_2=\pt_0\sigma_0 g=g=gpu_2
     \end{split}
\end{equation*}
from which (\ref{app.eq4}) follows.

Since $gp=\pt_0\gamma$ it follows from Lemma \ref{oper.lem2} a) that there is a
morphism $(H\amalg H)_p\rw(A_1)_{\pt_0}$ in $\cat\dcl$; that is,
\begin{diagram}
    \cdots &\rTo & \ker p\ttil(H\amalg H) & \pile{\rTo\\ \rTo\\ \lTo} & H\amalg
    H\\
    && \dTo_{(\gamma_|,\gamma)}&& \dTo_{\gamma}\\
    \cdots &\rTo & \ker \pt_0\ttil A_1 & \pile{\rTo\\ \rTo\\ \lTo} & A_1,
\end{diagram}
where the structural maps are as in Lemma \ref{oper.lem2} a) and where we have
denoted by $\gamma_|$ the restriction of $\gamma$ to $\ker p$. We are now going
to build an internal functor $(r,\pt_1):(A_1)_{\pt_0}\rw\abb$ as follows:
\begin{diagram}[h=2.8em]
    \cdots &\rTo & \ker \pt_0\ttil A_1 & \pile{\rTo^{d_0}\\ \rTo^{d_1}\\ \lTo_{s_0}} & A_1\\
    && \dTo_{r}&& \dTo_{\pt_1}\\
    \cdots &\rTo &  A_1 & \pile{\rTo^{\pt_0}\\ \rTo^{\pt_1}\\ \lTo_{\sigma_0}} & A_0.
\end{diagram}
We define $r$ by $r(x,y)=x(\sigma_0\pt_1 y)$. To show that $(r,\pt_1)$ is an
internal functor we need to check that $r$ is a morphism in $\dcl$ and that
$(r,\pt_1)$ is compatible with the structural maps. For this recall, from
Section \ref{oper} that, since $\dcl$ is a category of groups with operations,
$[\ker\pt_0,\ker\pt_1]=1$. Hence, for each $(x,y),(x',y')\in\ker\pt_0\ttil A_1$
and $\omega\in\Omega'_1$, we compute:
\begin{equation*}
    \begin{split}
       & r((x,y)(x',y'))=r(xyx'y^{-1},yy')=xyx'y^{-1}(\sigma_0\pt_1 yy')=
       \\
       & =xyx'y^{-1}(\sigma_0\pt_1 y)(\sigma_0\pt_1y')= xyy^{-1}(\sigma_0\pt_1 y)x'(\sigma_0\pt_1y')=\\
       & =r(x,y)r(x',y'); \\
    \end{split}
\end{equation*}
\begin{equation*}
    r(\omega(x,y))=r(\omega x,\omega y)=\omega x(\sigma_0\pt_1\omega y)=\omega
    x\;\;\omega\sigma_0\pt_1 y=\omega r(x,y).
\end{equation*}
This shows that $r$  is a morphism in $\dcl$. As for the compatibility with the
structural maps, recalling from Lemma \ref{oper.lem2} that $d_0(x,y)=y$,
$d_1(x,y)=xy$ and $s_0(z)=(1,z)$, we have
\begin{equation*}
    \begin{split}
       & \pt_0r(x,y)=\pt_0 x(\pt_0\sigma_0\pt_1y)=\pt_1y=\pt_1d_0(x,y), \\
       & \pt_1r(x,y)=\pt_1 x(\pt_1\sigma_0\pt_1y)=\pt_1x\pt_1y=\pt_1(xy)=\pt_1d_1(x,y),  \\
        & rs_0(z)=r(1,z)=\sigma_0\pt_1z.
    \end{split}
\end{equation*}
Compatibility with the composition maps holds automatically. Hence $(r,\pt_1)$
is a morphism in $\cat\dcl$. Let $\alpha(h,g):\lcl H\rw \abb$ be the composite
\begin{equation*}
    \lcl H=(H\amalg H)_p\supar{} (A_1)_{\pt_0}\supar{(r,\pt_1)}\abb.
\end{equation*}
Since, by (\ref{app.eq4}), $\pt_1\gamma=f$, we have
\begin{equation}\label{app.eq5}
    \alpha(h,g)=(r(\gamma_|,\gamma),f)
\end{equation}
where $\gamma$ and $f$ are as in (\ref{app.eq3}), and $r(x,y)=x(\sigma_0\pt_1
y)$ for $(x,y)\in\ker\pt_0\ttil A_1$.

Conversely, suppose we are given a morphism $\;(v,k):\lcl H\rw \abb$ in
$\cat\dcl$,
\begin{diagram}[h=4em]
    \cdots &\rTo & \ker p\ttil(H\amalg H) & \pile{\rTo^{d_0}\\ \rTo^{d_1}\\ \lTo_{s_0}} & H\amalg
    H\\
    && \dTo_{v}&& \dTo_{k}\\
    \cdots &\rTo &  A_1 & \pile{\rTo^{\pt_0}\\ \rTo^{\pt_1}\\ \lTo_{\sigma_0}} & A_0.
\end{diagram}
We are going to build a morphism $\beta(v,k):H\rw\ker\pt_0\times A_0$. If $z\in
u_1(H)$ then, by definition of $p$, $z\in\ker p$ and also $\pt_0 v(z,1)=k
d_0(z,1)=1$, so that $v(z,1)\in\ker \pt_0$. Thus we define $h:H\rw\ker \pt_0$
by $h(x)=v(u_1(x),1)$. We also take $g:H\rw A_0$ to be $g=k u_2$ and finally
let $\beta(v,k)=(h,g)$. That is, for each $x\in H$,
\begin{equation}\label{app.eq6}
    \beta(v,k)(x)=(v(u_1(x),1),k u_2(x)).
\end{equation}
To prove the proposition it remains to show that $\beta\alpha=\id$ and
$\alpha\beta=\id$. From (\ref{app.eq3}), (\ref{app.eq5}) and (\ref{app.eq6}) we
have, for each $x\in H$,
\begin{equation*}
\begin{split}
   & \beta\alpha(h,g)(x)=\beta(r(\gamma_|,\gamma),f)(x)=(r(\gamma u_1(x),1),f
    u_2(x))= \\
    & =(\gamma u_1(x),f u_2(x))=(ih(x),g(x)).
\end{split}
\end{equation*}
Hence $\beta\alpha=\id$. From (\ref{app.eq5}) and (\ref{app.eq6}) we calculate
\begin{equation}\label{app.eq7}
     \alpha\beta(v,k)=\alpha(v(u_1(\mi),1),ku_2)=(r(\gamma_|,\gamma),f)
\end{equation}
where $\gamma$ and $f$ are determined by
\begin{equation}\label{app.eq8}
\begin{split}
   & \gamma u_1(x)=v(u_1(x),1)\qquad \gamma u_2=\sigma_0 k u_2=vs_0u_2\quad \\
    & fu_1(x)=\pt_1 i v(u_1(x),1)\qquad f u_2=k u_2.
\end{split}
\end{equation}
and $r$ is given by $r(x,y)=x(\sigma_0\pt_1 y)$.

Since $(v,k)$ is an internal functor, $\pt_1 v=kd_1$; therefore, by
(\ref{app.eq8}), $fu_1(x)=\pt_1 iv(u_1(x),1)=kd_1(u_1(x),1)=ku_1(x)$. In
conclusion, $fu_1=ku_1$; also, by (\ref{app.eq8}), $fu_2=ku_2$, so that, in
conclusion, $f=k$.

From (\ref{app.eq7}) in order to show that $\alpha\beta=\id$ it remains to
prove that
\begin{equation}\label{app.eq9}
    r(\gamma_|,\gamma)=v.
\end{equation}
For each $(x,y)\in\ker p\ttil(H\amalg H)$, using the fact proved in
(\ref{app.eq4}) that $\pt_1\gamma=f$, the definition of $r$ and the fact proved
above that $f=k$, we compute
\begin{equation}\label{app.eq10}
    r(\gamma_|,\gamma)(x,y)=\gamma(x)\sigma_0\pt_1\gamma(y)=\gamma(x)\sigma_0
    f(y)=\gamma(x)\sigma_0 k(y).
\end{equation}
On the other hand, since $v$ is a group homomorphism
\begin{equation}\label{app.eq11}
    v(x,y)=v(x,1)v(1,y),
\end{equation}
we are going to show that
\begin{equation}\label{app.eq12}
    v(x,1)=\gamma(x)\quad \text{and}\quad \sigma_0 k(y)=v(1,y)
\end{equation}
for each $x\in\ker p$, $y\in H\amalg H$. Then (\ref{app.eq10}),
(\ref{app.eq11}), (\ref{app.eq12}) immediately imply (\ref{app.eq9}). The
second identity in (\ref{app.eq12}) is trivial since $v(1,y)=vs_0(y)=\sigma_0 k
y$. The proof that $v(x,1)=\gamma(x)$ for each $x\in\ker p$ is more elaborate
and needs some preliminary observations.

The first observation we need is that every element of the underlying group of
$\ker p$ is a product of elements of the form $wu_1(z)w^{-1}$, where $z\in H$
and $w\in H\amalg H$. In fact, let $v\in\ker p$. By Lemma \ref{oper.lem3}, $v$
has the form (\ref{app.eq1}); say for instance $v=u_1(h_1)u_2(h'_1)\ldots
u_1(h_n)u_2(h'_n)$. Then, since $pu_1=1$ and $pu_2=\id$, we have $1=p(v)=h'_1
h'_2\ldots h'_n$ so that $u_2(h'_1)=u_2((h'_2\ldots h'_n)^{-1})$. Hence if we
put, for each $i=2,\ldots,n$ $w_i=u_2(h'_i\ldots
h'_n)^{-1}u_1(h_i)u_2(h'_i\ldots h'_n)$, we have $v=u_1(h_1)w_2w_3\ldots w_n$.
The case where $v$ is of the other form in (\ref{app.eq1}) is similar.

The second observation we need is the trivial fact that, since $v$ and $\gamma$
are group homomorphisms, if the identity $v(x,1)=\gamma(x)$ is satisfied for
$x=x_1\in\ker p$ and $x=x_2\in\ker p$, then it is also satisfied for
$x=x_1x_2$. It then follows from above that in order to show that
$v(x,1)=\gamma(x)$ for all $x\in\ker p$, it is sufficient to show that
\begin{equation}\label{app.eq13}
    v(wu_1(z)w^{-1},1)=\gamma(wu_1(z)w^{-1})
\end{equation}
for each $z\in H$, $w\in H\amalg H$. To show (\ref{app.eq13}) we first observe
that, if (\ref{app.eq13}) holds for $w=w_1$ it also holds for $w=u_1(x)w_1$ and
for $w=u_2(x)w_1$. To prove this, recall from (\ref{app.eq8}) that $\gamma
u_1(x)=v(u_1(x),1)$ and $\gamma u_2(x)=vs_0u_2(x)=v(1,u_2(x))$; thus, using the
expression of the group multiplication in the semidirect product $(\ker
p)\ttil(H\amalg H)$ as given in Lemma \ref{oper.lem2}, we calculate:
\begin{equation*}
\begin{split}
   & v(u_1(x)w_1u_1(z)w_1^{-1}u_1(x_1^{-1}),1)=v((u_1(x),1)(w_1u_1(z)w_1^{-1},1)(u_1(x^{-1}),1))= \\
    & =\gamma u_1(x)\gamma(w_1u_1(z)w_1^{-1})\gamma
    u_1(x^{-1})=\gamma(u_1(x)w_1u_1(z)w_1^{-1}u_1(x^{-1})), \text{ and}
\end{split}
\end{equation*}
\begin{equation*}
    \begin{split}
       & v(u_2(x)w_1u_1(z)w_1^{-1}u_2(x^{-1}),1)=v((1,u_2(x))(w_1u_1(z)w_1^{-1},1)(1,u_2(x^{-1}))= \\
       & =\gamma u_2(x)\gamma(w_1 u_1(z)w_1^{-1})\gamma
       u_2(x^{-1})=\gamma(u_2(x)w_1u_1(z)w_1^{-1}u_2(x^{-1})).
    \end{split}
\end{equation*}
Hence (\ref{app.eq13}) holds for $w=u_1(x)w_1$ and for $w=u_2(x)w_1$ as
claimed. Since, by Lemma \ref{oper.lem3} $w$ has the form (\ref{app.eq1}) and,
by (\ref{app.eq8}), (\ref{app.eq13}) holds for $w=1$, it follows from above
that it holds for each $w\in H\amalg H$ and $z\in H$, as required.\enpr\bigskip

\bigskip

\emph{\textbf{Proof of Lemma \ref{adj.lem1}}}\medskip

By induction on $n$. For $n=1$, recall from the proof of Theorem \ref{adj.the0}
that $\ucl_1:\ccl^1\gcl\rw\set$ is the composite
$\ccl^1\gcl\supar{\beta_1}\cat(\gp)\supar{\nu_1}\gp\supar{U}\set$ and that,
from Lemma \ref{oper.lem1}, given $\gcl=(G,d,t)\in\ccl^1\gcl$, $\beta_1\gcl$ is
\begin{equation*}
    \gcl\tiund{\im d}\;\gcl\rTo \gcl\pile{\rTo^{d}\\ \rTo{t}\\ \lTo_{i}}\im d.
\end{equation*}
Hence $\ucl_1\gcl=U(\ker d\times \im d)$. On the other hand, since $d$ has a
section $i$, it is $G\cong\ker d\rtimes \im d$. It follows that $\ucl_1
\gcl=U(\ker d \times \im d)=U(\ker d \rtimes \im d)\cong UG= R_1\gcl$.
Naturality is immediate.

Inductively, suppose that $\ucl_\nm\gcl\cong R_\nm\gcl$ for each
$\gcl\in\cngm$. Recall from the proof of Theorem \ref{adj.the0} that
$\ucl_n\!=\!\ucl_\nm\vcl_n\beta^{(k)}_n$ and let
$\gcl=(G,d_1,...,d_n,t_1,...,t_n)$ $\in\cng$. Put
$\gcl'=(G,d_1,...,d_{k-1}\;d_{k+1}\ldots,d_n,t_1,...,t_{k-1}\;t_{k+1}...,t_n)\in\cngm$.
Then $\beta^{(k)}_n\gcl$ is
\begin{equation*}
    \gcl'\tiund{\im d_k}\;\gcl'\rTo \gcl'\pile{\rTo^{d_k}\\ \rTo{t_k}\\ \lTo_{i}}\im d_k.
\end{equation*}
Therefore by the inductive hypothesis
\begin{equation}\label{app.eq14}
    \ucl_n\gcl=\ucl_\nm\ker d_k\times\ucl_\nm\im d_k\cong R_\nm\ker
    d_k\times R_\nm\im d_k.
\end{equation}
On the other hand, since $d_k$ has a section $i$, $\gcl'$ is isomorphic to the
semidirect product $\gcl'\cong\ker d_k\ttil\im d_k$. From the description of
$\cng$ as a category of groups with operations, the functor $R_n$ is precisely
the underlying set functor; recall from the general discussion about semidirect
products in a category of groups with operations in Section \ref{oper} that the
underlying set of the semidirect product object is the product of the
underlying sets. Hence
\begin{equation}\label{app.eq15}
    R_\nm\gcl'\cong R_\nm(\ker d_k\ttil\im d_k)\cong R_\nm\ker
    d_k\times R_\nm\im d_k.
\end{equation}
It follows from (\ref{app.eq14}) and (\ref{app.eq15}) that $\ucl_n\gcl\cong
R_\nm\gcl' = UG= R_n\gcl$, as required. Naturality is immediate. \enpr\bigskip

\emph{\textbf{Proof of Corollary \ref{adj.cor1}}}\medskip

From Lemma \ref{adj.lem1} it is enough to show that $R_n:\cng\rw\set$ reflects
regular epis.

Let $\gcl=(G,d_1,\ldots,d_n,$ $t_1,\ldots,t_n)$ and
$\gcl'=(G',d'_1,\ldots,d'_n,t'_1,\ldots,t'_n)$ be two objects in $\cng$ and
$f:\gcl\rw\gcl'$ a morphism in $\cng$ such that $R_n f$ is a regular epi. Since
$\set$ is a regular category \cite{borc}, $R_n f$ is the coequalizer of its
kernel pair, hence the following is a coequalizer
\begin{equation*}
    U(\tens{G}{G'})=\tens{UG}{UG'}\pile{\rTo^{p_1}\\
    \rTo_{p_2}}UG\rTo^{Uf} UG'.
\end{equation*}
Since $U:\gp\rw\set$ is monadic, it reflects regular epis, so the following is
a coequalizer in $\gp$:
\begin{equation}\label{app.eq16}
    \tens{G}{G'}\pile{\rTo^{p_1}\\ \rTo_{p_2}}G\rTo^{f} G'.
\end{equation}
We claim that the following is a coequalizer in $\cng$:
\begin{equation*}
    \tens{\gcl}{\gcl'}\pile{\rTo^{p_1}\\ \rTo_{p_2}}\gcl\rTo^{f} \gcl'.
\end{equation*}
To see this, let $\gcl''=(G'',d''_1,\ldots,d''_n,t''_1,\ldots,t''_n)$ and let
$h:\gcl\rw\gcl''$ be a morphism in $\cng$ such that $hp_1=hp_2$. Since
(\ref{app.eq16}) is a coequalizer in groups, there is a map of groups $r:G'\rw
G''$ such that $rf=h$, as maps of groups. We claim that $r$ is in fact a map in
$\cng$. In fact, given $y\in G'$, since $UG\supar{Uf}UG'$ is surjective, there
is $x\in UG$ such that $f(x)=y$. Hence, for each $1\leq i\leq n$,
\begin{equation*}
    rd'_i(y)=rd'_i f(x)=rfd_i(x)=hd_i(x)=d''_ih(x)=d''_irf(x)=d''_ir(y).
\end{equation*}
Hence $rd'_i=d''_ir$. Similarly one proves that $rt'_i=t''_ir$. This shows that
$r$ is a map in $\cng$, as claimed. Hence $f:\gcl\rw\gcl'$ is a coequalizer in
$\cng$, so $R_n$ reflects regular epis.\enpr\bigskip

\emph{\textbf{Proof of Proposition \ref{catn.pro1}}}\medskip

We use induction on $n$. For the case $n=1$, let $f:(G,d,t)\rw(G',d',t')$ be a
morphism in $\ccl^1\gcl$. Then $\ovl{\ncl}f$ is the map of simplicial groups:
\begin{diagram}
    \cdots\tens{G}{\im d}\tiund{\im d}G & \pile{\rTo\\ \rTo\\ \rTo\\ \rTo\\} &
    \tens{G}{\im d} & \pile{\rTo\\ \rTo\\ \rTo\\} & G & \pile{\rTo^d\\ \rTo^t\\
    \lTo_i\\} & \im d\\
    \dTo_{(f,f,f)}&& \dTo_{(f,f)}&& \dTo_f && \dTo_{f_{|}}\\
    \cdots\tens{G'}{\im d'}\tiund{\im d'}G & \pile{\rTo\\ \rTo\\ \rTo\\ \rTo\\} &
    \tens{G'}{\im d'} & \pile{\rTo\\ \rTo\\ \rTo\\} & G' & \pile{\rTo^{d'}\\ \rTo^{t'}\\
    \lTo_i\\} & \im d'.
\end{diagram}
Recall, as in the proof of Lemma \ref{adj.lem1}, that $\ucl_1(G,d,t)=U\ker
d\times U\im d$. Thus, since $\ucl_1 f$ is surjective, both $U\ker d\rw U\ker
d'$ and $U\im d\rw U\im d'$ are surjective. Also, since the map $d:G\rw\im d$
has a section, $G$ is isomorphic to the semidirect product $G\cong\ker d
\rtimes\im d$, and therefore, from above, the map $f:G\rw G'$ is surjective. It
remains to show that for each $k\geq 2\;$ the map $\pro{G}{\im d}{k}\rw$
$\pro{G'}{\im d'}{k}$ is surjective. Consider the case $k=2$. Let
$(x',y')\in\tens{G'}{\im d'}$, so $t'x'=d'y'$. Since $f$ is surjective, there
exists $x,y\in G$ with $f(x)=x'$, $f(y)=y'$. Then
$(x,t(x)d(y^{-1})y)\in\tens{G}{\im d}$ and
\begin{equation*}
\begin{split}
   & (f,f)(x,t(x) d(y^{-1})y)=(f(x),ft(x)fd(y^{-1})f(y))= \\
    & =(x',t'(x')d'(y'^{-1})y')=(x',y').
\end{split}
\end{equation*}
Thus $\tens{G}{\im d}\rw\tens{G'}{\im d'}$ is surjective. The case $k>2$ is
similar.

Suppose, inductively, that the lemma holds for $n-1$. Let
$\gcl=(G,d_1,\ldots,d_n,t_1,$ $\ldots,t_n)$ and
$\hcl=(H,d'_1,\ldots,d'_n,t'_1,\ldots,t'_n)$ be in $\cng$ and $f:\gcl\rw\hcl$ a
morphism in $\cng$ with $\ucl_n f$ surjective. Recall from the proof of Theorem
\ref{adj.the0} that $\ucl_n=\ucl_\nm \vcl_n\beta^{(k)}_n$
 and that $\beta^{(k)}_n\gcl$ is given by
\begin{equation*}
    \tens{G'}{\im d_k}  \rTo  G'  \pile{\rTo^{d_k}\\ \rTo^{t_k}\\
    \lTo_i\\}  \im d_k
\end{equation*}
as in the proof of Lemma \ref{adj.lem1}. Therefore $\ucl_n\gcl=\ucl_\nm\ker
d_k\times\ucl_\nm\im d_k$; hence the hypothesis that $\ucl_n f$ is surjective
implies that both maps $\ucl_\nm\ker d_k\rw\ucl_\nm\ker d'_k$ and $\ucl_\nm\im
d_k\rw\ucl_\nm\im d'_k$ are surjective. On the other hand, as a morphism of
simplicial objects in $[\Delta^{\nm^{op}},\gp]$, $\ovl{\ncl} f$ is given by
\begin{diagram}
    \cdots \tens{\ncl\gcl'}{\ncl\im d_k} & \pile{\rTo\\ \rTo\\ \rTo\\} & \ncl\gcl' & \pile{\rTo\\ \rTo\\
    \lTo\\} & \ncl\im d_k\\
     \dTo_{(\ovl{\ncl}f)_2}&& \dTo_{(\ovl{\ncl}f)_1} && \dTo_{(\ovl{\ncl}f)_0}\\
    \cdots \tens{\ncl\hcl'}{\ncl\im d'_k} & \pile{\rTo\\ \rTo\\ \rTo\\} & \ncl\hcl' & \pile{\rTo\\ \rTo\\
    \lTo\\} & \ncl\im d'_k.
\end{diagram}
Showing that $\ovl{\ncl}f$ is levelwise surjective is equivalent to showing
that for each $i\geq 0$, $(\ovl{\ncl}f)_i$ is levelwise surjective. Since\,
by\, Corollary\, \ref{adj.cor1}, $\;\ucl_\nm$ reflects\, regular\, epis\, and\,
as\, seen\, above, $\;\ucl_\nm\im d_k\rw\ucl_\nm\im d'_k$ is surjective, then
$\im d_k\rw\im d'_k$ is a regular epi in $\cngm$; thus, by the induction
hypothesis, $(\ovl{\ncl}f)_0$ is levelwise surjective.

By Lemma \ref{adj.lem1}, $\ucl_n\gcl\cong UG\cong\ucl_\nm\gcl'$, hence the
hypothesis that $\ucl_n f$ is surjective implies that
$\ucl_\nm\gcl'\rw\ucl_\nm\hcl'$ is surjective. Since $\ucl_\nm$ reflects
regular epis, $\gcl'\rw\hcl'$ is a regular epi in $\cngm$ hence, by the
induction hypothesis, $(\ovl{\ncl}f)_1:\ncl\gcl'\rw\ncl\hcl'$ is levelwise
surjective.

It remains to show that $(\ovl{\ncl}f)_i$ is levelwise surjective for all
$i\geq 2$. Consider the case $i=2$. Let $\,(x',y')\,\in\,\tens{H}{\im d'_k}$\,
so $\,t'_k x'\,=\,d'_ky'$. Since, from above, $UG\rw UH$ is surjective,
$x'=f(x)$, $y'=f(y)$ for $x,y\in G$. Then
$(x,t_k(x)d_k(y^{-1})y)\in\tens{G}{\im d_k}$ and is mapped to $(x',y')$ by
$(f,f)$. This shows that $\;U(\tens{G}{\im d_k})\,\rw\, U(\tens{H}{\im d_k})\,$
is\, surjective. \,But, \,by \,Lemma \;\ref{adj.lem1},
$\;\ucl_\nm\;(\tens{\gcl'}{\im d_k})\;\cong \;U(\tens{G}{\im d_k})$ and
similarly for $\tens{\hcl'}{\im d'_k}$; hence we conclude that
$\ucl_\nm(\tens{\gcl'}{\im d_k})\rw\ucl_\nm(\tens{\hcl'}{\im d'_k})$ is
surjective; since $\ucl_\nm$ reflects regular epis this implies
$\tens{\gcl'}{\im d_k}\rw\tens{\hcl'}{\im d'_k}$ is a regular epi and hence, by
the induction hypothesis, $(\ncl f)_2$ is levelwise surjective. The case $i> 2$
is similar.\enpr


\end{document}